\documentclass[12pt,reqno]{amsart}
\usepackage{amsmath,amssymb,amsthm,mathrsfs}
\usepackage[numbers,sort&compress]{natbib}
\usepackage[a4paper,margin=28mm]{geometry}
\usepackage{xcolor}

\usepackage{hyperref}
\hypersetup{colorlinks=true,linkcolor=blue,citecolor=blue,urlcolor=blue,
	pdftitle={Blow-up suppression by Taylor-Couette flow in an exterior disk},
	pdfauthor={Lili Wang}}
\allowdisplaybreaks[4]
\newtheorem{theorem}{Theorem}[section]
\newtheorem{proposition}{Proposition}[section]
\newtheorem{lemma}{Lemma}[section]
\theoremstyle{remark}
\newtheorem{remark}{Remark}[section]
\numberwithin{equation}{section}

\DeclareMathOperator{\curl}{curl}
\DeclareMathOperator{\diver}{div}

\ExplSyntaxOn
\cs_new_protected:Npn \EditorialReviewNormalizeTitle
{ \tl_replace_all:cnn { @title } { – } { -- } }
\AddToHook{cmd/maketitle/before}{\EditorialReviewNormalizeTitle}
\ExplSyntaxOff
\begin{document}
	\title[Blow-up suppression by Taylor--Couette flow]{Blow-up suppression for the Patlak--Keller--Segel--Navier--Stokes system by Taylor--Couette flow outside a disk}
	\author{Lili~Wang}
	\address[Lili~Wang]{State Key Laboratory of Mathematical Sciences, Academy of Mathematics $\&$ Systems Science, The Chinese Academy of Sciences, Beijing, 100190, China}
	\email{wllmath@163.com}
	\date{\today}
	\begin{abstract}
		Taylor--Couette devices have been investigated for
		enzymatic reactions and cell culture, where cylinder
		rotation provides a means of controlling mixing
		and mass transfer (Schrimpf et al., AIChE J., 2021).
		Motivated by these applications, this paper studies
		blow-up suppression in the two-dimensional
		Patlak--Keller--Segel--Navier--Stokes system outside
		the unit disk. We prove that sufficiently strong Taylor--Couette flow
		yields a unique global strong solution whose density
		and concentration are uniformly bounded in time,
		without smallness assumptions on the initial cell mass
		or the fixed initial velocity perturbation.
		The main difficulty is that the angular shear tends
		to zero at infinity, so decay estimates on bounded
		regions do not directly control the full solution.
		To overcome this difficulty, radially weighted energy
		estimates exploit the radial velocity to control
		the density and vorticity at large radii.
		Combining these bounds with linear decay
		estimates and an energy comparison with the radial
		chemotaxis solution yields convergence of the density
		in $L^2$ to the radial solution with angularly averaged
		initial density, uniformly on fixed finite time
		intervals away from $t=0$ as the rotation amplitude
		tends to infinity.
		The decay of the radial solution then makes the full
		density small in $L^2$ at a finite time.
		An energy estimate preserves this smallness and
		gives global continuation.
	\end{abstract}
	\maketitle
	\tableofcontents
	
	\section{Introduction}\label{sec:model}
	\subsection{Background and related results}
	We consider the two-dimensional
	Patlak--Keller--Segel--Navier--Stokes (PKS--NS) system in the exterior
	of the unit disk $\Omega$:
	\begin{equation}\label{eq:main-system}
		\begin{cases}
			\partial_t n+v\cdot\nabla n
			=\Delta n-\nabla\cdot(n\nabla c),\\
			-\Delta c+c=n,\\
			\partial_t v+v\cdot\nabla v+\nabla P
			=\Delta v+n\nabla r,\\
			\nabla\cdot v=0,
		\end{cases}
	\end{equation}
	where
	$\Omega=\{x\in\mathbb R^2:|x|>1\}$ and $r=|x|$.
	The unknowns $n$, $c$, $v$, and $P$ denote the cell density,
	chemical concentration, fluid velocity, and pressure,
	respectively.
	
	The Patlak--Keller--Segel (PKS) model describes the competition
	between diffusion and chemotactic attraction
	\cite{Patlak1953,KellerSegel1970}.
	For the classical parabolic--elliptic PKS model on
	$\mathbb R^2$ with $-\Delta c=n$, the critical mass $8\pi$
	separates the subcritical and supercritical regimes
	\cite{BlanchetDolbeaultPerthame2006}.
	Since chemotactic processes often occur in moving fluids,
	a natural question is whether fluid transport can prevent
	the concentration responsible for blow-up.
	Although incompressible transport alone preserves
	$L^p$ norms, mixing can create small spatial scales on
	which diffusion acts more efficiently.
	Kiselev and Xu \cite{KiselevXu2016} showed that suitable
	incompressible flows can suppress chemotactic blow-up
	in two and three dimensions.
	For prescribed shear flows, Bedrossian and He
	\cite{BedrossianHe2017} established suppression on the torus
	without a mass restriction in two dimensions.
	Their argument combines enhanced dissipation with
	a reduction of the effective dimension of the density
	dynamics.
	He \cite{He2018} extended the study to the
	parabolic--parabolic system under strictly monotone shear
	on $\mathbb T\times\mathbb R$
	(with concentration equation
	$\partial_t c+V\cdot\nabla c=\Delta c+n-c$,
	where $V$ is the prescribed velocity).
	Iyer, Xu and Zlato\v{s} \cite{IyerXuZlatos2021} obtained
	a suppression criterion in terms of dissipation time.
	Other mechanisms on unbounded domains include fast
	splitting and spreading
	\cite{HeTadmor2019,HeTadmorZlatos2022}.
	Suppression by prescribed Poiseuille flow was studied
	by Li and Xiang \cite{LiXiang2026}.
	
	For the coupled PKS--NS system, the density also drives
	the fluid, so estimates for the velocity must be obtained
	together with control of chemotactic aggregation.
	Zeng, Zhang and Zi \cite{ZengZhangZi2021} proved suppression
	near Couette flow on $\mathbb T\times\mathbb R$ for both
	the parabolic--elliptic and parabolic--parabolic systems.
	Their parabolic--parabolic result requires an additional
	smallness condition on the nonzero modes of the initial
	chemical gradient.
	Subsequent studies treated Poiseuille flow in a finite channel
	and on $\mathbb T\times\mathbb R$
	\cite{CuiWang2024,LiXiangXu2025}, as well as models with
	nonlinear coupling between the cells and the fluid
	\cite{He2023,LiPengXiang2025}.
	Quantitative conditions on the flow amplitude and suppression
	by Couette flow on $\mathbb R^2$ were investigated in
	\cite{ChenWangYang2025,WangWangWu2026}.
	Related results concern Couette stability thresholds
	\cite{Hao2026} and Poiseuille--Couette flow \cite{Ding2026}.
	Supercritical chemotaxis models with volume filling were considered in
	\cite{WangWangZhang2024}, and models with nonlinear
	chemotactic sensitivity in
	\cite{ChenWangYangZhang2026,Wang2026Sensitivity}.
	More recently, Chen, Cui, Wang and Zhang
	\cite{ChenCuiWangZhang2026} removed the additional
	smallness condition on the initial chemical gradient
	for the parabolic--parabolic PKS--NS system near Couette flow.
	
	For related suppression results in three dimensions under
	prescribed flows, see
	\cite{BedrossianHe2017,BedrossianHe2018Erratum,He2025Alternating}.
	Results for the coupled PKS--NS system are given in
	\cite{CuiWangWang2026,CuiWangWangWei2025}.
	Besides imposed shear flows, buoyancy generated by the
	density itself can prevent chemotactic concentration,
	as shown for Stokes--Boussinesq and viscous fluid couplings
	\cite{HuKiselev2024,Hu2024}, and for a Darcy coupling
	\cite{HuKiselevYao2025}.
	These developments lead to the question of suppression
	by circular shear in a domain with an unbounded radial direction.
	
	\subsection{Taylor--Couette flow}
	Taylor--Couette flow is a classical circular shear flow
	associated with the relative rotation of coaxial cylinders
	\cite{Taylor1923}.
	The use of rotating cylinder devices in chemical
	and biological processes is reviewed by
	Schrimpf et al.~\cite{SchrimpfEtAl2021}.
	For applications to extracorporeal blood treatment,
	Ameer et al.~\cite{AmeerEtAl1999} evaluated a
	Taylor--Couette device that combined plasma separation
	with an immobilized enzyme to remove heparin.
	The aim was to maintain anticoagulation within the circuit
	while reducing the heparin returned to the systemic
	circulation.
	K\"orfer, Klaus and Mottaghy
	\cite{KorferKlausMottaghy2003} compared platelet activation
	and protein adsorption under laminar flow and Taylor
	vortex flow.
	These studies illustrate how flow geometry affects
	the transport of biological components and their
	interaction with artificial surfaces.
	The present two-dimensional model examines how circular
	shear counteracts chemotactic concentration when
	the fluid velocity also responds to the density.
	
	For the PKS--NS system in a bounded annulus, Cui, Wang
	and Wang \cite{CuiWangWang2025} proved that sufficiently
	strong Taylor--Couette flow suppresses blow-up without
	requiring the initial cell mass or the fixed initial
	velocity perturbation to be small.
	The exterior problem removes the outer boundary and leaves an unbounded radial direction. After removal of the solid rotation, the angular transport coefficient is $Ar^{-2}$, where $A$ denotes the rotation amplitude. For each fixed $A$, this coefficient tends to zero as $r\to\infty$. Linear enhanced dissipation in this geometry was studied by Li, Zhang and Zhang \cite{LZZ2026}. The nonlinear problem also requires estimates for the density and the evolving velocity perturbation. Theorem~\ref{thm:main} establishes global existence for sufficiently large $A$ and convergence to the radial solution on each fixed finite interval away from $t=0$. Remark~\ref{rem:uniform-positive-time-convergence} strengthens this convergence to all $t\geq\delta$ for every fixed $\delta>0$. The proof combines control at large radii, decay of the nonzero angular modes, and an energy comparison with the radial solution, as explained in Section~\ref{sec:main-proof}. The background velocity and the rotating coordinates are specified below.
	
	In polar coordinates
	$x=(r\cos\theta,r\sin\theta)$, with $r>1$ and
	$\theta\in\mathbb T=\mathbb R/(2\pi\mathbb Z)$, write
	\[
	e_r=(\cos\theta,\sin\theta),\qquad
	e_\theta=(-\sin\theta,\cos\theta).
	\]
	We take the stationary Taylor--Couette flow
	\[
	U_A=A(r+r^{-1})e_\theta,\qquad A\geq1,
	\]
	as the background velocity.
	Its vorticity is $\curl U_A=2A$, and its associated
	radial pressure $P_A$ satisfies
	$P_A'(r)=r^{-1}|U_A|^2$.
	Thus $(0,0,U_A,P_A)$ is a steady solution of
	\eqref{eq:main-system}.
	Define the velocity, pressure and vorticity perturbations by
	\[
	u=v-U_A,\qquad p=P-P_A,\qquad
	w=\curl u=\partial_{x_1}u_2-\partial_{x_2}u_1.
	\]
	Note that the background transport is
	$U_A\cdot\nabla=A(1+r^{-2})\partial_\theta$.
	The term $A\partial_\theta$ describes rigid rotation
	with the same angular velocity at every radius and
	does not contribute to shear mixing.
	We therefore use coordinates rotating with angular
	velocity $A$ to remove this term, leaving
	$Ar^{-2}\partial_\theta$, whose dependence on $r$
	produces angular mixing.
	More precisely, write $u=u_r e_r+u_\theta e_\theta$
	and set
	\begin{equation}\label{eq:rotating coordinates}
		\widetilde f(t,r,\theta)
		=f(t,r,\theta+At)\qquad\text{for}\qquad
		f\in\{n,c,p,w,u_r,u_\theta\}.
	\end{equation}
	For scalar functions, this change replaces
	$\partial_t+A\partial_\theta$ by $\partial_t$.
	For the velocity, it acts on the polar components $\widetilde u_r$ and $\widetilde u_\theta$.
	
	Substituting $v=U_A+u$ into \eqref{eq:main-system}, 
	subtracting the steady equation for $U_A$, and
	applying \eqref{eq:rotating coordinates}, we obtain,
	after omitting tildes,
	\begin{equation}\label{eq:rotated-perturbation-velocity}
		\begin{aligned}
			&	\partial_t u+Ar^{-1}e_\theta\cdot\nabla u
			-\Delta u+\nabla p
			\\=~&-Au\cdot\nabla(r^{-1}e_\theta)
			-2A(-u_2,u_1)-u\cdot\nabla u+ne_r.
		\end{aligned}
	\end{equation}
	The term $2A(-u_2,u_1)$ comes from the rotation of
	the vector coordinates together with
	$u\cdot\nabla(Ar e_\theta)=A(-u_2,u_1)$.
	Since $\nabla\cdot u=0$, we have
	\[
	\curl\bigl(2A(-u_2,u_1)\bigr)
	=2A\nabla\cdot u=0,
	\qquad
	\curl(ne_r)
	=-r^{-1}\partial_\theta n
	=-\nabla\cdot(ne_\theta).
	\]
	Then taking the curl of
	\eqref{eq:rotated-perturbation-velocity} and applying
	the same coordinate change to the density and
	concentration equations in \eqref{eq:main-system} gives
	\begin{equation}\label{eq:physical-scalars-moments}
		\begin{cases}
			\partial_t n+Ar^{-2}\partial_\theta n
			+u\cdot\nabla n
			=\Delta n-\nabla\cdot(n\nabla c),\\
			-\Delta c+c=n,\\
			\partial_t w+Ar^{-2}\partial_\theta w
			+u\cdot\nabla w
			=\Delta w-\nabla\cdot(ne_\theta),\\
			\nabla\cdot u=0.
		\end{cases}
	\end{equation}
	From now on, all unknowns are expressed in these
	rotating coordinates, while the time variable
	$t$ remains unchanged.
	Moreover, we supplement \eqref{eq:physical-scalars-moments}
	with the boundary conditions
	\begin{equation}\label{eq:main-boundary}
		n=c=w=0,\qquad u\cdot e_r=0
		\quad\text{on}\quad r=1,
	\end{equation}
	and the initial conditions
	$(n,u)|_{t=0}=(n_{\mathrm{in}},u_{\mathrm{in}})$.
	Here $u\cdot e_r=0$ is the impermeability condition,
	equivalent to $v\cdot e_r=0$ because $U_A\cdot e_r=0$.
	The condition $w=0$ prescribes zero perturbation
	vorticity at the boundary.
	The angular change of variables preserves these
	boundary conditions and is the identity at $t=0$,
	so the initial data are unchanged.
	
	\subsection{Main results}
	For the local existence argument, the equations for
	$n$, $c$, and $u$ in the rotating coordinates are
	written in operator form as
	\begin{align}
		&\partial_t n+A_Dn=-Ar^{-1}e_\theta\cdot\nabla n-u\cdot\nabla n
		-\nabla\cdot(n\nabla c),
		\qquad(1+A_D)c=n,
		\label{eq:local-rotated-scalar}\\
		&\partial_t u+Su
		=\mathbb P\bigl[-Ar^{-1}e_\theta\cdot\nabla u
		-Au\cdot\nabla(r^{-1}e_\theta)-2A(-u_2,u_1)
		-u\cdot\nabla u+ne_r\bigr].
		\label{eq:local-rotated-velocity}
	\end{align}
	Here $A_D=-\Delta$ is the scalar Dirichlet Laplacian,
	and $\mathbb P$ is the
	$L^2(\Omega;\mathbb R^2)$ Helmholtz projection.
	The operator $S$ is the Stokes operator with boundary
	conditions
	$u\cdot\mathbf n=0$ and $\operatorname{curl}u=0$
	on $r=1$, where $\mathbf n=-e_r$ is the outward unit
	normal. These conditions agree with
	\eqref{eq:main-boundary}.
	
	Proposition~\ref{prop:local-strong-theory} gives local
	existence and uniqueness for each fixed $A\geq1$.
	Its proof is given in Appendix~\ref{app:local},
	where the operator domains and the required linear
	estimates are also established.
	\begin{proposition}[Local existence and uniqueness]
		\label{prop:local-strong-theory}
		Fix $A\geq1$. Let $0\leq n_{\mathrm{in}}\in H_0^1(\Omega)$ and
		$u_{\mathrm{in}}\in H^2(\Omega;\mathbb R^2)$ satisfy
		$\nabla\cdot u_{\mathrm{in}}=0$ in $\Omega$,
		$u_{\mathrm{in}}\cdot\mathbf n=0$, and
		$\operatorname{curl}u_{\mathrm{in}}=0$ on $\partial\Omega$.
		Then the following statements hold.
		\begin{enumerate}
			\item[(i)]
			\textbf{Local existence, uniqueness and regularity.}
			There exists $\tau>0$ such that the system
			\eqref{eq:local-rotated-scalar}--\eqref{eq:local-rotated-velocity}
			has a unique strong solution with initial data
			$(n_{\mathrm{in}},u_{\mathrm{in}})$ and
			\begin{align*}
				\qquad	n&\in C([0,\tau];H_0^1(\Omega))
				\cap L^2(0,\tau;H^2(\Omega)),\qquad\quad\quad
				\partial_tn\in L^2(0,\tau;L^2(\Omega)),\\
				\qquad	u&\in C([0,\tau];H^2(\Omega;\mathbb R^2))
				\cap L^2(0,\tau;H^3(\Omega;\mathbb R^2)),\quad
				\partial_tu\in L^2(0,\tau;H^1(\Omega;\mathbb R^2)).
			\end{align*}
			Moreover, $n\geq0$, the velocity satisfies
			$\nabla\cdot u=0$ in $\Omega$ and
			$u\cdot\mathbf n=\operatorname{curl}u=0$ on $\partial\Omega$,
			and
			\[
			c=(1+A_D)^{-1}n
			\in C([0,\tau];H^3(\Omega)\cap H_0^1(\Omega)).
			\]
			
			\item[(ii)]
			\textbf{Uniform time of existence.}
			For every $R>0$, there exists
			$\tau=\tau(A,\Omega,R)>0$ such that the conclusions
			in {\rm(i)} hold for all initial data satisfying
			the above assumptions and
			\[
			\|n_{\mathrm{in}}\|_{H^1(\Omega)}
			+\|u_{\mathrm{in}}\|_{H^2(\Omega)}\leq R.
			\]
			For data prescribed at any time $t_0\geq0$ and
			satisfying the same assumptions and norm bound,
			the same existence time applies, with the solution
			defined on $[t_0,t_0+\tau]$.
			
			\item[(iii)]
			\textbf{$L^1$ bound.}
			If additionally $n_{\mathrm{in}}\in L^1(\Omega)$, then
			\[
			\|n(t)\|_{L^1(\Omega)}
			\leq\|n_{\mathrm{in}}\|_{L^1(\Omega)}\qquad\text{for}\qquad 
			0\leq t\leq\tau.
			\]
		\end{enumerate}
	\end{proposition}
	For initial data satisfying the additional regularity and
	weighted assumptions below,
	sufficiently strong rotation gives global existence and
	convergence of the density to the radial solution.
	\begin{theorem}\label{thm:main}
		Assume that 
		\begin{equation}\label{eq:initial-regularity}
			0\leq n_{\mathrm{in}}
			\in L^1(\Omega)\cap L^\infty(\Omega)\cap H_0^1(\Omega),
			\qquad
			u_{\mathrm{in}}\in H^2(\Omega;\mathbb R^2),
		\end{equation}
		and
		\begin{equation}\label{eq:initial-moments}
			r^2n_{\mathrm{in}}\in L^2(\Omega),
			\qquad
			r^2w_{\mathrm{in}}\in L^2(\Omega),
			\qquad
			w_{\mathrm{in}}=\curl u_{\mathrm{in}}.
		\end{equation}
		Suppose also that
		$\nabla\cdot u_{\mathrm{in}}=0$ in $\Omega$ and
		$u_{\mathrm{in}}\cdot e_r=w_{\mathrm{in}}=0$ on $r=1$.
		Then the following statements hold.
		
		\begin{enumerate}
			\item[(i)]
			\textbf{Global existence and uniform boundedness.}
			There exist constants $A_0\geq1$ and $M_0>0$,
			depending only on $\Omega$ and the fixed initial data,
			such that, for every $A\geq A_0$, the system
			\eqref{eq:physical-scalars-moments}, with the boundary
			conditions \eqref{eq:main-boundary} and initial data
			$(n,u)|_{t=0}=(n_{\mathrm{in}},u_{\mathrm{in}})$,
			admits a unique global strong solution $(n,c,u,w)$.
			Moreover,
			\[
			\sup_{t\geq0}
			\left(
			\|n(t)\|_{L^2(\Omega)}
			+\|n(t)\|_{L^\infty(\Omega)}
			+\|c(t)\|_{W^{1,\infty}(\Omega)}
			\right)
			\leq M_0.
			\]
			For every $0<T<\infty$, the density and velocity
			perturbation satisfy
			\[
			\begin{aligned}
				n&\in C([0,T];H_0^1(\Omega))
				\cap L^2(0,T;H^2(\Omega)),\\
				u&\in C([0,T];H^2(\Omega;\mathbb R^2))
				\cap L^2(0,T;H^3(\Omega;\mathbb R^2)).
			\end{aligned}
			\]
			\item[(ii)]
			\textbf{Convergence to the radial solution.}
			Let $(N,C)$ denote the unique global radial solution of
			\begin{equation}\label{eq:radial-reference-problem}
				\begin{cases}
					\displaystyle
					\partial_t N
					=\frac1r\partial_r(r\partial_r N)
					-\frac1r\partial_r(rN\partial_r C),
					& t>0,\ r>1,\\[1ex]
					\displaystyle
					-\frac1r\partial_r(r\partial_r C)+C=N,
					& t>0,\ r>1,\\[1ex]
					N(t,1)=C(t,1)=0,
					& t>0,\\[1ex]
					\displaystyle
					N(0,r)=\frac1{2\pi}\int_0^{2\pi}
					n_{\mathrm{in}}(r,\theta) d\theta,
					& r>1,
				\end{cases}
			\end{equation}
			where $C(t)\in H_0^1(\Omega)$.
			Denote by $n^A$ the density obtained in part~{\rm(i)}
			for the parameter $A$. Then, for every
			$0<\delta<T<\infty$, there holds
			\begin{equation}\label{eq:main-averaging}
				\lim_{A\to\infty}
				\sup_{\delta\leq t\leq T}
				\|n^A(t)-N(t)\|_{L^2(\Omega)}
				=0.
			\end{equation}
		\end{enumerate}
	\end{theorem}
	\begin{remark}\label{rem:comparison-earlier-results}
		Cui, Wang and Wang \cite{CuiWangWang2025} established
		blow-up suppression for the PKS--NS system under
		sufficiently strong Taylor--Couette flow in a bounded
		annulus. Theorem~\ref{thm:main} gives the corresponding result
		in the exterior of the unit disk.
		The unbounded radial direction introduces an additional
		difficulty: for each fixed $A$, the angular transport
		coefficient $Ar^{-2}$ in
		\eqref{eq:physical-scalars-moments} and its radial
		derivative tend to zero as $r\to\infty$.
		Consequently, estimates on bounded annuli must
		be supplemented by control of the density and
		vorticity at large radii.
		The weighted assumptions \eqref{eq:initial-moments}
		are used to obtain this control in
		Lemma~\ref{lem:finite-time-uniform-moments}.
		In addition, Theorem~\ref{thm:main} shows that
		the density converges to the radial solution of
		\eqref{eq:radial-reference-problem} in $L^2(\Omega)$
		as $A\to\infty$, uniformly on $[\delta,T]$
		for every $0<\delta<T<\infty$.
	\end{remark}
	\begin{remark}\label{rem:exterior-difficulties}
		The radial velocity plays two roles in the proof.
		First, radial weights reduce the fluid transport terms
		in the weighted energy identities to terms involving
		only $u_r$.
		Under the density bound
		\eqref{eq:finite-time-L2-bootstrap}, the estimate in
		Lemma~\ref{lem:exterior-radial-biot-savart}
		yields a differential inequality that is at most
		linear in the weighted energy.
		The coefficient has a time integral bounded
		independently of $A$ on each fixed finite interval.
		This gives the estimates in
		Lemma~\ref{lem:finite-time-uniform-moments}
		without requiring the fixed initial velocity
		perturbation to be small.
		Second, incompressibility and
		\eqref{eq:main-boundary} imply that $u_r$ has zero
		angular mean.
		Hence the additional term in the averaged density
		equation is a radial derivative of averaged products
		of nonzero angular modes.
		Decay of these modes gives
		\eqref{eq:interaction-flux-vanishes}, which controls
		the contribution of these products to the energy
		estimate for $n_0^A-N$ in
		Proposition~\ref{prop:finite-time-averaging}.
	\end{remark}
	
	\begin{remark}\label{rem:uniform-positive-time-convergence}
		The convergence in \eqref{eq:main-averaging} extends to all
		times bounded away from zero: for every fixed $\delta>0$,
		\begin{equation}\label{eq:main-averaging-uniform-time}
			\lim_{A\to\infty}\sup_{t\geq\delta}
			\|n^A(t)-N(t)\|_{L^2(\Omega)}=0.
		\end{equation}
		Indeed, fix $0<\epsilon<\frac{\rho_*}{2}$, where $\rho_*$ is given by
		Lemma~\ref{lem:small-L2-density}.
		Proposition~\ref{prop:radial-reference-decay} provides
		$T>\delta$ such that
		$\sup_{t\geq T}\|N(t)\|_{L^2(\Omega)}<\epsilon$.
		By \eqref{eq:main-averaging}, for all sufficiently large $A$,
		the error on $[\delta,T]$ is less than $\epsilon$ and
		$\|n^A(T)\|_{L^2(\Omega)}<2\epsilon<\rho_*$.
		Lemma~\ref{lem:small-L2-density} preserves this bound for
		all later times, and hence
		\[
		\sup_{t\geq T}\|n^A(t)-N(t)\|_{L^2(\Omega)}
		\leq\sup_{t\geq T}\|n^A(t)\|_{L^2(\Omega)}
		+\sup_{t\geq T}\|N(t)\|_{L^2(\Omega)}
		\leq3\epsilon.
		\]
		Combining the two time intervals and letting
		$\epsilon\to 0$ proves
		\eqref{eq:main-averaging-uniform-time}.
		The initial data are fixed throughout this limit.
		For a nonradial initial density, $t=0$ cannot be included,
		since $n^A(0)-N(0)=n_{\mathrm{in}}-N(0)$ is independent of $A$
		and does not vanish.
	\end{remark}
	
	\noindent\textbf{Notation.}
	\begin{itemize}
		\item
		Spatial integration is with respect to $dx=rdrd\theta$.
		Unless another domain is specified, spatial norms are taken
		over $\Omega$. In particular, $\|f\|_{L^p}$ and
		$\|f\|_{H^s}$ denote $\|f\|_{L^p(\Omega)}$ and
		$\|f\|_{H^s(\Omega)}$, respectively.
		For a vector field $u=(u_1,u_2)$, the notation
		$u\in H^s(\Omega;\mathbb R^2)$ means that
		$u_1,u_2\in H^s(\Omega)$, with
		\[
		\|u\|_{H^s(\Omega;\mathbb R^2)}^2
		=\|u_1\|_{H^s(\Omega)}^2+\|u_2\|_{H^s(\Omega)}^2.
		\]
		The notation $L^p(\Omega;\mathbb R^2)$ likewise denotes
		the space of vector fields with components in $L^p(\Omega)$.
		We omit $\mathbb R^2$ in vector norms when no ambiguity arises.
		For norms involving both time and space, the time interval
		is also specified, as in
		$\|f\|_{L^q(0,T;L^p(\Omega))}$.
		\item
		At each fixed time, the Fourier coefficients of $f$ in
		the angular variable are defined by
		\[
		\widehat f_k(r)
		=\frac1{2\pi}\int_0^{2\pi}
		f(r,\theta)e^{-ik\theta} d\theta,
		\qquad k\in\mathbb Z.
		\]
		The corresponding Fourier expansion is
		$f(r,\theta)=\sum_{k\in\mathbb Z}\widehat f_k(r)e^{ik\theta}$,
		with convergence in $L^2(\Omega)$ when $f\in L^2(\Omega)$.
		\item
		We denote the zero and nonzero angular modes by
		\[
		f_0(r)=\widehat f_0(r)
		=\frac1{2\pi}\int_0^{2\pi}f(r,\theta) d\theta,
		\qquad
		f_{\neq}(r,\theta)=f(r,\theta)-f_0(r),
		\]
		respectively. For a vector field
		$u=u_r e_r+u_\theta e_\theta$, these projections are applied
		to $u_r$ and $u_\theta$ separately. Thus
		$u_0=(u_r)_0e_r+(u_\theta)_0e_\theta$ and
		$u_{\neq}=u-u_0$.
		\item
		The radial and angular components of $u$ are
		$u_r=u\cdot e_r$ and $u_\theta=u\cdot e_\theta$.
		We use $\nabla^\perp=(-\partial_{x_2},\partial_{x_1})$.
		
		\item
		The superscript $A$ indicates the dependence of a
		solution on the rotation parameter.
		
		\item
		The chemical concentrations in the coupled system and the radial
		problem are denoted by $c$ and $C$, respectively.
		
		\item
		The letter $K$ denotes a positive estimate constant whose value
		may change from line to line. Its dependence on the initial data
		and other parameters is stated in each result or proof.
		Subscripts indicate additional dependence when needed. 
	\end{itemize}

	\section{Main ideas and proof of Theorem~\ref{thm:main}}
	\label{sec:bootstrap}\label{sec:main-proof}
	Proposition~\ref{prop:density-continuation} reduces global
	existence to an $L^2(\Omega)$ bound for the density.
	This section explains how comparison with the radial
	solution provides such a bound.
	The first subsection presents the proof strategy and
	the key estimates. The second proves
	Theorem~\ref{thm:main}.
	\subsection{Proof strategy and key estimates}
	Let $(N,C)$ be the global solution of
	\eqref{eq:radial-reference-problem}, with
	$N(0)=(n_{\mathrm{in}})_0$, given by
	Proposition~\ref{prop:radial-reference-decay}.
		The finite-time comparison must be established before
	global existence is known.
	For this purpose, set
	\begin{equation}\label{eq:bootstrap-level}
		B=4\left(
		1+\|n_{\mathrm{in}}\|_{L^{2}(\Omega)}
		+\sup_{t\geq0}\|N(t)\|_{L^{2}(\Omega)}
		\right).
	\end{equation}
	This constant is finite by
	Proposition~\ref{prop:radial-reference-decay}
	and is independent of $A$.
	For fixed $T>0$, let $[0,\tau_A)$, with
	$0<\tau_A\leq T$, be the maximal initial interval on which
	the strong solution exists and
	$\|n^A(t)\|_{L^2(\Omega)}<B$.
	Under this bound, the constants $K_T$ below may depend
	on $T$, $B$, and the fixed initial data, but are independent
	of $A$, $\tau_A$, and $\delta$.
	The decomposition
	\[
	n^A-N=(n_0^A-N)+n_{\neq}^A
	\]
	requires estimates for both the nonzero angular modes
	and the difference between the zero mode and $N$.
	The discussion below is organized into four parts.
	The first three establish the finite-time comparison,
	and the fourth uses radial decay to obtain global existence.
	
	\noindent\textbf{$\bullet$ Uniform control at large radii.}
	The first difficulty is that, for each fixed $A$,
	the angular transport coefficient $Ar^{-2}$ tends to zero
	as $r\to\infty$.
	Decay estimates on bounded annuli must therefore be
	combined with control at large radii.
	For the nonlinear problem, this control must also remain
	uniform in $A$ despite the evolving velocity perturbation.
	The proof estimates the weighted energy
	\begin{equation}\label{eq:bootstrap-energy}
		\int_\Omega r^4
		\bigl(n^A(t,x)^2+w^A(t,x)^2\bigr)dx.
	\end{equation}
	The background angular transport has zero integral
	against the radial weight.
	For the transport by $u^A$, incompressibility and
	\eqref{eq:main-boundary} give
	\[
	\begin{aligned}
		-\int_\Omega r^4
		\bigl(n^Au^A\cdot\nabla n^A
		+w^Au^A\cdot\nabla w^A\bigr)dx=
		2\int_\Omega r^3u_r^A
		\bigl((n^A)^2+(w^A)^2\bigr)dx.
	\end{aligned}
	\]
	Thus only the radial velocity needs to be estimated.
	Lemma~\ref{lem:exterior-radial-biot-savart},
	interpolation, and the unweighted energy bounds yield
	\[
	\begin{aligned}
		\left|\int_\Omega r^3u_r^A
		\bigl((n^A)^2+(w^A)^2\bigr)dx\right|
		&\leq
		K_T\|\nabla w^A\|_{L^2(\Omega)}^{\frac13}
		\left(
		\|r^2n^A\|_{L^2(\Omega)}^2
		+\|r^2w^A\|_{L^2(\Omega)}^2
		\right)^{\frac{11}{12}}\\
		&\leq
		K_T\bigl(1+\|\nabla w^A\|_{L^2(\Omega)}^2\bigr)
		\left(
		1+\|r^2n^A\|_{L^2(\Omega)}^2
		+\|r^2w^A\|_{L^2(\Omega)}^2
		\right).
	\end{aligned}
	\]
	Since $\frac{11}{12}<1$, the resulting bound is linear
	in one plus the weighted energy.
	Moreover, \eqref{eq:finite-time-g-bound} gives
	\[
	\int_0^{\tau_A}
	\bigl(1+\|\nabla w^A(t)\|_{L^2(\Omega)}^2\bigr)dt
	\leq K_T.
	\]
	After estimating the diffusion, chemotactic, and source
	terms, Gronwall's inequality therefore bounds
	\eqref{eq:bootstrap-energy} independently of $A$.
	No smallness assumption on the fixed initial velocity
	perturbation is required.
	In particular,
	\[
	\sup_{0\leq t<\tau_A}
	\int_{\{r>R\}}
	\bigl((n^A)^2+(w^A)^2\bigr)dx
	\leq K_TR^{-4}.
	\]
	Lemma~\ref{lem:finite-time-uniform-moments}
	proves these estimates and the corresponding spatial
	regularity and bounds at large radii for
	$u^An^A+n^A\nabla c^A$.
	
	\noindent\textbf{$\bullet$ Decay of the nonzero angular modes.}
	The preceding bounds provide the spatial control needed
	to use the linear estimates in the density equation,
	written as
	\[
	\partial_t n^A+Ar^{-2}\partial_\theta n^A-\Delta n^A
	=-\nabla\cdot(u^An^A+n^A\nabla c^A).
	\]
	The linear estimates of Li, Zhang and Zhang
	\cite{LZZ2026} are used for both the homogeneous and
	the forced equations.
	Their use here requires an approximation argument:
	the bounds above control the vector field
	$u^An^A+n^A\nabla c^A$, but do not directly control
	its divergence in the weighted norm required by the
	linear estimate.
	
	The spatial regularity on bounded annuli and the bounds
	at large radii allow this vector field to be approximated,
	uniformly in $A$, by vector fields that are smooth and
	compactly supported in space.
	For each fixed approximation, the linear estimates give
	decay of the nonzero modes of the corresponding solution
	with zero initial data.
	The error between this solution and the solution with
	the original forcing is bounded by
	the $L^2(0,\tau_A;L^2(\Omega))$ approximation error,
	using an energy estimate independent of $A$.
	Combining this argument with the homogeneous linear
	decay and the weighted bound above gives
	\[
	\sup_{\delta\leq t<\tau_A}
	\|n_{\neq}^A(t)\|_{L^2(\Omega)}
	\to 0
	\qquad\text{as}\qquad A\to\infty
	\]
	for every fixed $0<\delta<T$.
	The supremum is taken to be zero if $\tau_A\leq\delta$.
	This passage from the linear estimates to an equation
	with forcing in divergence form is proved in
	Lemma~\ref{lem:finite-time-fast-rotation}.
	
	\noindent\textbf{$\bullet$ Comparison of the density zero mode with $N$.}
	To complete the comparison, it remains to estimate
	$n_0^A-N$.
	Angular averaging removes the rotation term, so the
	linear decay estimates do not apply directly to this
	difference.
	The proof therefore uses an energy estimate for the
	difference between the two density equations.
	
	Incompressibility and \eqref{eq:main-boundary} imply
	\[
	\partial_r(r(u_r^A)_0)=0,\qquad
	(u_r^A)_0(1)=0,
	\]
	and hence $(u_r^A)_0=0$.
	Consequently, the equation for $n_0^A$ has the same
	radial diffusion and chemotactic terms as the equation
	for $N$, together with
	\[
	-\frac1r\partial_r\left[
	r\bigl(
	n_{\neq}^A\partial_r c_{\neq}^A
	+u_r^A n_{\neq}^A
	\bigr)_0
	\right].
	\]
	The significance of this structure is that both averaged
	products contain the decaying factor $n_{\neq}^A$.
	The uniform estimates give
	\[
	\begin{aligned}
		\left(\int_0^{\tau_A}\int_1^\infty r
		\left|
		(n_{\neq}^A\partial_r c_{\neq}^A
		+u_r^A n_{\neq}^A)_0
		\right|^2 drdt\right)^{\frac12}\leq
		K_T\sqrt{\delta}
		+K_T\sup_{\delta\leq t<\tau_A}
		\|n_{\neq}^A(t)\|_{L^2(\Omega)}.
	\end{aligned}
	\]
	The first term bounds the contribution from the initial
	interval of length at most $\delta$, where decay of the
	nonzero modes is not available.
	For fixed $\delta$, the second tends to zero by the
	convergence established in the preceding part.
	Taking $A\to\infty$ and then $\delta\to 0$
	proves \eqref{eq:interaction-flux-vanishes}.
	Since $n_0^A$ and $N$ have the same initial data,
	an energy estimate for $n_0^A-N$, the radial elliptic
	estimate for $c_0^A-C$, and Gronwall's inequality give
	\eqref{eq:mean-convergence-varying}.
	Thus the zero mode converges to $N$ uniformly on
	$[0,\tau_A)$, while the nonzero modes converge to zero
	uniformly away from $t=0$.
	
	Together with the density energy estimate near $t=0$,
	these limits improve the assumed density bound to
	\[
	\sup_{0\leq t<\tau_A}
	\|n^A(t)\|_{L^2(\Omega)}\leq\frac B2
	\]
	for sufficiently large $A$.
	Continuity and
	Proposition~\ref{prop:density-continuation}
	then give $\tau_A=T$ and existence through time $T$.
	These estimates establish the finite-time comparison
	stated in Proposition~\ref{prop:finite-time-averaging}.
	The detailed argument is given in its proof.
	
	\noindent\textbf{$\bullet$ From radial decay to global existence.}
	The finite-time comparison can now be applied at a time
	when the radial density is small in $L^2(\Omega)$.
	The required decay is proved in
	Proposition~\ref{prop:radial-reference-decay}.
	In the radial problem, the zero density condition at
	$r=1$ allows cell mass to leave the domain, and the total
	mass lost cannot exceed the initial mass.
	Using local regularity and the positivity of the radial
	derivative at the inner boundary for nonzero solutions,
	the proof shows that failure of decay on a bounded annulus
	would force a fixed positive loss of mass on infinitely
	many disjoint time intervals.
	This contradicts the bound on the total mass loss.
	Combining decay on bounded annuli with the uniform
	estimates at large radii then gives
	$\|N(t)\|_{L^2(\Omega)}\to0$.
	
	Let $\rho_*>0$ be the constant in
	Lemma~\ref{lem:small-L2-density}.
	Choose $T_0>0$ such that
	$\|N(T_0)\|_{L^2(\Omega)}<\frac{\rho_*}{4}$.
	The finite-time comparison gives, for all sufficiently
	large $A$,
	\[
	\|n^A(T_0)\|_{L^2(\Omega)}
	\leq
	\|N(T_0)\|_{L^2(\Omega)}
	+\|n^A(T_0)-N(T_0)\|_{L^2(\Omega)}
	<\frac{\rho_*}{2}.
	\]
	Lemma~\ref{lem:small-L2-density} preserves this smallness
	for subsequent times.
	The density is therefore bounded in $L^2(\Omega)$
	throughout its maximal interval of existence, and
	Proposition~\ref{prop:density-continuation}
	gives global existence.
	The proof below makes the choice of the rotation
	threshold precise and derives the uniform bounds for
	the density and concentration.
	
	\subsection{Proof of Theorem~\ref{thm:main}}
	We now apply the finite-time comparison and the
	smallness estimate for the density to prove
	Theorem~\ref{thm:main}.
	\begin{proof}
		Throughout this proof, $n^A$ and $c^A$ denote the
		density $n$ and concentration $c$ corresponding
		to the rotation amplitude $A$.
		
		\noindent\textbf{(i) Global existence and uniform boundedness.}
		Let $\rho_*>0$ be the constant in
		Lemma \ref{lem:small-L2-density}.
		By Proposition \ref{prop:radial-reference-decay}, choose
		$T_0>0$ such that
		$\|N(T_0)\|_{L^{2}}<\frac{\rho_*}{4}$.
		Apply Proposition \ref{prop:finite-time-averaging} with
		$T=T_0$ and $\delta=\frac{T_0}{2}$.
		The two convergence statements in this proposition give
		$A_0\geq A_{T_0}$ such that, for every $A\geq A_0$,
		the solution exists on $[0,T_0]$ and
		\[
		\|n^A(T_0)-N(T_0)\|_{L^{2}}<\frac{\rho_*}{4}.
		\]
		Consequently,
		\begin{equation}\label{eq:main-proof-small-L2-at-T0}
			\|n^A(T_0)\|_{L^{2}}
			\leq\|N(T_0)\|_{L^{2}}
			+\|n^A(T_0)-N(T_0)\|_{L^{2}}
			<\frac{\rho_*}{2}.
		\end{equation}
		Applying Lemma~\ref{lem:small-L2-density} at time $T_0$ and using \eqref{eq:main-proof-small-L2-at-T0}, we obtain
		\[
		\|n^A(t)\|_{L^2(\Omega)}\leq\|n^A(T_0)\|_{L^2(\Omega)}<\frac{\rho_*}{2}
		\]
		for every $t\geq T_0$ belonging to the maximal interval of existence.
		On $[0,T_0]$, Proposition \ref{prop:finite-time-averaging}
		gives
		$\sup_{0\leq t\leq T_0}\|n^A(t)\|_{L^{2}}\leq\frac{B}{2}$.
		These two estimates bound the density in $L^2$ throughout
		its maximal interval of existence.
		Proposition \ref{prop:density-continuation} therefore
		excludes a finite maximal existence time.
		In particular,
		\begin{equation}\label{eq:main-proof-global-L2}
			\sup_{t\geq0}\|n^A(t)\|_{L^{2}}
			\leq\max\left\{\frac{B}{2},\frac{\rho_*}{2}\right\}
			\leq\frac{B+\rho_*}{2}.
		\end{equation}
		By Lemma \ref{lem:exterior-moser} and
		\eqref{eq:main-proof-global-L2}, there holds
		\begin{equation}\label{eq:main-proof-global-Linfty}
			\sup_{t\geq0}\|n^A(t)\|_{L^{\infty}}
			\leq K
			\bigl(\|n_{\mathrm{in}}\|_{L^{\infty}}+B+\rho_*\bigr)
			\bigl(1+(B+\rho_*)^2\bigr),
		\end{equation}
		where $K$ depends only on $\Omega$.
		Lemma \ref{lem:screened-regularity}, Sobolev embedding,
		and interpolation give
		\begin{equation}\label{eq:main-proof-chemical-bound}
			\|c^A(t)\|_{W^{1,\infty}}
			\leq K\|n^A(t)\|_{L^{4}}
			\leq K\|n^A(t)\|_{L^{2}}^{\frac12}
			\|n^A(t)\|_{L^{\infty}}^{\frac12}.
		\end{equation}
		Combining \eqref{eq:main-proof-global-L2},
		\eqref{eq:main-proof-global-Linfty}, and
		\eqref{eq:main-proof-chemical-bound}, we obtain the
		uniform bound in part~{\rm(i)}.
		More precisely, we may take
		\[
		M_0
		=K\bigl(1+\|n_{\mathrm{in}}\|_{L^{\infty}}\bigr)
		\bigl(1+B+\rho_*\bigr)^3,
		\]
		with a sufficiently large fixed $K$ depending only on
		$\Omega$. Since $B$ and $\rho_*$ are independent of $A$,
		the constant $M_0$ is independent of $A\geq A_0$
		and $t\geq0$.
		
		For each fixed $A$, the local theory and continuation
		argument in Appendix~\ref{app:local} give the stated
		strong regularity on every finite time interval.
		Uniqueness follows from
		Proposition \ref{prop:local-strong-theory}.
		
		\noindent\textbf{(ii) Convergence to the radial solution.}
		Fix $0<\delta<T<\infty$ and take
		$A\geq\max\{A_0,A_T\}$.
		By uniqueness, the global solution constructed in
		part~{\rm(i)} coincides on $[0,T]$ with the solution
		provided by Proposition~\ref{prop:finite-time-averaging}.
		The convergence estimates in that proposition therefore
		apply to the global solution.
		Using $n^A-N=(n_0^A-N)+n_{\neq}^A$, we obtain
		\begin{align*}
			\sup_{\delta\leq t\leq T}
			\|n^A(t)-N(t)\|_{L^{2}}\leq\sqrt{2\pi}
			\sup_{0\leq t\leq T}
			\left(
			\int_1^\infty
			r|n_0^A(t,r)-N(t,r)|^2 dr
			\right)^{\frac12}
			+\sup_{\delta\leq t\leq T}
			\|n_{\neq}^A(t)\|_{L^{2}}.
		\end{align*}
		The first term tends to zero by
		\eqref{eq:finite-time-averaging}, and the second by
		\eqref{n omega nonzero}.
		This proves \eqref{eq:main-averaging}.
	\end{proof}
	\section{Elliptic and product estimates}
	\label{app:scalar}
	The concentration equation requires elliptic estimates on
	the exterior domain, including a radial gradient bound.
	A Sobolev product estimate controls the nonlinear terms
	in the proof of Proposition~\ref{prop:finite-time-averaging}.
	\begin{lemma}\label{lem:screened-regularity}
		Let $c\in H_0^1(\Omega)$ solve $-\Delta c+c=f$.
		If $2\leq p<\infty$ and $f\in L^2(\Omega)\cap L^p(\Omega)$, then
		\begin{equation}\label{eq:screened-W2p}
			\|c\|_{W^{2,p}(\Omega)}\leq K_p\|f\|_{L^p(\Omega)}.
		\end{equation}
		If $f\in H^1(\Omega)$, then $\|c\|_{H^3(\Omega)}\leq K\|f\|_{H^1(\Omega)}$. 	Here $K_p$ depends on $p$ and $\Omega$, and $K$ depends only on $\Omega$.
	\end{lemma}
	\begin{proof}
		Since $c\in H_0^1(\Omega)$, its zero extension belongs to
		$H^1(\mathbb R^2)$. The two-dimensional Sobolev embedding
		therefore gives
		$c\in L^p(\Omega)$ for every $2\leq p<\infty$.
		We first estimate $\|c\|_{L^p(\Omega)}$.
		For $M>0$, the function
		$c\min\{|c|,M\}^{p-2}$ belongs to $H_0^1(\Omega)$,
		where it is understood as $c$ when $p=2$.
		Using this function to test $-\Delta c+c=f$, discarding
		a nonnegative gradient term on $\{|c|\geq M\}$, and
		applying H\"older's inequality, we obtain
		\begin{align*}
			&(p-1)\int_{\{|c|<M\}}
			|c|^{p-2}|\nabla c|^2dx
			+\int_\Omega |c|^2\min\{|c|,M\}^{p-2}dx\\
			\leq~&
			\int_\Omega fc\min\{|c|,M\}^{p-2}dx\leq
			\int_\Omega |f||c|^{p-1}dx
			\leq\|f\|_{L^p(\Omega)}
			\|c\|_{L^p(\Omega)}^{p-1}.
		\end{align*}
		Letting $M\to\infty$ by monotone convergence gives
		\[
		(p-1)\int_\Omega |c|^{p-2}|\nabla c|^2 dx
		+\|c\|_{L^p(\Omega)}^p
		\leq
		\|f\|_{L^p(\Omega)}
		\|c\|_{L^p(\Omega)}^{p-1}.
		\]
		Consequently,
		\begin{equation}\label{eq:screened-Lp-bound}
			\|c\|_{L^p(\Omega)}
			\leq \|f\|_{L^p(\Omega)}.
		\end{equation}
		To estimate the spatial derivatives, apply the local
		Dirichlet estimates in \cite[Chapter~9]{GT2001} near
		the unit circle. Since $c=0$ on $r=1$, these estimates give
		\begin{equation}\label{eq:screened-boundary-W2p}
			\|c\|_{W^{2,p}(\{1<r<2\})}
			\leq K_p\Bigl(
			\|f\|_{L^p(\{1<r<3\})}
			+\|c\|_{L^p(\{1<r<3\})}
			\Bigr).
		\end{equation}
		For the region $\{r\geq2\}$, choose balls
		$B_{\frac14}(x_j)$ covering this region such that
		$B_{\frac12}(x_j)\subset\Omega$ and the larger balls have
		uniformly bounded overlap. The interior estimates from
		the same chapter in \cite{GT2001} yield
		\begin{equation}\label{eq:screened-interior-W2p}
			\|c\|_{W^{2,p}(B_{\frac14}(x_j))}^p
			\leq K_p\Bigl(
			\|f\|_{L^p(B_{\frac12}(x_j))}^p
			+\|c\|_{L^p(B_{\frac12}(x_j))}^p
			\Bigr),
		\end{equation}
		where $K_p$ is independent of $j$.
		Summing \eqref{eq:screened-interior-W2p} and using
		the bounded overlap, we obtain
		\[
		\begin{aligned}
			\|c\|_{W^{2,p}(\{r>2\})}^p
			&\leq
			\sum_j\|c\|_{W^{2,p}(B_{\frac14}(x_j))}^p\leq K_p\sum_j\Bigl(
			\|f\|_{L^p(B_{\frac12}(x_j))}^p
			+\|c\|_{L^p(B_{\frac12}(x_j))}^p
			\Bigr)\\
			&\leq K_p\Bigl(
			\|f\|_{L^p(\Omega)}^p
			+\|c\|_{L^p(\Omega)}^p
			\Bigr).
		\end{aligned}
		\]
		Combining this inequality with
		\eqref{eq:screened-boundary-W2p} and
		\eqref{eq:screened-Lp-bound} proves
		\eqref{eq:screened-W2p}. In particular,
		\begin{equation}\label{eq:screened-H2-bound}
			\|c\|_{H^2(\Omega)}
			\leq K\|f\|_{L^2(\Omega)}.
		\end{equation}
		Assume now that $f\in H^1(\Omega)$.
		Away from the boundary, differentiation in the sense
		of distributions gives
		\[
		(1-\Delta)\partial_{x_i}c=\partial_{x_i}f,
		\qquad i\in\{1,2\}.
		\]
		Applying the interior $H^2$ estimate to these equations
		on the preceding pairs of balls and summing, we obtain
		\begin{equation}\label{eq:screened-interior-H3}
			\begin{aligned}
				\|c\|_{H^3(\{r>2\})}^2
				&\leq K\sum_j\Bigl(
				\|f\|_{H^1(B_{\frac12}(x_j))}^2
				+\|c\|_{H^2(B_{\frac12}(x_j))}^2
				\Bigr)\\&\leq K\Bigl(
				\|f\|_{H^1(\Omega)}^2
				+\|c\|_{H^2(\Omega)}^2
				\Bigr).
			\end{aligned}
		\end{equation}
		It remains to estimate the third derivatives near $r=1$.
		Choose a smooth radial cutoff $\chi$ satisfying
		$\chi=1$ for $1\leq r\leq2$ and $\chi=0$ for $r\geq3$.
		Angular translations preserve $H_0^1(\Omega)$, so
		\[
		\chi(r)\frac{c(r,\theta+h)-c(r,\theta)}{h}
		\in H_0^1(\Omega),\qquad h\neq0.
		\]
		Since $c\in H^2(\Omega)$ and $\chi$ has bounded support,
		these functions converge to $\chi\partial_\theta c$
		in $H^1(\Omega)$ as $h\to0$.
		The closedness of $H_0^1(\Omega)$ therefore gives
		$\chi\partial_\theta c\in H_0^1(\Omega)$.
		Differentiating $-\Delta c+c=f$ with respect to $\theta$
		and using
		$\Delta(\partial_\theta c)=\partial_\theta(\Delta c)$,
		we obtain
		$(1-\Delta)\partial_\theta c=\partial_\theta f$.
		The product rule then gives
		\[
		(1-\Delta)(\chi\partial_\theta c)
		=\chi\partial_\theta f
		-2\nabla\chi\cdot\nabla\partial_\theta c
		-(\Delta\chi)\partial_\theta c.
		\]
		The cutoff is supported in a fixed bounded region.
		Thus \eqref{eq:screened-H2-bound}, applied to this equation,
		yields
		\begin{equation}\label{eq:screened-angular-H2}
			\begin{aligned}
				\|\chi\partial_\theta c\|_{H^2(\Omega)}
				&\leq K\Bigl(
				\|\chi\partial_\theta f\|_{L^2(\Omega)}
				+\|\nabla\chi\cdot\nabla\partial_\theta c\|_{L^2(\Omega)}
				+\|(\Delta\chi)\partial_\theta c\|_{L^2(\Omega)}
				\Bigr)\\
				&\leq K\Bigl(
				\|f\|_{H^1(\Omega)}
				+\|c\|_{H^2(\Omega)}
				\Bigr).
			\end{aligned}
		\end{equation}
		In particular, all third polar derivatives containing
		an angular derivative satisfy
		\[
		\begin{aligned}
			\|\partial_r^2\partial_\theta c\|_{L^2(\{1<r<2\})}
			+\|\partial_r\partial_\theta^2c\|_{L^2(\{1<r<2\})}
			+\|\partial_\theta^3c\|_{L^2(\{1<r<2\})}\leq
			K\|\chi\partial_\theta c\|_{H^2(\Omega)}.
		\end{aligned}
		\]
		For the remaining derivative, writing the elliptic equation
		in polar coordinates and differentiating with respect to $r$, we have
		\[
		\partial_r^3c
		=\partial_rc-\partial_rf
		+r^{-2}\partial_rc-r^{-1}\partial_r^2c
		+2r^{-3}\partial_\theta^2c
		-r^{-2}\partial_r\partial_\theta^2c.
		\]
		Since $1<r<2$, this identity and
		\eqref{eq:screened-angular-H2} imply
		\[
		\begin{aligned}
			\|\partial_r^3c\|_{L^2(\{1<r<2\})}
			\leq K\Bigl(
			\|f\|_{H^1(\Omega)}
			+\|c\|_{H^2(\Omega)}
			+\|\chi\partial_\theta c\|_{H^2(\Omega)}
			\Bigr)\leq K\Bigl(
			\|f\|_{H^1(\Omega)}
			+\|c\|_{H^2(\Omega)}
			\Bigr).
		\end{aligned}
		\]
		Combining the preceding bounds with
		\eqref{eq:screened-interior-H3} and then using
		\eqref{eq:screened-H2-bound}, we conclude that
		\[
		\|c\|_{H^3(\Omega)}
		\leq K\Bigl(
		\|f\|_{H^1(\Omega)}+\|c\|_{H^2(\Omega)}
		\Bigr)
		\leq K\|f\|_{H^1(\Omega)}.
		\]
	\end{proof}
	The radial estimate below controls $C_r$ in
	\eqref{eq:radial-reference-uniform-bounds} and
	$\partial_r(c_0^A-C)$ in \eqref{eq:mean-chemical-difference}.
	\begin{lemma}[Radial elliptic estimate]
		\label{lem:exterior-radial-screened}
		Assume that $\int_1^\infty r|f(r)|^2 dr<\infty$. The problem
		\begin{equation}\label{eq:radial-elliptic-problem}
			-c''-r^{-1}c'+c=f,\qquad c(1)=0,
			\qquad \int_1^\infty r\bigl(|c|^2+|c'|^2\bigr)dr<\infty
		\end{equation}
		has a unique weak solution. It satisfies
		\begin{equation}\label{eq:radial-elliptic-estimates}
			\max\left\{\int_1^\infty r|c|^2 dr,
			\int_1^\infty r|c'|^2 dr\right\}
			\leq\int_1^\infty r|f|^2 dr,\quad
			\|c'\|_{L^\infty(1,\infty)}^2
			\leq 6\int_1^\infty r|f|^2 dr.
		\end{equation}
	\end{lemma}
	\begin{proof}
		We seek a radial function $c\in H_0^1(\Omega)$ such that
		\begin{equation}\label{eq:radial-elliptic-weak-form}
			\int_1^\infty r(c'\varphi'+c\varphi)dr
			=\int_1^\infty rf\varphi dr
		\end{equation}
		for every radial test function $\varphi\in H_0^1(\Omega)$.
		The radial subspace of $H_0^1(\Omega)$ is closed:
		rotations are isometries in $H^1(\Omega)$, so a limit
		in $H^1(\Omega)$ of radial functions remains radial.
		Equipped with the $H^1(\Omega)$ norm, this subspace
		is therefore a Hilbert space.
		For radial functions $c,\varphi\in H_0^1(\Omega)$,
		the Cauchy--Schwarz inequality gives
		\[
		\begin{aligned}
			\left|\int_1^\infty r(c'\varphi'+c\varphi)dr\right|
			&\leq
			\left(\int_1^\infty r(|c'|^2+|c|^2)dr\right)^{\frac12}
			\left(\int_1^\infty
			r(|\varphi'|^2+|\varphi|^2)dr\right)^{\frac12}\\
			&=\frac1{2\pi}
			\|c\|_{H^1(\Omega)}\|\varphi\|_{H^1(\Omega)}.
		\end{aligned}
		\]
		Moreover,
		\[
		\int_1^\infty r(|c'|^2+|c|^2)dr
		=\frac1{2\pi}\|c\|_{H^1(\Omega)}^2.
		\]
		Thus the bilinear form on the left-hand side of
		\eqref{eq:radial-elliptic-weak-form} is bounded and coercive
		on the radial subspace.
		For the right-hand side, the Cauchy--Schwarz inequality yields
		\[
		\begin{aligned}
			\left|\int_1^\infty rf\varphi dr\right|
			\leq
			\left(\int_1^\infty r|f|^2 dr\right)^{\frac12}
			\left(\int_1^\infty r|\varphi|^2 dr\right)^{\frac12}
			\leq\frac1{\sqrt{2\pi}}
			\left(\int_1^\infty r|f|^2 dr\right)^{\frac12}
			\|\varphi\|_{H^1(\Omega)}.
		\end{aligned}
		\]
		The assumption on $f$ therefore makes this expression
		a bounded linear functional.
		The Lax--Milgram theorem gives a unique radial
		$c\in H_0^1(\Omega)$ satisfying
		\eqref{eq:radial-elliptic-weak-form}. 
		Taking $\varphi\in C_c^\infty(1,\infty)$ in
		\eqref{eq:radial-elliptic-weak-form}, we obtain
		\[
		-(rc')'+rc=rf
		\qquad\text{in the sense of distributions on }(1,\infty).
		\]
		Noting that $c\in H_0^1(\Omega)$ is radial, its trace satisfies
		$c(1)=0$, and
		\[
		\int_1^\infty r(|c'|^2+|c|^2)dr<\infty.
		\]
		Hence $c$ solves \eqref{eq:radial-elliptic-problem}.
		Conversely, any radial solution in $H_0^1(\Omega)$
		satisfies \eqref{eq:radial-elliptic-weak-form},
		by the density of smooth, compactly supported radial functions
		in the radial subspace of $H_0^1(\Omega)$.
		The uniqueness given by the Lax--Milgram theorem therefore
		also proves uniqueness for
		\eqref{eq:radial-elliptic-problem}.
		
		Since $c$ belongs to the test-function space, we may choose
		$\varphi=c$ in \eqref{eq:radial-elliptic-weak-form}.
		Young's inequality then gives
		\[
		\begin{aligned}
			\int_1^\infty r(|c'|^2+|c|^2)dr
			=\int_1^\infty rfc dr\leq\frac12\int_1^\infty r|f|^2 dr
			+\frac12\int_1^\infty r|c|^2 dr.
		\end{aligned}
		\]
		Hence
		\begin{equation}\label{eq:radial-elliptic-energy-bound}
			\int_1^\infty r|c'|^2 dr
			+\frac12\int_1^\infty r|c|^2 dr
			\leq\frac12\int_1^\infty r|f|^2 dr.
		\end{equation}
		In particular, both
		$\int_1^\infty r|c|^2 dr$ and
		$\int_1^\infty r|c'|^2 dr$
		are bounded by $\int_1^\infty r|f|^2 dr$.
		This proves the first estimate in
		\eqref{eq:radial-elliptic-estimates}.
		
		We next estimate the derivative of $c'$.
		Noting that $r\geq1$, the assumption on $f$ and
		\eqref{eq:radial-elliptic-energy-bound} imply that
		$f,c,c'\in L^2(1,\infty)$.
		Equation \eqref{eq:radial-elliptic-problem} gives
		$c''=c-f-r^{-1}c'$
		in the sense of distributions. Its right-hand side belongs
		to $L^2(1,\infty)$, and
		\begin{equation}\label{eq:radial-elliptic-second-derivative}
			\begin{aligned}
				\int_1^\infty |c''|^2 dr
				\leq3\int_1^\infty
				\bigl(|c|^2+|f|^2+r^{-2}|c'|^2\bigr)dr\leq9\int_1^\infty r|f|^2 dr.
			\end{aligned}
		\end{equation}
		Thus $c'\in H^1(1,\infty)$.
		The finiteness of $\int_1^\infty|c'|^2 dr$ allows us to
		choose $R_j\to\infty$ such that $c'(R_j)\to0$.
		For each fixed $r\geq1$ and $R_j>r$, the fundamental theorem
		of calculus gives
		\[
		|c'(r)|^2
		=|c'(R_j)|^2
		-2\int_r^{R_j}c'(\rho)c''(\rho) d\rho.
		\]
		By the Cauchy--Schwarz inequality, $c'c''$ belongs to
		$L^1(1,\infty)$.
		We may therefore pass to the limit as $j\to\infty$ and obtain
		\[
		\begin{aligned}
			|c'(r)|^2
			=-2\int_r^\infty c'(\rho)c''(\rho)d\rho\leq
			2\left(\int_1^\infty|c'|^2 dr\right)^{\frac12}
			\left(\int_1^\infty|c''|^2 dr\right)^{\frac12}.
		\end{aligned}
		\]
		Using $r\geq1$,
		\eqref{eq:radial-elliptic-energy-bound} and
		\eqref{eq:radial-elliptic-second-derivative}, we conclude that
		\[
		|c'(r)|^2\leq6\int_1^\infty r|f|^2 dr.
		\]
		Taking the supremum over $r\geq1$ proves the second estimate
		in \eqref{eq:radial-elliptic-estimates}.
	\end{proof}
	The nonlinear terms $un$ and $n\nabla c$ also require an
	estimate for products of $H^1(\Omega)$ functions.
	\begin{lemma}\label{lem:exterior-product}
		For every $0<s<1$ and $f,g\in H^1(\Omega)$, there holds
		\[
		\|fg\|_{H^s(\Omega)}
		\leq K_s\|f\|_{H^1(\Omega)}\|g\|_{H^1(\Omega)},
		\]
		where $K_s$ depends only on $s$ and $\Omega$.
	\end{lemma}
	\begin{proof}
		Reflection across the circle near $r=1$, followed by a
		smooth cutoff inside the disk, gives a bounded extension
		operator
		$E:H^1(\Omega)\to H^1(\mathbb R^2)$ satisfying
		\begin{equation}\label{eq:exterior-product-extension}
			(Ef)|_\Omega=f,\qquad
			\|Ef\|_{H^1(\mathbb R^2)}
			\leq K\|f\|_{H^1(\Omega)}.
		\end{equation}
		It therefore suffices to establish the corresponding
		product estimate on $\mathbb R^2$.
		First let $f,g\in C_c^\infty(\mathbb R^2)$.
		In this proof, $\widehat f(\xi)$ denotes the Fourier
		transform of $f$ with respect to $x\in\mathbb R^2$.
		For each fixed $\xi\in\mathbb R^2$, the Fourier convolution
		formula gives
		\[
		|\widehat{fg}(\xi)|^2
		\leq K\left|
		\int_{\mathbb R^2}
		\widehat f(\eta)\widehat g(\xi-\eta)d\eta
		\right|^2.
		\]
		Writing
		\[
		\begin{aligned}
			\widehat f(\eta)\widehat g(\xi-\eta)
			=
			\frac{1}
			{\sqrt{1+|\eta|^2}\sqrt{1+|\xi-\eta|^2}}
			\left[
			\sqrt{1+|\eta|^2}\,\widehat f(\eta)
			\sqrt{1+|\xi-\eta|^2}\,\widehat g(\xi-\eta)
			\right]
		\end{aligned}
		\]
		and applying the Cauchy--Schwarz inequality in $\eta$,
		we obtain
		\begin{equation}\label{eq:exterior-product-fourier}
			\begin{aligned}
				|\widehat{fg}(\xi)|^2
				&\leq K
				\left(
				\int_{\mathbb R^2}
				\frac{d\eta}
				{(1+|\eta|^2)(1+|\xi-\eta|^2)}
				\right)\\
				&\quad\times
				\left(
				\int_{\mathbb R^2}
				(1+|\eta|^2)|\widehat f(\eta)|^2
				(1+|\xi-\eta|^2)
				|\widehat g(\xi-\eta)|^2 d\eta
				\right).
			\end{aligned}
		\end{equation}
		To estimate the first integral in
		\eqref{eq:exterior-product-fourier}, consider first
		$|\xi|\leq2$. By the Cauchy--Schwarz inequality, there holds
		\[
		\begin{aligned}
			\int_{\mathbb R^2}
			\frac{d\eta}
			{(1+|\eta|^2)(1+|\xi-\eta|^2)}\leq
			\left(
			\int_{\mathbb R^2}\frac{d\eta}{(1+|\eta|^2)^2}
			\right)^{\frac12}
			\left(
			\int_{\mathbb R^2}
			\frac{d\eta}{(1+|\xi-\eta|^2)^2}
			\right)^{\frac12}
			=\pi .
		\end{aligned}
		\]
		For $|\xi|\geq2$, set $R=|\xi|$.
		On $\{|\eta|\leq \frac{R}{2}\}$, one has
		$|\xi-\eta|\geq \frac{R}{2}$, and hence
		\[
		\begin{aligned}
			&\int_{\{|\eta|\leq \frac{R}{2}\}}
			\frac{d\eta}
			{(1+|\eta|^2)(1+|\xi-\eta|^2)}\\
			\leq~&
			\frac4{R^2}
			\int_{\{|\eta|\leq \frac{R}{2}\}}
			\frac{d\eta}{1+|\eta|^2}
			=\frac{4\pi}{R^2}\log\left(1+\frac{R^{2}}{4}\right)
			\leq \frac{K\log(2+R)}{R^2}.
		\end{aligned}
		\]
		The same bound holds on $\{|\xi-\eta|\leq \frac{R}{2}\}$,
		by the change of variables $\eta\mapsto\xi-\eta$.
		On the remaining region, both $|\eta|$ and $|\xi-\eta|$
		are greater than $\frac{R}{2}$. Its part with $|\eta|\leq2R$
		satisfies
		\[
		\int_{\substack{\frac{R}{2}<|\eta|\leq2R\\|\xi-\eta|>\frac{R}{2}}}
		\frac{d\eta}
		{(1+|\eta|^2)(1+|\xi-\eta|^2)}
		\leq
		\frac{16}{R^4}\int_{\{|\eta|\leq2R\}}d\eta
		\leq \frac{K}{R^2}.
		\]
		If $|\eta|>2R$, then $|\xi-\eta|\geq\frac{|\eta|}{2}$, so
		\[
		\int_{\{|\eta|>2R\}}
		\frac{d\eta}
		{(1+|\eta|^2)(1+|\xi-\eta|^2)}
		\leq
		4\int_{\{|\eta|>2R\}}|\eta|^{-4}d\eta
		\leq \frac{K}{R^2}.
		\]
		Since $0<s<1$, we have
		\[
		\begin{aligned}
			\log(2+|\xi|)
			&=\int_1^{2+|\xi|}t^{-1} dt
			\leq\int_1^{2+|\xi|}t^{1-2s}dt=\frac{(2+|\xi|)^{2(1-s)}-1}{2(1-s)}
			\leq K_s(1+|\xi|^2)^{1-s},
		\end{aligned}
		\]
		where the last step uses
		$(2+|\xi|)^2\leq5(1+|\xi|^2)$.
		Combining this inequality with the preceding estimates
		for $|\xi|\leq2$ and $|\xi|\geq2$, we obtain
		\begin{equation}\label{eq:exterior-product-kernel}
			\int_{\mathbb R^2}
			\frac{d\eta}
			{(1+|\eta|^2)(1+|\xi-\eta|^2)}
			\leq K\frac{\log(2+|\xi|)}{1+|\xi|^2}
			\leq K_s(1+|\xi|^2)^{-s}.
		\end{equation}
		Multiplying \eqref{eq:exterior-product-fourier} by
		$(1+|\xi|^2)^s$ and using
		\eqref{eq:exterior-product-kernel}, Tonelli's theorem,
		and Plancherel's identity, we obtain
		\begin{align*}
			\|fg\|_{H^s(\mathbb R^2)}^2
			&\leq K_s
			\int_{\mathbb R^2}\int_{\mathbb R^2}
			(1+|\eta|^2)|\widehat f(\eta)|^2
			(1+|\xi-\eta|^2)
			|\widehat g(\xi-\eta)|^2 d\eta d\xi\\
			&=K_s
			\left(
			\int_{\mathbb R^2}
			(1+|\eta|^2)|\widehat f(\eta)|^2 d\eta
			\right)
			\left(
			\int_{\mathbb R^2}
			(1+|\zeta|^2)|\widehat g(\zeta)|^2 d\zeta
			\right).
		\end{align*}
		Thus
		\begin{equation}\label{eq:exterior-product-whole-space}
			\|fg\|_{H^s(\mathbb R^2)}
			\leq K_s
			\|f\|_{H^1(\mathbb R^2)}
			\|g\|_{H^1(\mathbb R^2)}
		\end{equation}
		for smooth, compactly supported functions.
		For general $f,g\in H^1(\mathbb R^2)$, choose
		$f_j,g_j\in C_c^\infty(\mathbb R^2)$ such that
		$f_j\to f$ and $g_j\to g$ in $H^1(\mathbb R^2)$.
		Applying \eqref{eq:exterior-product-whole-space} to
		$f_jg_j-f_\ell g_\ell
		=(f_j-f_\ell)g_j+f_\ell(g_j-g_\ell)$
		gives
		\[
		\begin{aligned}
			\|f_jg_j-f_\ell g_\ell\|_{H^s(\mathbb R^2)}
			\leq K_s
			\|f_j-f_\ell\|_{H^1(\mathbb R^2)}
			\|g_j\|_{H^1(\mathbb R^2)}+K_s
			\|f_\ell\|_{H^1(\mathbb R^2)}
			\|g_j-g_\ell\|_{H^1(\mathbb R^2)}
			\to 0
		\end{aligned}
		\]
		as $j,\ell\to\infty$.
		Hence $f_jg_j$ converges in $H^s(\mathbb R^2)$.
		Moreover, the embedding
		$H^1(\mathbb R^2)\hookrightarrow L^4(\mathbb R^2)$
		and H\"older's inequality imply
		\[
		\begin{aligned}
			\|f_jg_j-fg\|_{L^2(\mathbb R^2)}
			\leq
			\|f_j-f\|_{L^4(\mathbb R^2)}
			\|g_j\|_{L^4(\mathbb R^2)}+
			\|f\|_{L^4(\mathbb R^2)}
			\|g_j-g\|_{L^4(\mathbb R^2)}
			\to 0.
		\end{aligned}
		\]
		The $H^s(\mathbb R^2)$ limit is therefore $fg$.
		Passing to the limit proves
		\eqref{eq:exterior-product-whole-space} for all
		$f,g\in H^1(\mathbb R^2)$.
		
		Finally, apply \eqref{eq:exterior-product-whole-space}
		to $Ef$ and $Eg$. Restriction to $\Omega$ and
		\eqref{eq:exterior-product-extension} give
		\[
		\begin{aligned}
			\|fg\|_{H^s(\Omega)}
			\leq K_s\|(Ef)(Eg)\|_{H^s(\mathbb R^2)}\leq K_s
			\|Ef\|_{H^1(\mathbb R^2)}
			\|Eg\|_{H^1(\mathbb R^2)}\leq K_s
			\|f\|_{H^1(\Omega)}
			\|g\|_{H^1(\Omega)}.
		\end{aligned}
		\]
	\end{proof}
	\section{Velocity estimates in the exterior of the unit disk}
	\label{app:velocity}
	The weighted energy estimates in
	Lemma~\ref{lem:finite-time-uniform-moments} require control
	of $u$ in $L^2(\Omega)$ and $u_r$ in $L^\infty(\Omega)$.
	The vorticity provides these bounds as well as the Sobolev
	estimates needed for local existence and continuation.
	\begin{lemma}
		\label{lem:exterior-radial-biot-savart}
		Suppose that $r^2w\in L^2(\Omega)$.
		There exists a unique $u\in L^2(\Omega;\mathbb R^2)$
		satisfying
		\begin{equation}\label{eq:exterior-velocity-conditions}
			\begin{cases}
				\operatorname{curl}u=w,\quad
				\nabla\cdot u=0
				&\text{\rm in}\quad\Omega,\\
				u\cdot e_r=0
				&\text{\rm on}\quad r=1.
			\end{cases}
		\end{equation}
		The curl equation holds in the sense of distributions.
		For $u\in L^2(\Omega;\mathbb R^2)$, the divergence and
		boundary conditions are imposed by requiring
		\begin{equation}\label{eq:exterior-velocity-weak-condition}
			\int_\Omega u\cdot\nabla\psi dx=0
			\quad\text{for every}\quad\psi\in H^1(\Omega).
		\end{equation}
		The velocity satisfies
		\begin{align}
			\|u\|_{L^2(\Omega)}
			&\leq K\|r^2w\|_{L^2(\Omega)},
			\label{eq:weighted-vorticity-velocity}\\
			\|u_r\|_{L^\infty(\Omega)}
			&\leq K
			\|w\|_{L^1(\Omega)}^{\frac13}
			\|w\|_{L^4(\Omega)}^{\frac23}
			\quad\text{\rm if}\quad w\in L^4(\Omega).
			\label{eq:radial-velocity-interpolation}
		\end{align}
		More generally, let $u\in L^2(\Omega;\mathbb R^2)$
		satisfy \eqref{eq:exterior-velocity-conditions}
		with $w\in L^2(\Omega)$.
		Then
		\begin{equation}\label{eq:exterior-hodge-H1}
			\|u\|_{H^1(\Omega)}
			\leq K\bigl(
			\|u\|_{L^2(\Omega)}+\|w\|_{L^2(\Omega)}
			\bigr).
		\end{equation}
		If additionally $w\in H_0^1(\Omega)$, then
		\begin{equation}\label{eq:exterior-hodge-H2}
			\|u\|_{H^2(\Omega)}
			\leq K\bigl(
			\|u\|_{L^2(\Omega)}+\|w\|_{H^1(\Omega)}
			\bigr).
		\end{equation}
		All constants in
		\eqref{eq:weighted-vorticity-velocity}--\eqref{eq:exterior-hodge-H2} depend only on $\Omega$.
	\end{lemma}
	\begin{proof}
		\noindent\textbf{Existence and uniqueness for \eqref{eq:exterior-velocity-conditions} and estimate
			\eqref{eq:weighted-vorticity-velocity}.}
		First suppose that $w\in C_c^\infty(\Omega)$.
		For $y\in\Omega$, set $y^*=\frac{y}{|y|^2}$ and define
		\begin{equation}\label{eq:exterior-velocity-green}
			G(x,y)=\frac1{2\pi}
			\log\frac{|x-y|}{|y|\,|x-y^*|},
			\qquad
			\varphi(x)=\int_\Omega G(x,y)w(y)dy.
		\end{equation}
		Using
		$\Delta_x((2\pi)^{-1}\log|x-y|)=\delta_y$,
		$|y^*|<1$, and
		$|x-y|=|y|\,|x-y^*|$ when $|x|=1$, we obtain
		\[
		\Delta_xG(x,y)=\delta_y
		\quad\text{in}\quad\Omega,
		\qquad
		G(x,y)=0
		\quad\text{on}\quad|x|=1.
		\]
		The corresponding stream function and velocity satisfy
		\[
		\Delta\varphi=w,\qquad
		\varphi|_{r=1}=0,\qquad
		u=\nabla^\perp\varphi.
		\]
		Choose $R>2$ such that
		$\operatorname{supp}w\subset\{1<|y|<R\}$.
		For $r=|x|\geq2R$ and $y\in\operatorname{supp}w$,
		we have
		\[
		|x-((1-s)y+sy^*)|\geq r-R\geq\frac r2
		\qquad\text{for}\quad 0\leq s\leq1.
		\]
		Differentiating \eqref{eq:exterior-velocity-green}
		and applying the mean value theorem, we obtain
		\[
		\begin{aligned}
			|\nabla_xG(x,y)|
			=\frac1{2\pi}\left|
			\frac{x-y}{|x-y|^2}
			-\frac{x-y^*}{|x-y^*|^2}
			\right|\leq \frac{K|y-y^*|}{r^2}
			\leq\frac{KR}{r^2}.
		\end{aligned}
		\]
		Consequently,
		\begin{equation}\label{eq:exterior-velocity-decay}
			|u(x)|
			\leq\frac{KR}{r^2}\|w\|_{L^1(\Omega)}\qquad\text{for}\quad
			r\geq2R,
		\end{equation}
		and hence
		\[
		\int_{\{r>2R\}}|u(x)|^2 dx
		\leq KR^2\|w\|_{L^1(\Omega)}^2
		\int_{2R}^\infty r^{-3}dr
		<\infty.
		\]
		Since $w\in C_c^\infty(\Omega)$, the potential $\varphi$
		and its gradient are smooth up to the circle $r=1$.
		Thus
		\[
		\int_{\{1<r<2R\}}|u(x)|^2\,dx<\infty,
		\]
		which proves $u\in L^2(\Omega;\mathbb R^2)$.
		Using $u=\nabla^\perp\varphi$,
		$\Delta\varphi=w$, and $\varphi(1,\theta)=0$, we also have
		\[
		\nabla\cdot u=0,\qquad
		\operatorname{curl}u=w
		\quad\text{in}\quad\Omega,
		\qquad
		u_r|_{r=1}
		=-r^{-1}\partial_\theta\varphi|_{r=1}=0.
		\]
		These identities verify
		\eqref{eq:exterior-velocity-conditions}.
		To verify \eqref{eq:exterior-velocity-weak-condition}, let
		$\psi\in H^1(\Omega)$ and choose radial cutoffs $\eta_R$
		with 
		\begin{equation}\label{eq:radial-cutoff} 
			\eta_R=1\quad\text{for}\quad r\leq R,\qquad \eta_R=0\quad\text{for}\quad r\geq2R,\qquad |\nabla\eta_R|\leq \frac{K}{R}.
		\end{equation}
		Integration by parts gives
		\[
		\int_\Omega\eta_Ru\cdot\nabla\psi dx
		=-\int_\Omega\psi u\cdot\nabla\eta_R dx.
		\]
		The term on the right satisfies
		\[
		\left|\int_\Omega\psi u\cdot\nabla\eta_R dx\right|
		\leq\frac KR\|\psi\|_{L^2(\{r>R\})}
		\|u\|_{L^2(\{r>R\})}.
		\]
		Letting $R\to\infty$ proves
		\eqref{eq:exterior-velocity-weak-condition}.
		
		We now prove \eqref{eq:weighted-vorticity-velocity}
		by estimating the zero and nonzero angular modes.
		For the nonzero modes, the angular Poincar\'e inequality gives
		\[
		\int_\Omega r^{-2}|\varphi_{\neq}|^2 dx
		\leq
		\int_\Omega r^{-2}|\partial_\theta\varphi_{\neq}|^2 dx
		\leq
		\|\nabla\varphi_{\neq}\|_{L^2(\Omega)}^2.
		\]
		Let $0\leq\eta_R\leq1$ be the radial cutoffs in
		\eqref{eq:radial-cutoff}.
		Since $\varphi_{\neq}$ has zero trace at $r=1$,
		we have $\eta_R^2\varphi_{\neq}\in H_0^1(\Omega)$.
		Testing $\Delta\varphi_{\neq}=w_{\neq}$ with
		$-\eta_R^2\varphi_{\neq}$ gives
		\[
		\begin{aligned}
			\int_\Omega\eta_R^2|\nabla\varphi_{\neq}|^2 dx
			=-\int_\Omega\eta_R^2w_{\neq}\varphi_{\neq}dx-2\int_\Omega\eta_R\varphi_{\neq}
			\nabla\eta_R\cdot\nabla\varphi_{\neq}dx.
		\end{aligned}
		\]
		The derivatives of $\eta_R$ are supported in
		$\{R<r<2R\}$, where $R^{-1}\leq2r^{-1}$.
		Hence \eqref{eq:radial-cutoff} and the
		Cauchy--Schwarz inequality imply
		\[
		\begin{aligned}
			\left|2\int_\Omega\eta_R\varphi_{\neq}
			\nabla\eta_R\cdot\nabla\varphi_{\neq}dx\right|&\leq K\int_{\{R<r<2R\}}
			r^{-1}|\varphi_{\neq}||\nabla\varphi_{\neq}|dx\\&\leq
			K\|r^{-1}\varphi_{\neq}\|_{L^2(\{r>R\})}
			\|\nabla\varphi_{\neq}\|_{L^2(\{r>R\})}
			\to 0
			\quad\text{as}\quad R\to\infty,
		\end{aligned}
		\]
		because $r^{-1}\varphi_{\neq}$ and
		$\nabla\varphi_{\neq}$ belong to $L^2(\Omega)$.
		Moreover,
		\[
		\int_\Omega|w_{\neq}\varphi_{\neq}|dx
		\leq
		\|rw_{\neq}\|_{L^2(\Omega)}
		\|r^{-1}\varphi_{\neq}\|_{L^2(\Omega)}
		<\infty.
		\]
		The dominated convergence theorem therefore permits
		passage to the limit as $R\to\infty$ in the tested equation.
		Using the angular Poincar\'e inequality above, we obtain
		\[
		\begin{aligned}
			\|\nabla\varphi_{\neq}\|_{L^2(\Omega)}^2
			=-\int_\Omega w_{\neq}\varphi_{\neq}dx\leq
			\|rw_{\neq}\|_{L^2(\Omega)}
			\|r^{-1}\varphi_{\neq}\|_{L^2(\Omega)}\leq
			\|rw_{\neq}\|_{L^2(\Omega)}
			\|\nabla\varphi_{\neq}\|_{L^2(\Omega)}.
		\end{aligned}
		\]
		Therefore,
		\[
		\|u_{\neq}\|_{L^2(\Omega)}
		=\|\nabla\varphi_{\neq}\|_{L^2(\Omega)}
		\leq\|rw_{\neq}\|_{L^2(\Omega)}.
		\]
		For the zero mode, incompressibility and the boundary
		condition give
		$\partial_r\bigl(r(u_r)_0\bigr)=0$ and $(u_r)_0(1)=0$,
		so $(u_r)_0=0$.
		Taking the angular mean of the curl equation in
		\eqref{eq:exterior-velocity-conditions}, written as
		$
		w=\frac1r\partial_r(ru_\theta)
		-\frac1r\partial_\theta u_r$,
		gives
		\begin{equation}\label{eq:radial}
			\partial_r\bigl(r(u_\theta)_0\bigr)=rw_0.
		\end{equation}
		Moreover, \eqref{eq:exterior-velocity-decay} implies
		\[
		\bigl|r(u_\theta)_0(r)\bigr|
		\leq\frac{r}{2\pi}\int_0^{2\pi}|u(r,\theta)|d\theta
		\leq\frac{KR}{r}\|w\|_{L^1(\Omega)}
		\to 0
		\quad\text{as}\quad r\to\infty.
		\]
		Integrating \eqref{eq:radial} over $[r,S]$, we obtain
		\[
		S(u_\theta)_0(S)-r(u_\theta)_0(r)
		=\int_r^S s w_0(s) ds.
		\]
		Letting $S\to\infty$ yields
		\[
		(u_\theta)_0(r)
		=-\frac1r\int_r^\infty s w_0(s)ds.
		\]
		By the Cauchy--Schwarz inequality, we get
		\[
		\begin{aligned}
			\left|\int_r^\infty s w_0(s)ds\right|^2
			\leq
			\left(\int_r^\infty s^{-3}ds\right)
			\left(\int_r^\infty s^5|w_0(s)|^2 ds\right)=\frac1{2r^2}
			\int_r^\infty s^5|w_0(s)|^2 ds.
		\end{aligned}
		\]
		Integration in $r$ and Tonelli's theorem give
		\[
		\begin{aligned}
			\int_1^\infty r|(u_\theta)_0(r)|^2 dr
			&\leq\frac12\int_1^\infty r^{-3}
			\int_r^\infty s^5|w_0(s)|^2 dsdr\\
			&=\frac14\int_1^\infty
			(1-s^{-2})s^5|w_0(s)|^2 ds\leq\frac14\int_1^\infty s^5|w_0(s)|^2 ds.
		\end{aligned}
		\]
		Multiplying by $2\pi$, we conclude that
		\[
		\begin{aligned}
			\|u\|_{L^2(\Omega)}^2
			=\|u_0\|_{L^2(\Omega)}^2
			+\|u_{\neq}\|_{L^2(\Omega)}^2\leq\frac14\|r^2w_0\|_{L^2(\Omega)}^2
			+\|rw_{\neq}\|_{L^2(\Omega)}^2\leq\|r^2w\|_{L^2(\Omega)}^2.
		\end{aligned}
		\]
		To extend the construction and estimate to general
		$r^2w\in L^2(\Omega)$,
		choose $w_j\in C_c^\infty(\Omega)$ such that
		$\|r^2(w_j-w)\|_{L^2(\Omega)}
		\to 0$ as $j\to\infty$.
		Using the Green function in
		\eqref{eq:exterior-velocity-green}, define
		\[
		u_j(x)=\int_\Omega
		\nabla_x^\perp G(x,y)w_j(y)dy.
		\]
		By linearity and the estimate
		\eqref{eq:weighted-vorticity-velocity} already proved
		for smooth vorticities, we obtain, as $j,\ell\to\infty$,
		\begin{equation}\label{eq:exterior-velocity-approximation}
			\begin{aligned}
				\|u_j-u_\ell\|_{L^2(\Omega)}&\leq K\|r^2(w_j-w_\ell)\|_{L^2(\Omega)}\\
				&\leq K\Bigl(
				\|r^2(w_j-w)\|_{L^2(\Omega)}
				+\|r^2(w_\ell-w)\|_{L^2(\Omega)}
				\Bigr)\to 0.
			\end{aligned}
		\end{equation}
		The completeness of $L^2(\Omega;\mathbb R^2)$ therefore
		gives $u\in L^2(\Omega;\mathbb R^2)$ such that
		$u_j\to u$ in $L^2(\Omega;\mathbb R^2)$
		as $j\to\infty$.
		Since $r\geq1$, we arrive at
		\[
		\|w_j-w\|_{L^2(\Omega)}
		\leq\|r^2(w_j-w)\|_{L^2(\Omega)}
		\to 0.
		\]
		Thus, for every $\eta\in C_c^\infty(\Omega)$,
		\[
		\begin{aligned}
			-\int_\Omega u\cdot\nabla^\perp\eta dx
			=\lim_{j\to\infty}
			\left(-\int_\Omega u_j\cdot\nabla^\perp\eta dx\right)=\lim_{j\to\infty}\int_\Omega w_j\eta dx
			=\int_\Omega w\eta dx.
		\end{aligned}
		\]
		Hence $\operatorname{curl}u=w$ in the sense of distributions.
		Moreover, for every $\psi\in H^1(\Omega)$,
		\[
		\int_\Omega u\cdot\nabla\psi dx
		=\lim_{j\to\infty}
		\int_\Omega u_j\cdot\nabla\psi dx=0.
		\]
		These identities establish
		\eqref{eq:exterior-velocity-conditions}, with the divergence
		and boundary conditions given by
		\eqref{eq:exterior-velocity-weak-condition}.
		Strong convergence also gives
		\[
		\begin{aligned}
			\|u\|_{L^2(\Omega)}
			=\lim_{j\to\infty}\|u_j\|_{L^2(\Omega)}\leq K\lim_{j\to\infty}
			\|r^2w_j\|_{L^2(\Omega)}
			=K\|r^2w\|_{L^2(\Omega)}.
		\end{aligned}
		\]
		Thus a velocity satisfying
		\eqref{eq:exterior-velocity-conditions} and
		\eqref{eq:weighted-vorticity-velocity} exists.
		
		It remains to prove uniqueness.
		For this purpose, let an arbitrary
		$u\in L^2(\Omega;\mathbb R^2)$ satisfy
		\eqref{eq:exterior-velocity-conditions}.
		We first construct a stream function for this velocity.
		Define
		\[
		U(x)=
		\begin{cases}
			u(x),&|x|>1,\\
			0,&|x|\leq1.
		\end{cases}
		\]
		For every $\eta\in C_c^\infty(\mathbb R^2)$,
		\eqref{eq:exterior-velocity-weak-condition} gives
		\[
		\int_{\mathbb R^2}U\cdot\nabla\eta dx
		=\int_\Omega u\cdot\nabla\eta dx=0.
		\]
		The field $(U_2,-U_1)$ belongs to
		$L^2_{\mathrm{loc}}(\mathbb R^2;\mathbb R^2)$ and satisfies
		$\operatorname{curl}(U_2,-U_1)
		=-\nabla\cdot U=0$
		in the sense of distributions. Since the plane is simply connected, there exists
		$\varphi\in H^1_{\mathrm{loc}}(\mathbb R^2)$ with
		$\nabla\varphi=(U_2,-U_1)$.
		In the unit disk, $U=0$, so $\nabla\varphi=0$ and
		$\varphi$ is constant there. Replacing $\varphi$ by
		$\varphi-\frac1\pi\int_{B_1(0)}\varphi(y)dy$
		leaves its gradient unchanged and gives
		$\varphi=0$ in $B_1(0)$.
		Consequently,
		\begin{equation}\label{eq:exterior-stream-function}
			\begin{gathered}
				u=\nabla^\perp\varphi,\qquad
				\Delta\varphi
				=\partial_{x_1}u_2-\partial_{x_2}u_1=w
				\quad\text{in}\quad\Omega,\qquad \varphi|_{r=1}=0,\\
				\|\nabla\varphi\|_{L^2(\mathbb R^2)}
				=\|U\|_{L^2(\mathbb R^2)}
				=\|u\|_{L^2(\Omega)}.
			\end{gathered}
		\end{equation}
		We now show that $w=0$ implies $u=0$.
		In this case, \eqref{eq:exterior-stream-function} gives
		\[
		\partial_r^2\varphi+\frac1r\partial_r\varphi
		+\frac1{r^2}\partial_\theta^2\varphi=0.
		\]
		Taking the $k$th angular Fourier coefficient and integrating
		by parts in $\theta$, we obtain
		\begin{equation}\label{eq:exterior-harmonic-mode-equation}
			\partial_r^2\widehat\varphi_k
			+\frac1r\partial_r\widehat\varphi_k
			-\frac{k^2}{r^2}\widehat\varphi_k=0,
			\qquad
			\widehat\varphi_k(1)=0.
		\end{equation}
		For $k\neq0$, substitution of $r^\lambda$ into this
		equation gives $\lambda^2-k^2=0$.
		For $k=0$, the equation becomes
		$\partial_r(r\partial_r\widehat\varphi_0)=0$.
		Consequently,
		\[
		\widehat\varphi_k(r)
		=\begin{cases}
			a_kr^{|k|}+b_kr^{-|k|},&\quad k\neq0,\\
			a_0+b_0\log r,&\quad k=0.
		\end{cases}
		\]
		On the other hand, Parseval's identity and
		\eqref{eq:exterior-stream-function} give
		\begin{equation}\label{eq:exterior-stream-fourier-energy}
			\begin{aligned}
				\sum_{k\in\mathbb Z}\int_1^\infty
				\left(
				|\partial_r\widehat\varphi_k|^2
				+\frac{k^2}{r^2}|\widehat\varphi_k|^2
				\right)rdr=
				\frac1{2\pi}\|\nabla\varphi\|_{L^2(\Omega)}^2
				=\frac1{2\pi}\|u\|_{L^2(\Omega)}^2<\infty.
			\end{aligned}
		\end{equation}
		For $k\neq0$, direct substitution yields
		\[
		\begin{aligned}
			r|\partial_r\widehat\varphi_k|^2
			+\frac{k^2}{r}|\widehat\varphi_k|^2=
			2k^2\left(
			|a_k|^2r^{2|k|-1}
			+|b_k|^2r^{-2|k|-1}
			\right).
		\end{aligned}
		\]
		By \eqref{eq:exterior-stream-fourier-energy},
		for every $R>1$, there holds
		\[
		2k^2|a_k|^2\int_1^R r^{2|k|-1}dr
		\leq\frac1{2\pi}\|u\|_{L^2(\Omega)}^2.
		\]
		Letting $R\to\infty$ gives $a_k=0$.
		The boundary condition in
		\eqref{eq:exterior-harmonic-mode-equation} then gives
		$b_k=0$.
		For $k=0$, \eqref{eq:exterior-stream-fourier-energy} implies
		\[
		|b_0|^2\log R
		=\int_1^R r|\partial_r\widehat\varphi_0|^2 dr
		\leq\frac1{2\pi}\|u\|_{L^2(\Omega)}^2.
		\]
		Letting $R\to\infty$ gives $b_0=0$, while
		$\widehat\varphi_0(1)=0$ gives $a_0=0$.
		Hence every angular Fourier coefficient vanishes,
		so $\varphi=0$ and $u=\nabla^\perp\varphi=0$.
		Applying this conclusion to the difference of two velocities
		satisfying \eqref{eq:exterior-velocity-conditions}
		with the same $w$ proves uniqueness.
		
		\noindent\textbf{Estimate \eqref{eq:radial-velocity-interpolation}.}
		First let $w\in C_c^\infty(\Omega)$.
		Write
		$x=(r\cos\theta,r\sin\theta)$ and
		$y=(\rho\cos\alpha,\rho\sin\alpha)$.
		Differentiating \eqref{eq:exterior-velocity-green}
		and using $u_r=-r^{-1}\partial_\theta\varphi$, we obtain
		\[
		u_r(x)=\int_\Omega K_r(x,y)w(y)dy,
		\]
		where
		\[
		K_r(x,y)
		=-\frac1r\partial_\theta G(x,y)
		=-\frac{\sin(\theta-\alpha)}{2\pi}
		\left(
		\frac{\rho}{|x-y|^2}
		-\frac{\rho^{-1}}{|x-y^*|^2}
		\right).
		\]
		For $x\neq y$, the identity
		$|x-y^*|^2-\rho^{-2}|x-y|^2=(r^2-1)(1-\rho^{-2})\geq0$
		implies
		$\frac{\rho^{-1}}{|x-y^*|^2}
		\leq\frac{\rho}{|x-y|^2}$,
		and hence
		\[
		0\leq
		\frac{\rho}{|x-y|^2}
		-\frac{\rho^{-1}}{|x-y^*|^2}
		\leq\frac{\rho}{|x-y|^2}.
		\]
		Moreover,
		$|x-y|^2
		=(r-\rho\cos(\theta-\alpha))^2
		+\rho^2\sin^2(\theta-\alpha)$,
		so $\rho|\sin(\theta-\alpha)|\leq|x-y|$.
		Consequently,
		\begin{equation}\label{eq:exterior-radial-kernel-bound}
			|K_r(x,y)|
			\leq
			\frac{\rho|\sin(\theta-\alpha)|}
			{2\pi|x-y|^2}
			\leq\frac1{2\pi|x-y|}.
		\end{equation}
		For any $L>0$, polar integration gives
		\[
		\int_{\{|z|<L\}}|z|^{-\frac43}dz
		=2\pi\int_0^L s^{-\frac13}ds
		=3\pi L^{\frac23}.
		\]
		Splitting the integral into $|x-y|<L$ and $|x-y|\geq L$,
		and applying H\"older's inequality, we obtain from
		\eqref{eq:exterior-radial-kernel-bound}
		\begin{align}\label{ineq:ur}
			|u_r(x)|
			&\leq
			K\int_{\Omega\cap\{|x-y|<L\}}
			\frac{|w(y)|}{|x-y|}dy
			+KL^{-1}\|w\|_{L^1(\Omega)}\notag\\
			&\leq
			K\left(
			\int_{\{|z|<L\}}|z|^{-\frac43}dz
			\right)^{\frac34}\|w\|_{L^4(\Omega)}
			+KL^{-1}\|w\|_{L^1(\Omega)}\\\notag
			&\leq
			K L^{\frac12}\|w\|_{L^4(\Omega)}
			+K L^{-1}\|w\|_{L^1(\Omega)}.
		\end{align}
		If $w=0$, uniqueness gives $u=0$.
		Otherwise, choose
		$L=
		\left(
		\frac{\|w\|_{L^1(\Omega)}}{\|w\|_{L^4(\Omega)}}
		\right)^{\frac23}$.
		For this choice, we have
		\[
		L^{\frac12}\|w\|_{L^4(\Omega)}
		=L^{-1}\|w\|_{L^1(\Omega)}
		=\|w\|_{L^1(\Omega)}^{\frac13}
		\|w\|_{L^4(\Omega)}^{\frac23}.
		\]
		Substitution into \eqref{ineq:ur} proves
		\eqref{eq:radial-velocity-interpolation}
		for smooth, compactly supported vorticities.
		
		For general $r^2w\in L^2(\Omega)$ with $w\in L^4(\Omega)$,
		choose $w_j\in C_c^\infty(\Omega)$ as in the proof of
		\eqref{eq:weighted-vorticity-velocity}, with the additional
		$L^4(\Omega)$ convergence, so that
		\[
		\|r^2(w_j-w)\|_{L^2(\Omega)}
		+\|w_j-w\|_{L^4(\Omega)}
		\to 0
		\quad\text{as}\quad j\to\infty.
		\]
		The inequality
		\[
		\|w_j-w\|_{L^1(\Omega)}
		\leq\sqrt{\pi}\|r^2(w_j-w)\|_{L^2(\Omega)}
		\]
		also gives $w_j\to w$ in $L^1(\Omega)$.
		Let $u_j$ be the corresponding velocities.
		By \eqref{eq:weighted-vorticity-velocity},
		$u_j\to u$ in $L^2(\Omega;\mathbb R^2)$.
		After passing to a subsequence,
		$(u_j)_r\to u_r$ almost everywhere in $\Omega$.
		The estimate already proved for smooth vorticities
		therefore yields
		\[
		\begin{aligned}
			\|u_r\|_{L^\infty(\Omega)}
			\leq\liminf_{j\to\infty}
			\|(u_j)_r\|_{L^\infty(\Omega)}\leq K\lim_{j\to\infty}
			\|w_j\|_{L^1(\Omega)}^{\frac13}
			\|w_j\|_{L^4(\Omega)}^{\frac23}=K\|w\|_{L^1(\Omega)}^{\frac13}
			\|w\|_{L^4(\Omega)}^{\frac23}.
		\end{aligned}
		\]
		
		\noindent\textbf{Estimate
			\eqref{eq:exterior-hodge-H1}.}
		Let $u\in L^2(\Omega;\mathbb R^2)$ satisfy
		\eqref{eq:exterior-velocity-conditions} with
		$w\in L^2(\Omega)$, and let $\varphi$ be the stream
		function in \eqref{eq:exterior-stream-function}.
		We first estimate $\varphi$ on $B_4(0)$.
		Since $\varphi=0$ in $B_1(0)$, the fundamental theorem
		of calculus and the Cauchy--Schwarz inequality give,
		for almost every $\theta$ and $1<r<4$,
		\[
		\begin{aligned}
			|\varphi(r,\theta)|^2
			&=\left|\int_1^r
			\partial_s\varphi(s,\theta)ds\right|^2\leq
			\left(\int_1^r s^{-1}ds\right)
			\left(\int_1^r
			s|\partial_s\varphi(s,\theta)|^2 ds\right)\\&\leq\log4\int_1^4
			s|\nabla\varphi(s,\theta)|^2 ds.
		\end{aligned}
		\]
		Integration with respect to $rdrd\theta$, together
		with \eqref{eq:exterior-stream-function}, yields
		\begin{equation}\label{eq:exterior-stream-local-bound}
			\|\varphi\|_{H^1(B_4(0))}
			\leq K\|\nabla\varphi\|_{L^2(B_4(0))}
			\leq K\|u\|_{L^2(\Omega)}.
		\end{equation}
		Choose a fixed smooth radial function $\chi$ such that
		\begin{equation}\label{define:kappa}
			0\leq\chi\leq1,\qquad \chi=0\quad\text{for}\quad r\leq2,\qquad \chi=1\quad\text{for}\quad r\geq3.
		\end{equation}
		We estimate $(1-\chi)u$ and $\chi u$ separately.
		Noting that $\varphi|_{r=1}=0$ and $1-\chi$ vanishes for
		$r\geq3$, we have $(1-\chi)\varphi\in H_0^1(\Omega)$.
		Using $\Delta\varphi=w$ in
		\eqref{eq:exterior-stream-function}, the product rule gives
		\begin{equation}\label{eq:exterior-local-stream-equation}
			\begin{aligned}
				(1-\Delta)((1-\chi)\varphi)
				=(1-\chi)(\varphi-w)+2\nabla\chi\cdot\nabla\varphi
				+(\Delta\chi)\varphi.
			\end{aligned}
		\end{equation}
		The derivatives of $\chi$ are bounded and supported
		in $\{2\leq r\leq3\}$.
		Thus \eqref{eq:exterior-local-stream-equation} and
		\eqref{eq:exterior-stream-local-bound} imply
		\[
		\begin{aligned}
			&\|(1-\Delta)((1-\chi)\varphi)\|_{L^2(\Omega)}\\
			\leq~&
			\|(1-\chi)\varphi\|_{L^2(\Omega)}
			+\|(1-\chi)w\|_{L^2(\Omega)}
			+2\|\nabla\chi\cdot\nabla\varphi\|_{L^2(\Omega)}
			+\|(\Delta\chi)\varphi\|_{L^2(\Omega)}\\
			\leq~&
			K\bigl(
			\|\varphi\|_{H^1(B_4(0))}+\|w\|_{L^2(\Omega)}
			\bigr)\leq
			K\bigl(
			\|u\|_{L^2(\Omega)}+\|w\|_{L^2(\Omega)}
			\bigr).
		\end{aligned}
		\]
		Applying \eqref{eq:screened-W2p} with $p=2$, we obtain
		\[
		\|(1-\chi)\varphi\|_{H^2(\Omega)}
		\leq K\bigl(
		\|u\|_{L^2(\Omega)}+\|w\|_{L^2(\Omega)}
		\bigr).
		\]
		Moreover, $u=\nabla^\perp\varphi$ gives
		$(1-\chi)u
		=\nabla^\perp((1-\chi)\varphi)
		+\varphi\nabla^\perp\chi$.
		Using the preceding estimate and
		\eqref{eq:exterior-stream-local-bound}, we conclude that
		\begin{equation}\label{eq:exterior-velocity-near-H1}
			\begin{aligned}
				\|(1-\chi)u\|_{H^1(\Omega)}
				&\leq
				K\|(1-\chi)\varphi\|_{H^2(\Omega)}
				+\|\varphi\nabla^\perp\chi\|_{H^1(\Omega)}\\
				&\leq
				K\bigl(
				\|(1-\chi)\varphi\|_{H^2(\Omega)}
				+\|\varphi\|_{H^1(B_4(0))}
				\bigr)\\
				&\leq K\bigl(
				\|u\|_{L^2(\Omega)}+\|w\|_{L^2(\Omega)}
				\bigr).
			\end{aligned}
		\end{equation}
		For $\chi u$, extend it by zero to $\mathbb R^2$.
		Since $\chi u$ vanishes near $r=1$, the product rule
		holds on $\mathbb R^2$ in the sense of distributions, which yields
		\[
		\nabla\cdot(\chi u)=\nabla\chi\cdot u,
		\qquad
		\operatorname{curl}(\chi u)
		=\chi w+\nabla^\perp\chi\cdot u.
		\]
		Consequently,
		\[
		\begin{aligned}
			\|\nabla\cdot(\chi u)\|_{L^2(\mathbb R^2)}
			+\|\operatorname{curl}(\chi u)\|_{L^2(\mathbb R^2)}\leq
			K\bigl(
			\|u\|_{L^2(\Omega)}+\|w\|_{L^2(\Omega)}
			\bigr).
		\end{aligned}
		\]
		Taking the Fourier transform with respect to
		$x\in\mathbb R^2$, we use the identity
		\[
		|\xi|^2|\widehat{\chi u}(\xi)|^2
		=
		|\xi\cdot\widehat{\chi u}(\xi)|^2
		+
		|\xi_1\widehat{\chi u}_2(\xi)
		-\xi_2\widehat{\chi u}_1(\xi)|^2.
		\]
		The Fourier identity and Plancherel's theorem give
		\[
		\begin{aligned}
			\int_{\mathbb R^2}
			(1+|\xi|^2)|\widehat{\chi u}(\xi)|^2 d\xi
			&\leq K\Bigl(
			\|\chi u\|_{L^2(\mathbb R^2)}^2
			+\|\nabla\cdot(\chi u)\|_{L^2(\mathbb R^2)}^2
			+\|\operatorname{curl}(\chi u)\|_{L^2(\mathbb R^2)}^2
			\Bigr)\\
			&\leq K\bigl(
			\|u\|_{L^2(\Omega)}^2
			+\|w\|_{L^2(\Omega)}^2
			\bigr)<\infty.
		\end{aligned}
		\]
		Hence $\chi u\in H^1(\mathbb R^2;\mathbb R^2)$, and
		\[
		\|\chi u\|_{H^1(\mathbb R^2)}
		\leq K\bigl(
		\|u\|_{L^2(\Omega)}
		+\|w\|_{L^2(\Omega)}
		\bigr).
		\]
		Combining this estimate with
		\eqref{eq:exterior-velocity-near-H1}, we obtain
		\[
		\begin{aligned}
			\|u\|_{H^1(\Omega)}
			\leq
			\|(1-\chi)u\|_{H^1(\Omega)}
			+\|\chi u\|_{H^1(\mathbb R^2)}\leq K\bigl(
			\|u\|_{L^2(\Omega)}+\|w\|_{L^2(\Omega)}
			\bigr).
		\end{aligned}
		\]
		
		\noindent\textbf{Estimate \eqref{eq:exterior-hodge-H2}.}
		Suppose additionally that $w\in H_0^1(\Omega)$.
		By \eqref{eq:exterior-hodge-H1},
		$u\in H^1(\Omega;\mathbb R^2)$.
		We retain the cutoff function $\chi$ defined in \eqref{define:kappa}.
		Using $\nabla\varphi=(u_2,-u_1)$ and
		\eqref{eq:exterior-stream-local-bound}, we obtain
		\begin{equation}\label{eq:exterior-stream-local-H2}
			\begin{aligned}
				\|\varphi\|_{H^2(\Omega\cap B_4(0))}
				&\leq K\Bigl(
				\|\varphi\|_{L^2(B_4(0))}
				+\|\nabla\varphi\|_{H^1(\Omega\cap B_4(0))}
				\Bigr)\\
				&\leq K\Bigl(
				\|\varphi\|_{L^2(B_4(0))}
				+\|u\|_{H^1(\Omega)}
				\Bigr)\leq K\|u\|_{H^1(\Omega)}.
			\end{aligned}
		\end{equation}
		Since $\chi$ is fixed and smooth, with its derivatives
		supported in $\{2\leq r\leq3\}$, the product rule in
		\eqref{eq:exterior-local-stream-equation} gives
		\[
		\begin{aligned}
			&\|(1-\Delta)((1-\chi)\varphi)\|_{H^1(\Omega)}\\
			\leq~&\|(1-\chi)(\varphi-w)\|_{H^1(\Omega)}
			+2\|\nabla\chi\cdot\nabla\varphi\|_{H^1(\Omega)}
			+\|(\Delta\chi)\varphi\|_{H^1(\Omega)}\\
			\leq~& K\Bigl(
			\|\varphi\|_{H^2(\Omega\cap B_4(0))}
			+\|w\|_{H^1(\Omega)}
			\Bigr)\leq K\bigl(
			\|u\|_{H^1(\Omega)}+\|w\|_{H^1(\Omega)}
			\bigr),
		\end{aligned}
		\]
		where the last inequality follows from
		\eqref{eq:exterior-stream-local-H2}.
		Applying the $H^3(\Omega)$ estimate in
		Lemma~\ref{lem:screened-regularity} to
		\eqref{eq:exterior-local-stream-equation}, we conclude that
		\[
		\|(1-\chi)\varphi\|_{H^3(\Omega)}
		\leq K\bigl(
		\|u\|_{H^1(\Omega)}+\|w\|_{H^1(\Omega)}
		\bigr).
		\]
		Furthermore,
		$(1-\chi)u
		=\nabla^\perp((1-\chi)\varphi)
		+\varphi\nabla^\perp\chi$.
		The preceding estimate and
		\eqref{eq:exterior-stream-local-H2} therefore yield
		\begin{equation}\label{eq:exterior-velocity-near-H2}
			\begin{aligned}
				\|(1-\chi)u\|_{H^2(\Omega)}
				&\leq
				K\|(1-\chi)\varphi\|_{H^3(\Omega)}
				+\|\varphi\nabla^\perp\chi\|_{H^2(\Omega)}\\
				&\leq K\Bigl(
				\|(1-\chi)\varphi\|_{H^3(\Omega)}
				+\|\varphi\|_{H^2(\Omega\cap B_4(0))}
				\Bigr)\\
				&\leq K\bigl(
				\|u\|_{H^1(\Omega)}+\|w\|_{H^1(\Omega)}
				\bigr).
			\end{aligned}
		\end{equation}
		For $\chi u$, the conditions in
		\eqref{eq:exterior-velocity-conditions} imply
		$\Delta u
		=\nabla(\nabla\cdot u)
		+\nabla^\perp(\operatorname{curl}u)
		=\nabla^\perp w$
		in the sense of distributions in $\Omega$.
		Since $\chi u$ vanishes near $r=1$, its zero extension
		to $\mathbb R^2$ satisfies
		\[
		\Delta(\chi u)
		=\chi\nabla^\perp w
		+2\sum_{i=1}^2
		(\partial_{x_i}\chi)\partial_{x_i}u
		+(\Delta\chi)u
		\]
		in the sense of distributions on $\mathbb R^2$.
		Consequently,
		\begin{equation}\label{eq:exterior-cutoff-laplacian-H2}
			\begin{aligned}
				\|\Delta(\chi u)\|_{L^2(\mathbb R^2)}
				&\leq
				\|\chi\nabla^\perp w\|_{L^2(\Omega)}
				+2\sum_{i=1}^2
				\|(\partial_{x_i}\chi)\partial_{x_i}u\|_{L^2(\Omega)}
				+\|(\Delta\chi)u\|_{L^2(\Omega)}\\
				&\leq K\bigl(
				\|\nabla w\|_{L^2(\Omega)}
				+\|u\|_{H^1(\Omega)}
				\bigr).
			\end{aligned}
		\end{equation}
		Taking the Fourier transform with respect to
		$x\in\mathbb R^2$, we have
		$\widehat{\Delta(\chi u)}(\xi)=-|\xi|^2\widehat{\chi u}(\xi)$.
		Plancherel's theorem and
		\eqref{eq:exterior-cutoff-laplacian-H2} give
		\[
		\int_{\mathbb R^2}
		|\xi|^4|\widehat{\chi u}(\xi)|^2 d\xi
		\leq K\|\Delta(\chi u)\|_{L^2(\mathbb R^2)}^2
		<\infty.
		\]
		Together with $\chi u\in L^2(\mathbb R^2;\mathbb R^2)$,
		this proves $\chi u\in H^2(\mathbb R^2;\mathbb R^2)$.
		Using $(1+|\xi|^2)^2\leq2(1+|\xi|^4)$, we also obtain
		\[
		\begin{aligned}
			\|\chi u\|_{H^2(\mathbb R^2)}
			\leq K\bigl(
			\|\chi u\|_{L^2(\mathbb R^2)}
			+\|\Delta(\chi u)\|_{L^2(\mathbb R^2)}
			\bigr)\leq K\bigl(
			\|u\|_{H^1(\Omega)}+\|w\|_{H^1(\Omega)}
			\bigr).
		\end{aligned}
		\]
		Combining this estimate with
		\eqref{eq:exterior-velocity-near-H2} and then using
		\eqref{eq:exterior-hodge-H1}, we conclude that
		\[
		\begin{aligned}
			\|u\|_{H^2(\Omega)}
			&\leq
			\|(1-\chi)u\|_{H^2(\Omega)}
			+\|\chi u\|_{H^2(\mathbb R^2)}\\
			&\leq K\bigl(
			\|u\|_{H^1(\Omega)}+\|w\|_{H^1(\Omega)}
			\bigr)\leq K\bigl(
			\|u\|_{L^2(\Omega)}+\|w\|_{H^1(\Omega)}
			\bigr).
		\end{aligned}
		\]
	\end{proof}
	\section{A priori estimates for the cell density}
	\label{sec:density}
	This section develops the density estimates used in
	Theorem~\ref{thm:main}. A uniform $L^2(\Omega)$ bound yields
	an $L^\infty(\Omega)$ bound. For the radial problem
	\eqref{eq:radial-reference-problem}, the density is globally
	bounded and decays, while an $L^2(\Omega)$ smallness
	condition ensures global existence for
	\eqref{eq:physical-scalars-moments} with the boundary conditions \eqref{eq:main-boundary}.
	\subsection{Uniform $L^\infty$ bounds for the density}
	An $L^2(\Omega)$ bound for the density yields
	an $L^\infty(\Omega)$ bound by Moser iteration.
	The constant in this estimate is independent of
	the length of the time interval and of the
	divergence-free velocity field in the density equation.
	\begin{lemma}\label{lem:exterior-moser}
		Let $0<T\leq\infty$, and let $n\geq0$ solve
		\begin{equation}\label{eq:exterior-moser-system}
			\begin{cases}
				\partial_tn+V\cdot\nabla n
				=\Delta n-\nabla\cdot(n\nabla c)
				&\text{\rm in}\quad(0,T)\times\Omega,\\
				-\Delta c+c=n
				&\text{\rm in}\quad(0,T)\times\Omega,\\
				n=c=0
				&\text{\rm on}\quad(0,T)\times\partial\Omega,\\
				n(0)=n_{\mathrm{in}}
				&\text{\rm in}~\quad\Omega.
			\end{cases}
		\end{equation}
		Assume that $c(t)\in H_0^1(\Omega)$ and that $n$ has
		the density regularity stated in
		Proposition~\ref{prop:local-strong-theory} on every interval
		$[0,\tau]$ with $0<\tau<T$.
		Suppose also that
		\[
		\nabla\cdot V=0,\quad
		V\cdot e_r=0\quad\text{\rm on}\quad r=1,\quad
		V\in L^\infty(0,\tau;L^\infty(\Omega))
		\quad\text{\rm for every}\quad 0<\tau<T.
		\]
		If
		$n_{\mathrm{in}}\in L^2(\Omega)\cap L^\infty(\Omega)$ and
		$$\sup_{0\leq t<T}\|n(t)\|_{L^2(\Omega)}\leq B_1$$
		for some $B_1>0$,
		then
		\[
		\sup_{0\leq t<T}\|n(t)\|_{L^\infty(\Omega)}
		\leq
		K\bigl(\|n_{\mathrm{in}}\|_{L^\infty(\Omega)}+B_1\bigr)
		(1+B_1^2),
		\]
		where $K$ depends only on $\Omega$.
	\end{lemma}
	
	\begin{proof}
		Fix $0<\tau<T$.
		By the assumed regularity and Sobolev embedding, we have
		$n\in C([0,\tau];L^p(\Omega))$ for $2\leq p<\infty$.
		In particular, these norms are bounded on $[0,\tau]$.
		Lemma~\ref{lem:screened-regularity} and Sobolev embedding
		give
		\begin{equation}\label{eq:exterior-elliptic-L4}
			\|c(t)\|_{H^2(\Omega)}
			\leq K\|n(t)\|_{L^2(\Omega)},
			\qquad
			\|\nabla c(t)\|_{L^4(\Omega)}\leq KB_1.
		\end{equation}
		For $h\in H_0^1(\Omega)$, the
		Gagliardo--Nirenberg inequality applied to its zero
		extension, together with interpolation, yields
		\begin{equation}\label{eq:exterior-GN-Nash}
			\begin{aligned}
				\|h\|_{L^4(\Omega)}^2
				\leq K\|h\|_{L^2(\Omega)}
				\|\nabla h\|_{L^2(\Omega)},\qquad
				\|h\|_{L^2(\Omega)}^4
				\leq K\|h\|_{L^1(\Omega)}^2
				\|\nabla h\|_{L^2(\Omega)}^2.
			\end{aligned}
		\end{equation}
		The second inequality applies when $h\in L^1(\Omega)$.
		Indeed, substituting the first inequality into
		$
		\|h\|_{L^2(\Omega)}
		\leq
		\|h\|_{L^1(\Omega)}^{\frac13}
		\|h\|_{L^4(\Omega)}^{\frac23}
		$
		gives the second.
		
		To estimate the density, fix $p\geq2$ and $M>0$,
		and define
		\[
		n_M=\min\{n,M\},
		\qquad
		\Phi_{p,M}(a)
		=p\int_0^a\min\{\sigma,M\}^{p-1}d\sigma.
		\]
		For $a\geq0$, a direct calculation gives
		\[
		\begin{gathered}
			0\leq\Phi_{p,M}(a)\leq a^p,
			\qquad
			\Phi_{p,M}(a)\leq\frac p2 M^{p-2}a^2,\qquad
			\Phi_{p,M}(a)\to a^p
			\quad\text{as}\quad M\to\infty.
		\end{gathered}
		\]
		For fixed $p$ and $M$, the functions
		$a\mapsto\min\{a,M\}^{p-1}$ and
		$a\mapsto\min\{a,M\}^{\frac{p}{2}}$ are Lipschitz on
		$[0,\infty)$ and vanish at zero.
		Since $n(t)\in H_0^1(\Omega)$, the Sobolev chain rule gives
		$pn_M^{p-1}\in H_0^1(\Omega)$, and $n_M^{\frac p2}\in H_0^1(\Omega)$.
		Moreover,
		\[
		\nabla\Phi_{p,M}(n)
		=pn_M^{p-1}\nabla n,
		\qquad
		|\nabla\Phi_{p,M}(n)|
		\leq pM^{p-2}|n||\nabla n|.
		\]
		The preceding bounds and the Cauchy--Schwarz inequality yield
		\[
		\begin{aligned}
			\|\Phi_{p,M}(n)\|_{L^1(\Omega)}
			\leq\frac p2M^{p-2}\|n\|_{L^2(\Omega)}^2,\quad
			\|\nabla\Phi_{p,M}(n)\|_{L^1(\Omega)}
			\leq pM^{p-2}\|n\|_{L^2(\Omega)}
			\|\nabla n\|_{L^2(\Omega)}.
		\end{aligned}
		\]
		Thus $\Phi_{p,M}(n)\in W^{1,1}(\Omega)$.
		Its trace is zero because $n$ has zero trace and
		$\Phi_{p,M}(0)=0$, so
		$\Phi_{p,M}(n)\in W_0^{1,1}(\Omega)$.
		Using $\nabla\cdot V=0$, the boundedness of $V$,
		and approximation in $W_0^{1,1}(\Omega)$ by functions
		in $C_c^\infty(\Omega)$, we obtain
		\[
		p\int_\Omega n_M^{p-1}V\cdot\nabla n dx
		=\int_\Omega V\cdot\nabla\Phi_{p,M}(n)dx
		=0
		\]
		for almost every $t\in(0,\tau)$.
		The assumed strong regularity gives
		$n\in H^1(0,\tau;L^2(\Omega))$.
		Since $\Phi_{p,M}'$ is Lipschitz and
		$\Phi_{p,M}'(0)=0$, the time chain rule implies that
		$t\mapsto\int_\Omega\Phi_{p,M}(n(t))dx$
		is absolutely continuous, with
		\[
		\frac{d}{dt}\int_\Omega\Phi_{p,M}(n)dx
		=p\int_\Omega n_M^{p-1}\partial_tn dx
		\]
		for almost every $t\in(0,\tau)$.
		Finally, the spatial chain rule gives
		\[
		\nabla n_M=\mathbf1_{\{n<M\}}\nabla n,
		\qquad
		\nabla(n_M^{\frac p2})
		=\frac p2n_M^{\frac p2-1}\nabla n_M.
		\]
		Testing the density equation in
		\eqref{eq:exterior-moser-system} with $pn_M^{p-1}$
		and using these identities therefore yields
		\begin{equation}\label{eq:exterior-moser-truncated-energy}
			\begin{aligned}
				\frac{d}{dt}\int_\Omega\Phi_{p,M}(n)dx
				+\frac{4(p-1)}p
				\|\nabla(n_M^{\frac p2})\|_{L^2(\Omega)}^2=
				2(p-1)\int_\Omega
				n_M^{\frac p2}\nabla c\cdot
				\nabla(n_M^{\frac p2})dx
			\end{aligned}
		\end{equation}
		for almost every $t\in(0,\tau)$.
		By H\"older's inequality,
		\eqref{eq:exterior-elliptic-L4},
		\eqref{eq:exterior-GN-Nash}, and Young's inequality, there holds
		\begin{align*}
			&2(p-1)\left|
			\int_\Omega
			n_M^{\frac{p}{2}}\nabla c\cdot\nabla(n_M^{\frac{p}{2}})dx
			\right|\\
			\leq~&
			2(p-1)\|n_M^{\frac{p}{2}}\|_{L^4(\Omega)}
			\|\nabla c\|_{L^4(\Omega)}
			\|\nabla(n_M^{\frac{p}{2}})\|_{L^2(\Omega)}\leq
			KpB_1\|n_M^{\frac{p}{2}}\|_{L^2(\Omega)}^{\frac12}
			\|\nabla(n_M^{\frac{p}{2}})\|_{L^2(\Omega)}^{\frac32}\\
			\leq~&
			\frac12\|\nabla(n_M^{\frac{p}{2}})\|_{L^2(\Omega)}^2
			+Kp^4B_1^4\|n_M^{\frac{p}{2}}\|_{L^2(\Omega)}^2\leq
			\frac12\|\nabla(n_M^{\frac{p}{2}})\|_{L^2(\Omega)}^2
			+Kp^4B_1^4\|n\|_{L^p(\Omega)}^p.
		\end{align*}
		Noting that $\frac{4(p-1)}{p}\geq2$, integration of
		\eqref{eq:exterior-moser-truncated-energy} yields
		\begin{equation}\label{Phi ineq}
			\begin{aligned}
				&\int_\Omega\Phi_{p,M}(n(t))dx
				+\int_a^t
				\|\nabla(n_M^{\frac{p}{2}})(\sigma)\|_{L^2(\Omega)}^2 d\sigma\\\leq~&
				\int_\Omega\Phi_{p,M}(n(a))dx
				+Kp^4B_1^4\int_a^t
				\|n(\sigma)\|_{L^p(\Omega)}^p d\sigma
			\end{aligned}
		\end{equation}
		for $0\leq a<t\leq\tau$.
		Since $n\in C([0,\tau];L^p(\Omega))$, the preceding
		inequality with $a=0$ and $t=\tau$, together with
		$n_M^p\leq n^p$ and
		$\Phi_{p,M}(n_{\mathrm{in}})\leq n_{\mathrm{in}}^p$,
		gives
		\[
		\begin{aligned}
			\sup_{M>0}
			\|n_M^{\frac{p}{2}}\|_{L^2(0,\tau;H^1(\Omega))}^2\leq
			\|n_{\mathrm{in}}\|_{L^p(\Omega)}^p
			+(1+Kp^4B_1^4)
			\int_0^\tau\|n(\sigma)\|_{L^p(\Omega)}^p d\sigma
			<\infty.
		\end{aligned}
		\]
		Moreover,
		$
		n_M^{\frac{p}{2}}\to n^{\frac{p}{2}}$ almost everywhere and
		$|n_M^{\frac{p}{2}}-n^{\frac{p}{2}}|^2
		\leq n^p\in L^1(0,\tau;L^1(\Omega))$.
		The dominated convergence theorem yields
		\[
		n_M^{\frac{p}{2}}\to n^{\frac{p}{2}}
		\qquad\text{in}\quad L^2(0,\tau;L^2(\Omega))
		\quad\text{as}\quad M\to\infty.
		\]
		Choose an increasing sequence $M_j\to\infty$.
		The uniform bound in $L^2(0,\tau;H_0^1(\Omega))$
		gives a weakly convergent subsequence, still indexed by $j$.
		The preceding strong convergence identifies its weak limit
		with $n^{\frac{p}{2}}$.
		Hence $n^{\frac{p}{2}}\in L^2(0,\tau;H_0^1(\Omega))$ and
		$n_{M_j}^{\frac{p}{2}}\rightharpoonup n^{\frac{p}{2}}$
		in $L^2(0,\tau;H_0^1(\Omega))$.
		For every $0\leq a<t\leq\tau$, weak lower semicontinuity gives
		\[
		\int_a^t
		\|\nabla(n^{\frac{p}{2}})(\sigma)\|_{L^2(\Omega)}^2 d\sigma
		\leq
		\liminf_{j\to\infty}
		\int_a^t
		\|\nabla(n_{M_j}^{\frac{p}{2}})(\sigma)\|_{L^2(\Omega)}^2 d\sigma.
		\]
		At the time endpoints, $\Phi_{p,M_j}(n(s))$
		increases to $n(s)^p$.
		The monotone convergence theorem therefore yields
		\[
		\int_\Omega\Phi_{p,M_j}(n(s))dx
		\to
		\|n(s)\|_{L^p(\Omega)}^p
		\quad\text{for}\quad s=a,t.
		\]
		Passing to the limit in \eqref{Phi ineq}, we obtain
		\begin{equation}\label{eq:exterior-moser-energy}
			\begin{aligned}
				&\|n(t)\|_{L^p(\Omega)}^p
				+\int_a^t
				\|\nabla(n^{\frac{p}{2}})(\sigma)\|_{L^2(\Omega)}^2 d\sigma\\
				\leq~&
				\|n(a)\|_{L^p(\Omega)}^p
				+Kp^4(1+B_1^4)\int_a^t
				\|n(\sigma)\|_{L^p(\Omega)}^p d\sigma.
			\end{aligned}
		\end{equation}
		We now take $p=2^j$, $j\geq2$, and fix $K\geq1$,
		depending only on $\Omega$, large enough for both
		\eqref{eq:exterior-GN-Nash} and
		\eqref{eq:exterior-moser-energy}.
		If
		$\sup_{0\leq t\leq\tau}
		\|n(t)\|_{L^{\frac{p}{2}}(\Omega)}=0$,
		then $n=0$ on $[0,\tau]$.
		Otherwise, the second inequality in
		\eqref{eq:exterior-GN-Nash}, applied to $n^{\frac{p}{2}}$,
		gives
		\[
		\|\nabla(n^{\frac{p}{2}})(t)\|_{L^2(\Omega)}^2
		\geq
		\frac{\|n(t)\|_{L^p(\Omega)}^{2p}}
		{K\left(
			\sup_{0\leq s\leq\tau}
			\|n(s)\|_{L^{\frac{p}{2}}(\Omega)}
			\right)^p}
		\]
		for almost every $t\in(0,\tau)$.
		Substitution into \eqref{eq:exterior-moser-energy} yields
		\begin{equation}\label{eq:exterior-moser-comparison}
			\begin{aligned}
				&\|n(t)\|_{L^p(\Omega)}^p
				-\|n(a)\|_{L^p(\Omega)}^p\\
				\leq~&
				\int_a^t\|n(\sigma)\|_{L^p(\Omega)}^p
				\left[
				Kp^4(1+B_1^4)
				-\frac{\|n(\sigma)\|_{L^p(\Omega)}^p}
				{K\left(
					\sup_{0\leq s\leq\tau}
					\|n(s)\|_{L^{\frac{p}{2}}(\Omega)}
					\right)^p}
				\right]d\sigma .
			\end{aligned}
		\end{equation}
		The integrand is nonpositive whenever
		\[
		\|n(\sigma)\|_{L^p(\Omega)}^p
		\geq
		K^2p^4(1+B_1^4)
		\left(
		\sup_{0\leq s\leq\tau}
		\|n(s)\|_{L^{\frac{p}{2}}(\Omega)}
		\right)^p.
		\]
		Thus continuity and
		\eqref{eq:exterior-moser-comparison} give
		\begin{equation}\label{eq:exterior-moser-recursion}
			\begin{aligned}
				\sup_{0\leq t\leq\tau}
				\|n(t)\|_{L^p(\Omega)}\leq
				\max\left\{
				\|n_{\mathrm{in}}\|_{L^p(\Omega)},
				[K^2p^4(1+B_1^4)]^{\frac{1}{p}}
				\sup_{0\leq t\leq\tau}
				\|n(t)\|_{L^{\frac{p}{2}}(\Omega)}
				\right\}.
			\end{aligned}
		\end{equation}
		To prove \eqref{eq:exterior-moser-recursion}, suppose
		that it fails. By continuity, there exist
		$0\leq a<b\leq\tau$ such that
		\[
		\begin{aligned}
			\|n(a)\|_{L^p(\Omega)}^p
			&=\max\left\{
			\|n_{\mathrm{in}}\|_{L^p(\Omega)}^p,\,
			K^2p^4(1+B_1^4)
			\left(
			\sup_{0\leq s\leq\tau}
			\|n(s)\|_{L^{\frac{p}{2}}(\Omega)}
			\right)^p
			\right\},\\
			\|n(t)\|_{L^p(\Omega)}^p
			&>\|n(a)\|_{L^p(\Omega)}^p,
			\qquad a<t\leq b.
		\end{aligned}
		\]
		Consequently, for $a<\sigma\leq b$,
		\[
		Kp^4(1+B_1^4)
		-\frac{\|n(\sigma)\|_{L^p(\Omega)}^p}
		{K\left(
			\sup_{0\leq s\leq\tau}
			\|n(s)\|_{L^{\frac{p}{2}}(\Omega)}
			\right)^p}
		\leq0.
		\]
		Applying \eqref{eq:exterior-moser-comparison} on $[a,b]$
		therefore gives
		$
		0<
		\|n(b)\|_{L^p(\Omega)}^p
		-\|n(a)\|_{L^p(\Omega)}^p
		\leq0$,
		a contradiction. This proves
		\eqref{eq:exterior-moser-recursion}.
		By interpolation and the assumed $L^2(\Omega)$ bound, there holds
		\[
		\begin{aligned}
			\|n_{\mathrm{in}}\|_{L^p(\Omega)}
			\leq
			\|n_{\mathrm{in}}\|_{L^\infty(\Omega)}^{1-\frac{2}{p}}
			\|n_{\mathrm{in}}\|_{L^2(\Omega)}^{\frac{2}{p}}\leq
			\|n_{\mathrm{in}}\|_{L^\infty(\Omega)}+B_1.
		\end{aligned}
		\]
		Iterating \eqref{eq:exterior-moser-recursion} from
		the given $L^2(\Omega)$ bound, we obtain
		\[
		\begin{aligned}
			\sup_{0\leq t\leq\tau}
			\|n(t)\|_{L^{2^j}(\Omega)}
			&\leq
			\bigl(\|n_{\mathrm{in}}\|_{L^\infty(\Omega)}+B_1\bigr)
			\prod_{\ell=2}^j
			[K^2 2^{4\ell}(1+B_1^4)]^{2^{-\ell}}\\
			&\leq
			64K
			\bigl(\|n_{\mathrm{in}}\|_{L^\infty(\Omega)}+B_1\bigr)
			(1+B_1^4)^{\frac12}.
		\end{aligned}
		\]
		The last inequality follows from
		$\sum_{\ell=2}^\infty2^{-\ell}=\frac12$ and
		$\sum_{\ell=2}^\infty\ell\,2^{-\ell}=\frac32$.
		Letting $j\to\infty$, using
		$(1+B_1^4)^{\frac12}\leq1+B_1^2$, we conclude that
		\[
		\sup_{0\leq t\leq\tau}
		\|n(t)\|_{L^\infty(\Omega)}
		\leq
		K\bigl(\|n_{\mathrm{in}}\|_{L^\infty(\Omega)}+B_1\bigr)
		(1+B_1^2).
		\]
		Letting $\tau\to T$ completes the proof.
	\end{proof}
	\subsection{Global boundedness and decay of the radial density}
	We next establish global bounds and decay for the density $N$
	in the radial problem \eqref{eq:radial-reference-problem}.
	The bounds follow from Lemma \ref{lem:exterior-radial-screened} and radial interpolation, while the proof of decay uses the mass loss at the boundary $r=1$.
	\begin{proposition}[Radial bounds and decay]
		\label{prop:radial-reference-decay}
		Assume that $n_{\mathrm{in}}\geq0$ and
		$n_{\mathrm{in}}\in L^1(\Omega)\cap L^\infty(\Omega)\cap H_0^1(\Omega)$.
		Then the radial problem
		\eqref{eq:radial-reference-problem} has a unique global nonnegative
		strong solution $(N,C)$. Its density is uniformly bounded and satisfies
		\[
		\lim_{t\to\infty}\|N(t)\|_{L^\infty(\Omega)}=0,
		\qquad \lim_{t\to\infty}\|N(t)\|_{L^2(\Omega)}=0.
		\]
	\end{proposition}
	\begin{proof}
		The case $N(0)=0$ is
		immediate, so assume that
		$\int_1^\infty rN(0,r) dr>0$.
		
		\noindent{\bf Step 1. Global existence and uniform bounds.}
		Recall that $$N(0,r)=\frac1{2\pi}\int_0^{2\pi}n_{\mathrm{in}}(r,\theta) d\theta$$
		by the initial condition in
		\eqref{eq:radial-reference-problem}.
		Jensen's inequality and integration in polar coordinates give
		\[
		\|N(0)\|_{L^p(\Omega)}
		\leq\|n_{\mathrm{in}}\|_{L^p(\Omega)},
		\qquad p\in\{1,2\},
		\]
		while
		$\|N(0)\|_{L^\infty(\Omega)}
		\leq\|n_{\mathrm{in}}\|_{L^\infty(\Omega)}$.
		Similarly, we obtain
		\[
		\begin{aligned}
			\|\nabla N(0)\|_{L^2(\Omega)}^2
			&=2\pi\int_1^\infty r
			\left|
			\frac1{2\pi}\int_0^{2\pi}
			\partial_r n_{\mathrm{in}}(r,\theta)d\theta
			\right|^2 dr\\&\leq\int_1^\infty\int_0^{2\pi}
			r|\partial_r n_{\mathrm{in}}|^2 d\theta dr
			\leq\|\nabla n_{\mathrm{in}}\|_{L^2(\Omega)}^2.
		\end{aligned}
		\]
		Since $n_{\mathrm{in}}\in H_0^1(\Omega)$, choose
		$n_{\mathrm{in},j}\in C_c^\infty(\Omega)$ such that
		$n_{\mathrm{in},j}\to n_{\mathrm{in}}$ in $H^1(\Omega)$.
		Their angular averages are smooth and vanish near $r=1$
		and for sufficiently large $r$, so
		$(n_{\mathrm{in},j})_0\in C_c^\infty(\Omega)$.
		The $L^2(\Omega)$ and gradient estimates for angular
		averaging give
		\[
		\|(n_{\mathrm{in},j})_0-N(0)\|_{H^1(\Omega)}
		\leq\|n_{\mathrm{in},j}-n_{\mathrm{in}}\|_{H^1(\Omega)}
		\to 0.
		\]
		Hence $N(0)\in H_0^1(\Omega)$.
		Together with the preceding $L^1$ and $L^\infty$
		bounds implies
		$N(0)\in L^1(\Omega)\cap L^\infty(\Omega)\cap H_0^1(\Omega)$.
		Proposition~\ref{prop:local-strong-theory} therefore gives
		a unique local strong solution of
		\eqref{eq:radial-reference-problem}, with
		\begin{equation}\label{eq:radial-reference-initial-mass-bound}
			N\geq0,\qquad
			\int_1^\infty rN(t,r)dr
			\leq\int_1^\infty rN(0,r)dr.
		\end{equation}
		Indeed, radial initial data and uniqueness preserve radial
		symmetry. Since $\partial_\theta N=\partial_\theta C=0$ and
		$N e_r=\nabla\left(\int_1^r N(t,\rho)d\rho\right)$,
		$u=0$ satisfies the fluid equation with a suitable pressure,
		and the scalar equations reduce to
		\eqref{eq:radial-reference-problem}. The calculations below are first justified
		for smooth compatible data. Approximation by the Dirichlet
		heat semigroup and Proposition~\ref{prop:local-strong-theory} then
		give the same estimates for the stated initial data.
		
		Testing the density equation in
		\eqref{eq:radial-reference-problem} with $N$ gives
		\[
		\begin{aligned}
			&\frac12\frac{d}{dt}\int_1^\infty rN^2 dr
			+\int_1^\infty r|N_r|^2 dr
			\\=~&\int_1^\infty rNN_rC_r dr\leq\frac12\int_1^\infty r|N_r|^2dr
			+\frac12\|C_r\|_{L^{\infty}}^2\int_1^\infty rN^2 dr.
		\end{aligned}
		\]
		By \eqref{eq:radial-elliptic-estimates} with $f=N$, we get
		$\|C_r\|_{L^{\infty}}^2\leq6\int_1^\infty rN^2 dr$. Hence
		\begin{equation}\label{eq:radial-reference-y-energy}
			\frac{d}{dt}\int_1^\infty rN^2 dr
			+\int_1^\infty r|N_r|^2 dr
			\leq12\left(\int_1^\infty rN^2 dr\right)^2.
		\end{equation}
		For each positive time, the fundamental theorem of calculus
		and the Cauchy--Schwarz inequality yield
		\[
		N(t,r)^2
		=-2\int_r^\infty N(t,\rho)N_r(t,\rho)d\rho
		\leq2\left(\int_1^\infty \rho N^2 d\rho\right)^{\frac12}
		\left(\int_1^\infty \rho|N_r|^2 d\rho\right)^{\frac12}.
		\]
		Together with
		\eqref{eq:radial-reference-initial-mass-bound}, it gives
		\begin{align*}
			\|N\|_{L^{\infty}}^4
			\leq4\left(\int_1^\infty rN^2 dr\right)
			\left(\int_1^\infty r|N_r|^2\,dr\right),~
			\int_1^\infty rN^2 dr
			\leq\left(\int_1^\infty rN(0,r) dr\right)\|N\|_{L^{\infty}}.
		\end{align*}
		Therefore,
		\[
		\left(\int_1^\infty rN^2 dr\right)^3
		\leq4\left(\int_1^\infty rN(0,r)dr\right)^4
		\int_1^\infty r|N_r|^2 dr.
		\]
		Substituting this into \eqref{eq:radial-reference-y-energy}, we get
		\begin{equation*}
			\frac{d}{dt}\int_1^\infty rN^2 dr
			+\frac{\left(\int_1^\infty rN^2 dr\right)^3}
			{4\left(\int_1^\infty rN(0,r)dr\right)^4}
			\leq12\left(\int_1^\infty rN^2 dr\right)^2.
		\end{equation*}
		This implies
		\[
		\frac{d}{dt}\int_1^\infty rN^2 dr\leq0
		\quad\text{whenever}\quad
		\int_1^\infty rN^2 dr
		\geq48\left(\int_1^\infty rN(0,r) dr\right)^4.
		\]
		Thus the integral cannot exceed the larger of its initial
		value and the preceding threshold. Consequently, for every
		$t$ in the interval of existence, there holds
		\begin{equation}\label{eq:radial-reference-L2-bound}
			\int_1^\infty rN(t,r)^2 dr
			\leq\max\left\{\int_1^\infty rN(0,r)^2 dr,~
			48\left(\int_1^\infty rN(0,r) dr\right)^4\right\}.
		\end{equation}
		Since $\|N(t)\|_{L^2(\Omega)}^2
		=2\pi\int_1^\infty rN(t,r)^2 dr$,
		this estimate bounds $\|N(t)\|_{L^2(\Omega)}$
		uniformly throughout the interval of existence.
		Proposition~\ref{prop:density-continuation}
		therefore gives global existence.
		Applying Lemma \ref{lem:exterior-moser} and using
		\eqref{eq:radial-reference-L2-bound}, we obtain
		\begin{equation}\label{eq:radial-reference-density-Linfty}
			\sup_{t\geq0}\|N(t)\|_{L^\infty}\leq K.
		\end{equation}
		We next estimate $C$ using the elliptic equation in
		\eqref{eq:radial-reference-problem}.
		For each fixed $t$, let $C_-=\max\{-C,0\}$.
		Since $C\in H_0^1(\Omega)$, both $C_-$ and
		$(C-\|N\|_{L^{\infty}})_+$ belong to $H_0^1(\Omega)$.
		Testing the elliptic equation with $-C_-$ and using
		$N\geq0$ gives
		\[
		\int_1^\infty r
		\bigl(|\partial_r C_-|^2+C_-^2\bigr)dr
		=-\int_1^\infty rNC_- dr\leq0.
		\]
		Hence $C_-=0$, which proves $C\geq0$.
		To obtain the upper bound, test the elliptic equation with
		$(C-\|N\|_{L^{\infty}})_+$.
		Using
		$C=\|N\|_{L^{\infty}}+(C-\|N\|_{L^{\infty}})_+$
		on the set where $C>\|N\|_{L^{\infty}}$ and $N\leq \|N\|_{L^{\infty}}$, we obtain
		\[
		\begin{aligned}
			&\int_1^\infty r
			\left(
			\left|\partial_r(C-\|N\|_{L^{\infty}})_+\right|^2
			+\left|(C-\|N\|_{L^{\infty}})_+\right|^2
			\right)dr\\=~&
			\int_1^\infty r(N-\|N\|_{L^{\infty}})
			(C-\|N\|_{L^{\infty}})_+ dr
			\leq0.
		\end{aligned}
		\]
		Thus $(C-\|N\|_{L^{\infty}})_+=0$, and consequently
		\begin{equation}\label{eq:radial-reference-chemical-pointwise}
			0\leq C(t,r)\leq\|N(t)\|_{L^{\infty}}.
		\end{equation}
		Moreover, \eqref{eq:radial-elliptic-estimates} and
		\eqref{eq:radial-reference-L2-bound} imply
		\[
		\|C_r(t)\|_{L^{\infty}}^2
		\leq6\int_1^\infty rN(t,r)^2 dr
		\leq K.
		\]
		Combining this estimate with
		\eqref{eq:radial-reference-density-Linfty} and
		\eqref{eq:radial-reference-chemical-pointwise} gives
		\begin{equation}\label{eq:radial-reference-uniform-bounds}
			\sup_{t\geq0}\bigl(\|N(t)\|_{L^{\infty}}
			+\|C(t)\|_{L^{\infty}}+\|C_r(t)\|_{L^{\infty}}\bigr)\leq K.
		\end{equation}
		
		\noindent{\bf Step 2. Uniform bounds at large radii.}
		Integrating \eqref{eq:radial-reference-y-energy} over unit
		intervals and using \eqref{eq:radial-reference-L2-bound},
		we obtain
		\begin{equation}\label{eq:radial-reference-D-window}
			\sup_{a\geq0}\int_a^{a+1}\int_1^\infty
			r|N_r(t,r)|^2 drdt\leq K.
		\end{equation}
		Using the elliptic equation in
		\eqref{eq:radial-reference-problem}, we rewrite
		the density equation in the same system as
		\[
		N_t=N_{rr}+r^{-1}N_r-C_rN_r+(N-C)N.
		\]
		Testing it with $-(N_{rr}+r^{-1}N_r)$ and using Young's
		inequality gives
		\begin{align*}
			&\frac12\frac{d}{dt}\int_1^\infty r|N_r|^2 dr
			+\int_1^\infty r\bigl|N_{rr}+r^{-1}N_r\bigr|^2 dr\\
			\leq~&\frac12\int_1^\infty
			r\bigl|N_{rr}+r^{-1}N_r\bigr|^2 dr
			+K\|C_r\|_{L^{\infty}}^2\int_1^\infty r|N_r|^2 dr+K\|N-C\|_{L^{\infty}}^2\int_1^\infty rN^2 dr.
		\end{align*}
		Here the boundary term is zero because $N_t(t,1)=0$.
		By \eqref{eq:radial-reference-L2-bound} and
		\eqref{eq:radial-reference-uniform-bounds}, we arrive at
		\begin{equation}\label{eq:radial-reference-D-differential}
			\frac{d}{dt}\int_1^\infty r|N_r|^2 dr
			\leq K\left(1+\int_1^\infty r|N_r|^2 dr\right).
		\end{equation}
		For $t\geq1$, by
		\eqref{eq:radial-reference-D-window} we may choose $a\in[t-1,t]$ such that
		$\int_1^\infty r|N_r(a,r)|^2 dr\leq K$.
		Gronwall's inequality applied to
		\eqref{eq:radial-reference-D-differential} on $[a,t]$ gives
		\begin{equation}\label{eq:radial-reference-gradient-bound}
			\sup_{t\geq1}\int_1^\infty r|N_r(t,r)|^2 dr\leq K.
		\end{equation}
		For each fixed $t\geq1$, the integral mean value theorem
		gives $R_j\in[j,j+1]$ such that
		\[
		R_jN(t,R_j)^2
		=\int_j^{j+1}\rho N(t,\rho)^2 d\rho
		\to 0\quad\text{as}\quad j\to\infty,
		\]
		where the convergence follows from
		\eqref{eq:radial-reference-L2-bound}.
		Integrating
		$\partial_\rho(\rho N^2)=N^2+2\rho NN_r$
		over $[r,R_j]$, we obtain
		\[
		rN(t,r)^2+\int_r^{R_j}N(t,\rho)^2 d\rho
		=R_jN(t,R_j)^2
		-2\int_r^{R_j}\rho N(t,\rho)N_r(t,\rho)d\rho.
		\]
		By \eqref{eq:radial-reference-L2-bound},
		\eqref{eq:radial-reference-gradient-bound} and the
		Cauchy--Schwarz inequality, we find that
		$\int_1^\infty\rho|NN_r| d\rho<\infty$.
		Letting $j\to\infty$ therefore gives
		\begin{align*}
			&rN(t,r)^2+\int_r^\infty N(t,\rho)^2 d\rho
			=-2\int_r^\infty\rho N(t,\rho)N_r(t,\rho) d\rho\\
			&\leq2\left(\int_1^\infty
			\rho N(t,\rho)^2 d\rho\right)^{\frac12}
			\left(\int_1^\infty
			\rho|N_r(t,\rho)|^2 d\rho\right)^{\frac12}\leq K.
		\end{align*}
		Consequently,
		\begin{equation}\label{eq:radial-reference-tail}
			\sup_{t\geq1}\sup_{r\geq R}N(t,r)
			\leq KR^{-\frac12},\qquad R>1.
		\end{equation}
		
		\noindent{\bf Step 3. Green function representation and mass balance.}
		To integrate the density equation over $(1,\infty)$, we first
		estimate $C_r$ at infinity. Set
		\[
		K_0(r)=\int_0^\infty e^{-r\cosh z}dz,
		\qquad
		h_0(r)=K_0(r)\int_1^r\frac{d\rho}{\rho K_0(\rho)^2}.
		\]
		Differentiation gives
		\[
		r(K_0''-K_0)+K_0'
		=-\int_0^\infty\partial_z
		\bigl(\sinh z e^{-r\cosh z}\bigr)dz=0,
		\quad h_0''+r^{-1}h_0'-h_0=0,
		\]
		and
		\begin{equation}\label{eq:radial-reference-kernel-identities}
			K_0>0,\qquad K_0'<0,\qquad h_0(1)=0,
			\qquad r(h_0'K_0-h_0K_0')=1.
		\end{equation}
		Moreover,
		\[
		h_0(r)>0\quad\text{for}\quad r>1,\qquad
		h_0'(1)=\frac1{K_0(1)}>0,\qquad
		(rh_0')'=rh_0>0,
		\]
		hence $h_0'>0$. The inequality $\cosh z\geq1+\frac{z^2}{2}$ yields
		\begin{equation}\label{eq:radial-reference-kernel-decay}
			K_0(r)\leq e^{-r}\int_0^\infty e^{-\frac{rz^2}{2}}dz
			=\left(\frac{\pi}{2r}\right)^{\frac12}e^{-r}.
		\end{equation}
		For the elliptic equation for $C$ in
		\eqref{eq:radial-reference-problem}, with $C(t,1)=0$ and
		$C(t)\in H_0^1(\Omega)$, the Green function and the
		representation of $C$ are
		\begin{equation}\label{eq:radial-reference-chemical-representation}
			G_0(r,\rho)=h_0(\min\{r,\rho\})K_0(\max\{r,\rho\}),
			\quad
			C(t,r)=\int_1^\infty G_0(r,\rho)N(t,\rho)\rho d\rho.
		\end{equation}
		We first verify the Green function and then justify the
		representation for $N(t)$.
		For each fixed $\rho>1$, we have
		\begin{equation}\label{eq:radial-reference-green-piecewise}
			G_0(r,\rho)=
			\begin{cases}
				h_0(r)K_0(\rho),&\quad 1\leq r\leq\rho,\\
				h_0(\rho)K_0(r),&\quad r\geq\rho.
			\end{cases}
		\end{equation}
		Since $h_0$ and $K_0$ solve $f''+r^{-1}f'-f=0$,
		\eqref{eq:radial-reference-green-piecewise} gives
		\[
		-\frac1r\partial_r\bigl(r\partial_rG_0(r,\rho)\bigr)
		+G_0(r,\rho)=0,\qquad r\neq\rho.
		\]
		The two expressions in
		\eqref{eq:radial-reference-green-piecewise} agree at $r=\rho$.
		Moreover, \eqref{eq:radial-reference-kernel-identities} and
		\eqref{eq:radial-reference-kernel-decay} imply
		\[
		G_0(1,\rho)=h_0(1)K_0(\rho)=0,\qquad
		\lim_{r\to\infty}G_0(r,\rho)
		=h_0(\rho)\lim_{r\to\infty}K_0(r)=0.
		\]
		Differentiating \eqref{eq:radial-reference-green-piecewise}
		on either side of $r=\rho$, we obtain
		\begin{equation}\label{eq:radial-reference-green-jump}
			\begin{aligned}
				\partial_rG_0(\rho+,\rho)
				-\partial_rG_0(\rho-,\rho)=h_0(\rho)K_0'(\rho)-h_0'(\rho)K_0(\rho)
				=-\frac1\rho,
			\end{aligned}
		\end{equation}
		where the last equality follows from
		\eqref{eq:radial-reference-kernel-identities}.
		Thus $G_0$ satisfies the homogeneous elliptic equation away
		from $r=\rho$, the Dirichlet condition at $r=1$, decay at
		infinity, and the derivative jump
		\eqref{eq:radial-reference-green-jump}.
		Fix $t$. We first verify
		\eqref{eq:radial-reference-chemical-representation}
		for radial $L^2(\Omega)$ data with bounded support,
		and then remove the support restriction by truncation.
		
		\noindent\textbf{\underline{Case 1. Radial data with bounded support.}}
		Assume that $N(t,r)=0$ for almost every $r\geq S$, for some
		$S>1$. Splitting the Green integral
		\eqref{eq:radial-reference-green-piecewise}, we obtain
		\begin{equation}\label{eq:radial-reference-chemical-split}
			C(t,r)
			=K_0(r)\int_1^r h_0(\rho)N(t,\rho)\rho d\rho
			+h_0(r)\int_r^\infty K_0(\rho)N(t,\rho)\rho d\rho.
		\end{equation}
		The integrands are integrable on bounded radial intervals,
		so the following differentiation identities hold almost
		everywhere. Differentiating
		\eqref{eq:radial-reference-chemical-split}, the two terms
		arising from the variable integration limits cancel. Thus
		\begin{equation}\label{eq:radial-reference-chemical-derivative}
			C_r(t,r)
			=K_0'(r)\int_1^r h_0(\rho)N(t,\rho)\rho d\rho
			+h_0'(r)\int_r^\infty K_0(\rho)N(t,\rho)\rho d\rho.
		\end{equation}
		Differentiating once more and using the equations for
		$h_0$ and $K_0$ together with
		\eqref{eq:radial-reference-kernel-identities}, we obtain
		\begin{equation}\label{eq:radial-reference-chemical-equation-check}
			\begin{aligned}
				C_{rr}+r^{-1}C_r-C=rN(t,r)\bigl(h_0(r)K_0'(r)-h_0'(r)K_0(r)\bigr)=-N(t,r).
			\end{aligned}
		\end{equation}
		It remains to verify that the function in
		\eqref{eq:radial-reference-chemical-split}
		belongs to $H_0^1(\Omega)$.
		On $[1,S]$, the coefficient functions in
		\eqref{eq:radial-reference-chemical-split} and
		\eqref{eq:radial-reference-chemical-derivative}
		are bounded, and the integrals are bounded by the
		Cauchy--Schwarz inequality. Hence
		\begin{equation}\label{eq:radial-reference-chemical-local-energy}
			\int_1^S r
			\bigl(|C(t,r)|^2+|C_r(t,r)|^2\bigr)dr<\infty.
		\end{equation}
		For $r\geq S$, the vanishing of $N(t,r)$ reduces the two
		formulas to
		\begin{equation}\label{eq:radial-reference-chemical-outside-support}
			\begin{aligned}
				C(t,r)=K_0(r)\int_1^S h_0(\rho)N(t,\rho)\rho d\rho,\quad
				C_r(t,r)=K_0'(r)\int_1^S h_0(\rho)N(t,\rho)\rho d\rho.
			\end{aligned}
		\end{equation}
		Since $(rK_0')'=rK_0$, we have
		$\partial_r(rK_0K_0')=r\bigl(|K_0'|^2+|K_0|^2\bigr)$.
		Integrating over $[S,R]$ and using $K_0>0>K_0'$ from
		\eqref{eq:radial-reference-kernel-identities}, we find
		\begin{align*}
			\int_S^R r\bigl(|K_0'|^2+|K_0|^2\bigr)dr
			=RK_0(R)K_0'(R)-SK_0(S)K_0'(S)\leq -SK_0(S)K_0'(S).
		\end{align*}
		Letting $R\to\infty$ and using
		\eqref{eq:radial-reference-chemical-outside-support} gives
		\begin{align*}
			&\int_S^\infty r
			\bigl(|C(t,r)|^2+|C_r(t,r)|^2\bigr)dr=
			\left(\int_1^S h_0(\rho)N(t,\rho)\rho d\rho\right)^2
			\int_S^\infty r\bigl(|K_0|^2+|K_0'|^2\bigr)dr\\&\leq
			-SK_0(S)K_0'(S)
			\left(\int_1^S h_0(\rho)N(t,\rho)\rho d\rho\right)^2
			<\infty.
		\end{align*}
		Combining this with
		\eqref{eq:radial-reference-chemical-local-energy}, we obtain
		\[
		\|C(t)\|_{H^1(\Omega)}^2
		=2\pi\int_1^\infty r
		\bigl(|C(t,r)|^2+|C_r(t,r)|^2\bigr)dr<\infty.
		\]
		Moreover, \eqref{eq:radial-reference-chemical-split} and
		$h_0(1)=0$ give $C(t,1)=0$. Thus the function defined by
		\eqref{eq:radial-reference-chemical-split}
		belongs to $H_0^1(\Omega)$ and satisfies the elliptic
		equation by \eqref{eq:radial-reference-chemical-equation-check}.
		Uniqueness in Lemma~\ref{lem:exterior-radial-screened}
		therefore proves
		\eqref{eq:radial-reference-chemical-representation}
		for data with bounded support.
		
		\noindent\textbf{\underline{Case 2. General radial data.}}
		For $S>1$, truncate the right-hand side to
		$N(t,r)\mathbf1_{\{r<S\}}$. By Case 1, the corresponding
		Dirichlet solution is
		\[
		C^{(S)}(t,r)
		=\int_1^S G_0(r,\rho)N(t,\rho)\rho d\rho.
		\]
		The difference $C^{(S)}-C$ satisfies
		\[
		\begin{cases}
			-\dfrac{1}{r}\partial_r\!\left(r\partial_r(C^{(S)}-C)\right)
			+(C^{(S)}-C)
			=-N(t,r)\mathbf1_{\{r\geq S\}}, & \text{a.e.}\quad r>1,\\[4pt]
			(C^{(S)}-C)(t,1)=0.
		\end{cases}
		\]
		Applying \eqref{eq:radial-elliptic-estimates} to
		$C^{(S)}-C$ yields
		\begin{equation}\label{eq:radial-reference-chemical-truncation}
			\begin{aligned}
				&\int_1^\infty r\left(
				|C^{(S)}(t,r)-C(t,r)|^2
				+|\partial_r C^{(S)}(t,r)-C_r(t,r)|^2
				\right)dr\\
				&\qquad\leq 2\int_S^\infty rN(t,r)^2 dr
				\to 0
				\qquad\text{as}\quad S\to\infty.
			\end{aligned}
		\end{equation}
		Since $r\geq1$, this estimate also gives convergence in
		$H^1(1,R)$ for every $R>1$. By the one-dimensional Sobolev
		embedding, $C^{(S)}(t,\cdot)$ therefore converges to
		$C(t,\cdot)$ uniformly on $[1,R]$.
		On the other hand, noting that $h_0$ is nonnegative and increasing
		and $K_0$ is positive and decreasing, we arrive at
		$0\leq G_0(r,\rho)\leq h_0(r)K_0(r)$.
		By monotone convergence and
		\eqref{eq:radial-reference-initial-mass-bound},
		for each $r\geq1$ we have
		\[
		\begin{aligned}
			C^{(S)}(t,r)
			\to
			\int_1^\infty G_0(r,\rho)N(t,\rho)\rho d\rho\leq h_0(r)K_0(r)
			\int_1^\infty \rho N(0,\rho)d\rho<\infty.
		\end{aligned}
		\]
		Comparing this with the limit given by
		\eqref{eq:radial-reference-chemical-truncation}, we obtain
		\eqref{eq:radial-reference-chemical-representation}
		for every $r\geq1$. Splitting the resulting integral at
		$\rho=r$ and differentiating as in Case 1 give
		\eqref{eq:radial-reference-chemical-split} and
		\eqref{eq:radial-reference-chemical-derivative}
		for general $N(t)$.
		
		We next estimate $C_r$ at large radii.
		It follows from \eqref{eq:radial-reference-kernel-identities} that
		\begin{equation}\label{eq:radial-reference-kernel-diagonal-bounds}
			\begin{gathered}
				rh_0'(r)K_0(r)+rh_0(r)|K_0'(r)|=1,\quad
				rh_0'(r)K_0(r)\leq1,\quad
				rh_0(r)|K_0'(r)|\leq1.
			\end{gathered}
		\end{equation}
		Using \eqref{eq:radial-reference-green-piecewise}, the
		monotonicity of $h_0$ and $K_0$, and
		\eqref{eq:radial-reference-kernel-diagonal-bounds}, we obtain
		\begin{equation*}
			r|\partial_rG_0(r,\rho)|=
			\begin{cases}
				rh_0'(r)K_0(\rho)
				\leq rh_0'(r)K_0(r)\leq1,&~r<\rho,\\
				rh_0(\rho)|K_0'(r)|
				\leq rh_0(r)|K_0'(r)|\leq1,&~r>\rho.
			\end{cases}
		\end{equation*}
		Combining this
		with \eqref{eq:radial-reference-chemical-derivative}
		and \eqref{eq:radial-reference-initial-mass-bound} gives
		\begin{equation}\label{eq:radial-reference-chemical-gradient}
			\begin{gathered}
				r|\partial_rG_0(r,\rho)|\leq1\quad\text{for}\quad r\neq\rho,\\
				|C_r(t,r)|
				\leq\frac1r\int_1^\infty\rho N(t,\rho)d\rho
				\leq\frac1r\int_1^\infty\rho N(0,\rho)d\rho.
			\end{gathered}
		\end{equation}
		We now use \eqref{eq:radial-reference-chemical-gradient}
		to justify integration of the density equation over
		$(1,\infty)$.
		For $R\geq2$, choose $\zeta_R\in C^\infty([1,\infty))$
		such that
		\begin{equation}\label{eq:radial-reference-mass-cutoff}
			\begin{gathered}
				0\leq\zeta_R\leq1,\qquad
				\zeta_R=1\quad \text{on}\quad[1,R],\qquad
				\zeta_R=0\quad \text{on}\quad[2R,\infty),\\
				|\zeta_R'|\leq \frac{K}{R},\qquad
				|\zeta_R''|\leq \frac{K}{R^2}.
			\end{gathered}
		\end{equation}
		Testing the density equation in
		\eqref{eq:radial-reference-problem} with $\zeta_R$ and
		using $N(t,1)=0$, we obtain
		\begin{equation}\label{eq:radial-reference-cutoff-mass-identity}
			\begin{aligned}
				\frac{d}{dt}\int_1^\infty N\zeta_Rr dr
				&=-N_r(t,1)
				-\int_1^\infty rN_r\zeta_R' dr
				+\int_1^\infty rNC_r\zeta_R'dr\\
				&=-N_r(t,1)+\int_1^\infty
				N(\zeta_R''+r^{-1}\zeta_R'+C_r\zeta_R')r dr.
			\end{aligned}
		\end{equation}
		The derivatives of $\zeta_R$ are supported in $[R,2R]$.
		Hence \eqref{eq:radial-reference-initial-mass-bound}, \eqref{eq:radial-reference-chemical-gradient} and \eqref{eq:radial-reference-mass-cutoff} imply
		\begin{equation*}
			\begin{aligned}
				&\left|\int_1^\infty
				N(\zeta_R''+r^{-1}\zeta_R'+C_r\zeta_R')rdr\right|\\\leq~&\frac{K}{R^2}
				\left(1+\int_1^\infty\rho N(0,\rho)d\rho\right)
				\int_R^{2R}rN(t,r)dr
				\leq\frac{K}{R^2}.
			\end{aligned}
		\end{equation*}
		Integrating \eqref{eq:radial-reference-cutoff-mass-identity}
		from $\varepsilon$ to $t$, with $0<\varepsilon<t$, and
		using the above inequality, we find
		\begin{equation}\label{eq:radial-reference-cutoff-mass-integrated}
			\begin{aligned}
				\bigg|\int_1^\infty rN(t,r)\zeta_R(r)dr
				-\int_1^\infty rN(\varepsilon,r)\zeta_R(r)dr+\int_\varepsilon^tN_r(\sigma,1)d\sigma\bigg|
				\leq\frac{K(t-\varepsilon)}{R^2}.
			\end{aligned}
		\end{equation}
		Since $\zeta_R(r)\to1$ for every fixed $r$,
		\eqref{eq:radial-reference-initial-mass-bound} permits
		dominated convergence in
		\eqref{eq:radial-reference-cutoff-mass-integrated}.
		Letting $R\to\infty$, we obtain
		\begin{equation}\label{eq:radial-reference-positive-time-mass}
			\int_1^\infty rN(t,r)dr
			+\int_\varepsilon^tN_r(\sigma,1)d\sigma
			=\int_1^\infty rN(\varepsilon,r)dr.
		\end{equation}
		It remains to pass to the initial time in
		\eqref{eq:radial-reference-positive-time-mass}.
		By Proposition \ref{prop:local-strong-theory},
		$N(\varepsilon)\to N(0)$ in $L^2(\Omega)$ as
		$\varepsilon\to 0$.
		For each fixed $R>1$, the Cauchy--Schwarz inequality gives
		\begin{equation}\label{eq:radial-reference-local-mass-continuity}
			\begin{aligned}
				\left|\int_1^R r\bigl(N(\varepsilon,r)-N(0,r)\bigr)dr\right|
				\leq~&\left(\int_1^R r\,dr\right)^{\frac12}
				\left(\int_1^R
				r|N(\varepsilon,r)-N(0,r)|^2 dr\right)^{\frac12}\\
				&\to 0\quad\text{as}\quad\varepsilon\to 0.
			\end{aligned}
		\end{equation}
		Using \eqref{eq:radial-reference-initial-mass-bound} and \eqref{eq:radial-reference-local-mass-continuity}, we deduce
		\begin{align*}
			\int_1^R rN(0,r)dr
			\leq\liminf_{\varepsilon\to 0}
			\int_1^\infty rN(\varepsilon,r)dr\leq\limsup_{\varepsilon\to 0}
			\int_1^\infty rN(\varepsilon,r)dr
			\leq\int_1^\infty rN(0,r)dr.
		\end{align*}
		Letting $R\to\infty$ yields
		\begin{equation}\label{eq:radial-reference-initial-mass-continuity}
			\lim_{\varepsilon\to 0}
			\int_1^\infty rN(\varepsilon,r)dr
			=\int_1^\infty rN(0,r)dr.
		\end{equation}
		Since $N\geq0$ and $N(t,1)=0$, we have $N_r(t,1)\geq0$
		for $t>0$.
		Thus \eqref{eq:radial-reference-initial-mass-continuity}
		and monotone convergence allow $\varepsilon\to 0$
		in \eqref{eq:radial-reference-positive-time-mass}, giving
		\begin{equation}\label{eq:radial-reference-mass-balance}
			\int_1^\infty rN(t,r) dr
			+\int_0^tN_r(\sigma,1)d\sigma
			=\int_1^\infty rN(0,r)dr.
		\end{equation}
		Finally, nonnegativity of $N$ and
		\eqref{eq:radial-reference-mass-balance} imply
		$\int_0^tN_r(\sigma,1)d\sigma
		\leq\int_1^\infty rN(0,r)dr$.
		Letting $t\to\infty$ by monotone convergence, we conclude that
		\begin{equation}\label{eq:radial-reference-flux-integrable}
			\int_0^\infty N_r(t,1)dt
			\leq\int_1^\infty rN(0,r)dr.
		\end{equation}
		
		\noindent\textbf{Step 4. Compactness and the limit system.}
		Let $t_j\to\infty$ be arbitrary. We establish local
		compactness of the time translates of $N$ and $N_r$ up to
		$r=1$ and identify the equations satisfied by the limit.
		
		Using the elliptic equation in
		\eqref{eq:radial-reference-problem}, we rewrite the density
		equation as
		\[
		N_t-N_{rr}=(r^{-1}-C_r)N_r+(N-C)N.
		\]
		For $a\geq2$ and $R>1$, set
		\[
		Q_\ell
		=\left(a-\frac{\ell}{3},a+1+\frac{\ell}{3}\right)\times(1,R+\ell),
		\qquad \ell\in\{0,1,2,3\}.
		\]
		For each $\ell\in\{1,2,3\}$, choose a smooth cutoff $\chi$ equal
		to one on $Q_{\ell-1}$ and vanishing near the time
		endpoints and the outer radial boundary of $Q_\ell$.
		The derivatives of $\chi$ are bounded independently of
		$a$. The function $\chi N$ has zero initial value and
		satisfies homogeneous Dirichlet conditions at both radial
		endpoints. Moreover,
		\[
		(\chi N)_t-(\chi N)_{rr}
		=\bigl[\chi(r^{-1}-C_r)-2\chi_r\bigr]N_r
		+\bigl[\chi(N-C)+\chi_t-\chi_{rr}\bigr]N.
		\]
		The Dirichlet $L^p$ estimate for the heat equation
		\cite[Chapter~VII]{Lieberman1996} and
		\eqref{eq:radial-reference-uniform-bounds} therefore give
		\begin{equation}\label{eq:radial-reference-local-heat-estimate}
			\begin{aligned}
				\|N\|_{W_p^{2,1}(Q_{\ell-1})}
				&\leq K_{p,R}
				\|(\chi N)_t-(\chi N)_{rr}\|_{L^p(Q_\ell)}\\
				&\leq K_{p,R}
				\bigl(\|N_r\|_{L^p(Q_\ell)}
				+\|N\|_{L^p(Q_\ell)}\bigr),
				\qquad 1<p<\infty.
			\end{aligned}
		\end{equation}
		Here $W_p^{2,1}$ requires the function, its first time
		derivative, and its first two spatial derivatives to
		belong to $L^p$.
		Since $a-1\geq1$, the gradient bound
		\eqref{eq:radial-reference-gradient-bound} gives
		\[
		\|N_r\|_{L^2(Q_3)}^2
		\leq 3\sup_{s\geq1}\int_1^\infty r|N_r(s,r)|^2 dr
		\leq K.
		\]
		Taking $p=2$ and $\ell=3$ in
		\eqref{eq:radial-reference-local-heat-estimate}, we obtain
		\[
		\|N\|_{W_2^{2,1}(Q_2)}\leq K_R.
		\]
		To estimate $N_r$ in $L^4(Q_2)$, apply the one-dimensional
		embedding $H^1(1,R+2)\hookrightarrow L^\infty(1,R+2)$
		at almost every time. Together with
		\eqref{eq:radial-reference-gradient-bound}, this yields
		\begin{align*}
			\|N_r\|_{L^4(Q_2)}^4
			&\leq
			\int_{a-\frac23}^{a+\frac53}
			\|N_r(t)\|_{L^2(1,R+2)}^2
			\|N_r(t)\|_{L^\infty(1,R+2)}^2 dt\\&\leq K_R\int_{Q_2}
			\bigl(|N_r|^2+|N_{rr}|^2\bigr)drdt
			\leq K_R.
		\end{align*}
		Taking $p=4$ and $\ell=2$ in
		\eqref{eq:radial-reference-local-heat-estimate} and using
		parabolic Sobolev embedding in one spatial dimension,
		we find
		\[
		\|N_r\|_{L^\infty(Q_1)}
		\leq K_R\|N\|_{W_4^{2,1}(Q_1)}
		\leq K_R.
		\]
		In particular,
		\[
		\|N_r\|_{L^p(Q_1)}+\|N\|_{L^p(Q_1)}
		\leq |Q_1|^{\frac{1}{p}}
		\bigl(\|N_r\|_{L^\infty(Q_1)}
		+\|N\|_{L^\infty(Q_1)}\bigr)
		\leq K_{p,R}.
		\]
		Thus \eqref{eq:radial-reference-local-heat-estimate}
		with $\ell=1$ gives
		\begin{equation}\label{eq:radial-reference-parabolic-regularity}
			\sup_{a\geq2}
			\|N\|_{W_p^{2,1}((a,a+1)\times(1,R))}
			\leq K_{p,R},
			\qquad 1<p<\infty.
		\end{equation}
		Fix $p>3$. For any $L>0$ and $R>1$, covering
		$(t_j-L,t_j+L)$ by finitely many unit intervals and using
		\eqref{eq:radial-reference-parabolic-regularity}, we obtain
		\[
		\|N(t_j+\cdot,\cdot)\|_{W_p^{2,1}((-L,L)\times(1,R))}
		\leq K_{p,L,R}
		\]
		for all sufficiently large $j$. For
		$0<\alpha<1-\frac{3}{p}$, parabolic Sobolev embedding gives
		\begin{align*}
			&\|N(t_j+\cdot,\cdot)\|_{C^{\frac{\alpha}{2},\alpha}([-L,L]\times[1,R])}
			+\|N_r(t_j+\cdot,\cdot)\|_{C^{\frac{\alpha}{2},\alpha}([-L,L]\times[1,R])}\\
			&\qquad\leq
			K_{p,L,R,\alpha}
			\|N(t_j+\cdot,\cdot)\|_{W_p^{2,1}((-L,L)\times(1,R))}
			\leq K_{p,L,R,\alpha}.
		\end{align*}
		By the Arzel\`a--Ascoli theorem and a diagonal argument,
		there is a subsequence, still indexed by $j$, such that
		\begin{equation}\label{eq:radial-reference-time-translates}
			\begin{aligned}
				N(t_j+\sigma,r)&\to N_*(\sigma,r),\\
				N_r(t_j+\sigma,r)&\to\partial_rN_*(\sigma,r)
			\end{aligned}
			\qquad\text{locally uniformly on }
			\mathbb R\times[1,\infty).
		\end{equation}
		The same subsequence may be taken to converge weakly
		in $W_p^{2,1}$ on bounded rectangles. Its weak limit
		agrees with the locally uniform limit, so
		$N_*\in W_p^{2,1}$ on these rectangles.
		By \eqref{eq:radial-reference-time-translates},
		$N_*\geq0$ and $N_*(\sigma,1)=0$. Fatou's lemma and
		\eqref{eq:radial-reference-initial-mass-bound} also give
		\begin{equation}\label{eq:radial-reference-limit-mass-bound}
			\int_1^\infty \rho N_*(\sigma,\rho)d\rho
			\leq \int_1^\infty \rho N(0,\rho)d\rho.
		\end{equation}
		The uniform $L^\infty$ bound passes to $N_*$, and hence
		\[
		\int_1^\infty rN_*(\sigma,r)^2 dr
		\leq \|N_*(\sigma,\cdot)\|_{L^\infty(1,\infty)}
		\int_1^\infty rN_*(\sigma,r)dr
		\leq K.
		\]
		For each $\sigma$, let $C_*(\sigma)\in H_0^1(\Omega)$
		be the Dirichlet solution with right-hand side
		$N_*(\sigma)$ given by
		Lemma~\ref{lem:exterior-radial-screened}. Since
		$N_*(\sigma)\in L^1(\Omega)\cap L^2(\Omega)$, the Green
		representation and its differentiated form give
		\[
		\begin{aligned}
			C_*(\sigma,r)
			=\int_1^\infty G_0(r,\rho)N_*(\sigma,\rho)\rho d\rho,\quad
			\partial_rC_*(\sigma,r)
			=\int_1^\infty
			\partial_rG_0(r,\rho)N_*(\sigma,\rho)\rho d\rho.
		\end{aligned}
		\]
		Fix $L>0$ and $S>R>1$. By
		\eqref{eq:radial-reference-green-piecewise}, for $1\leq r\leq R,\ 1\leq\rho\leq S,\ \rho\neq r$, we have
		\[
		|G_0(r,\rho)|+|\partial_rG_0(r,\rho)|
		\leq K_{R,S}.
		\]
		For $\rho>S$, we have $r<\rho$, so the same formula gives
		\[
		|G_0(r,\rho)|+|\partial_rG_0(r,\rho)|
		=\bigl(|h_0(r)|+|h_0'(r)|\bigr)K_0(\rho)
		\leq K_RK_0(S).
		\]
		Subtracting the Green representations and splitting the
		difference at $\rho=S$, we obtain, for $|\sigma|\leq L$
		and $1\leq r\leq R$,
		\begin{align*}
			&|C(t_j+\sigma,r)-C_*(\sigma,r)|
			+|C_r(t_j+\sigma,r)-\partial_rC_*(\sigma,r)|\\
			\leq~&
			K_{R,S}\sup_{1\leq\rho\leq S}
			|N(t_j+\sigma,\rho)-N_*(\sigma,\rho)|
			+K_RK_0(S)\int_S^\infty
			\rho\bigl(N(t_j+\sigma,\rho)+N_*(\sigma,\rho)\bigr)d\rho\\
			\leq~&
			K_{R,S}\sup_{1\leq\rho\leq S}
			|N(t_j+\sigma,\rho)-N_*(\sigma,\rho)|
			+2K_RK_0(S)\int_1^\infty \rho N(0,\rho)d\rho,
		\end{align*}
		where the last inequality follows from
		\eqref{eq:radial-reference-initial-mass-bound} and
		\eqref{eq:radial-reference-limit-mass-bound}.
		For fixed $S$, the first term tends to zero uniformly
		in $|\sigma|\leq L$ by
		\eqref{eq:radial-reference-time-translates}. The second
		term tends to zero as $S\to\infty$ by
		\eqref{eq:radial-reference-kernel-decay}. Consequently,
		\begin{equation}\label{eq:radial-reference-chemical-translates}
			\begin{aligned}
				C(t_j+\sigma,r)&\to C_*(\sigma,r),\\
				C_r(t_j+\sigma,r)&\to\partial_rC_*(\sigma,r)
			\end{aligned}
			\qquad\text{locally uniformly on }
			\mathbb R\times[1,\infty).
		\end{equation}
		Let $\varphi\in C_c^\infty(\mathbb R\times(1,\infty))$.
		For sufficiently large $j$, the support of $\varphi$
		lies in $\{\sigma>-t_j\}$. Setting $t=t_j+\sigma$ in the density equation in
		\eqref{eq:radial-reference-problem}, we obtain
		\[
		\partial_\sigma N(t_j+\sigma,r)
		=\frac{1}{r}\partial_r\!\left[
		r(N_r-NC_r)(t_j+\sigma,r)
		\right],
		\qquad \sigma>-t_j,\quad r>1.
		\]
		Multiplying by $r\varphi$ and integrating by parts in
		$\sigma$ and $r$, we obtain
		\[
		\int_{-t_j}^{\infty}\int_1^\infty r
		\left[
		-N(t_j+\sigma,r)\partial_\sigma\varphi
		+(N_r-NC_r)(t_j+\sigma,r)\partial_r\varphi
		\right]drd\sigma=0.
		\]
		Passing to the limit by
		\eqref{eq:radial-reference-time-translates} and
		\eqref{eq:radial-reference-chemical-translates}, we obtain
		\[
		\int_{\mathbb R}\int_1^\infty r
		\left[
		-N_*\partial_\sigma\varphi
		+\bigl(\partial_rN_*-N_*\partial_rC_*\bigr)
		\partial_r\varphi
		\right]drd\sigma=0.
		\]
		Together with 
		\eqref{eq:radial-elliptic-problem} for
		$c=C_*(\sigma)$ and $f=N_*(\sigma)$, we get
		the limit system
		\begin{equation}\label{eq:radial-reference-limit-system}
			\begin{cases}
				\partial_\sigma N_*
				=r^{-1}\partial_r
				\bigl(r\partial_rN_*-rN_*\partial_rC_*\bigr),\\
				-r^{-1}\partial_r(r\partial_rC_*)+C_*=N_*,\\
				N_*(\sigma,1)=C_*(\sigma,1)=0.
			\end{cases}
		\end{equation}
		In particular, \eqref{eq:radial-reference-time-translates}
		gives local uniform convergence of $N_r(t_j+\sigma,1)$
		to $\partial_rN_*(\sigma,1)$, which allows us to use
		\eqref{eq:radial-reference-flux-integrable} in the next step.
		
		\noindent\textbf{Step 5. Long-time decay of the density.}
		Since $N$ and $N_*$ are nonnegative and vanish at $r=1$,
		their inward radial derivatives at $r=1$ are nonnegative.
		For each $L>0$, the convergence at $r=1$ in
		\eqref{eq:radial-reference-time-translates} and the
		integrability in \eqref{eq:radial-reference-flux-integrable}
		give
		\begin{align*}
			0	\leq\int_{-L}^L\partial_rN_*(\sigma,1)d\sigma=\lim_{j\to\infty}
			\int_{t_j-L}^{t_j+L}N_r(t,1)dt\leq\lim_{j\to\infty}
			\int_{t_j-L}^\infty N_r(t,1)dt=0.
		\end{align*}
		Since $\partial_rN_*(\cdot,1)$ is continuous and
		nonnegative and $L$ is arbitrary, we obtain
		\begin{equation}\label{eq:radial-reference-zero-limit-flux}
			\partial_rN_*(\sigma,1)=0
			\qquad\text{for every}~\sigma\in\mathbb R.
		\end{equation}
		We claim that $N_*\equiv0$. Suppose instead that
		$N_*(\sigma_0,r_0)>0$ for some $r_0>1$. Using the
		elliptic equation in \eqref{eq:radial-reference-limit-system},
		we rewrite the density equation as
		\[
		\partial_\sigma N_*-\partial_{rr}N_*
		-(r^{-1}-\partial_rC_*)\partial_rN_*
		=(N_*-C_*)N_*.
		\]
		Fix $\sigma_1>\sigma_0$ and $R>r_0$. On the cylinder
		$(\sigma_0-1,\sigma_1+1)\times(1,R)$, choose $K>0$
		such that $N_*-C_*-K\leq0$. Then
		\[
		\left[
		\partial_\sigma-\partial_{rr}
		-(r^{-1}-\partial_rC_*)\partial_r
		\right](e^{-K\sigma}N_*)
		=(N_*-C_*-K)e^{-K\sigma}N_*.
		\]
		The transformed function is a nonnegative
		$W_p^{2,1}$ strong solution with $p>3$, and its radial
		derivative is continuous up to $r=1$. The coefficients
		are bounded and continuous, and the zeroth-order
		coefficient $N_*-C_*-K$ is nonpositive. We may therefore
		apply the strong maximum principle and the parabolic
		boundary point lemma
		\cite[Chapter~II]{Lieberman1996}.
		The strong maximum principle gives
		\[
		e^{-K\sigma}N_*(\sigma,r)>0
		\qquad\text{for}\qquad
		\sigma_0<\sigma<\sigma_1+1,\quad 1<r<R.
		\]
		Apply the boundary point lemma at $(\sigma_1,1)$ on
		the cylinder
		$
		\left(\frac{\sigma_0+\sigma_1}{2},\sigma_1+1\right)
		\times(1,R)$,
		where the transformed function is positive in the
		interior and vanishes on $r=1$. Since the inward
		direction at $r=1$ is the positive radial direction, we get
		\[
		0<\partial_r(e^{-K\sigma}N_*)(\sigma_1,1)
		=e^{-K\sigma_1}\partial_rN_*(\sigma_1,1),
		\]
		contradicting \eqref{eq:radial-reference-zero-limit-flux}.
		Thus every limit obtained in Step~4 is identically zero.
		
		Fix $R>1$.
		If $N(t,\cdot)$ did not converge uniformly to zero on $[1,R]$,
		there would exist $\varepsilon>0$ and a sequence
		$t_j\to\infty$ such that
		$\sup_{1\leq r\leq R}N(t_j,r)\geq\varepsilon$
		for every $j$.
		By the compactness established in Step~4,
		a subsequence, still indexed by $j$, satisfies
		\eqref{eq:radial-reference-time-translates}.
		In particular,
		\[
		\sup_{1\leq r\leq R}
		|N(t_j,r)-N_*(0,r)|
		\to 0
		\qquad\text{as}\quad j\to\infty.
		\]
		Since $N_*\equiv0$, this gives
		$
		\varepsilon
		\leq\lim_{j\to\infty}\sup_{1\leq r\leq R}N(t_j,r)
		=0$,
		contrary to $\varepsilon>0$.
		Hence
		\begin{equation}\label{eq:radial-reference-local-density-decay}
			\lim_{t\to\infty}\sup_{1\leq r\leq R}N(t,r)=0
			\qquad\text{for every}\quad R>1.
		\end{equation}
		For $t\geq1$, 
		\eqref{eq:radial-reference-tail} gives
		\[
		\|N(t)\|_{L^\infty(\Omega)}
		\leq\max\left\{
		\sup_{1\leq r\leq R}N(t,r),\,KR^{-\frac12}
		\right\}.
		\]
		Taking the upper limit as $t\to\infty$ and using
		\eqref{eq:radial-reference-local-density-decay}, we obtain
		\[
		\limsup_{t\to\infty}\|N(t)\|_{L^\infty(\Omega)}
		\leq KR^{-\frac12}
		\qquad\text{for every }R>1.
		\]
		Letting $R\to\infty$ proves that
		$\|N(t)\|_{L^\infty(\Omega)}\to0$. Finally,  \eqref{eq:radial-reference-mass-balance} yields
		\[
		\|N(t)\|_{L^2(\Omega)}^2
		\leq \|N(t)\|_{L^\infty(\Omega)}\|N(t)\|_{L^1(\Omega)}
		\leq 2\pi\|N(t)\|_{L^\infty(\Omega)}
		\int_1^\infty rN(0,r)dr
		\to 0
		\]
		as $t\to\infty$.
	\end{proof}
	\subsection{An $L^2$ smallness criterion for global existence}
	The decay of $N$ motivates a criterion that applies when
	the density becomes small at some time. Once the density
	is sufficiently small in $L^2(\Omega)$, diffusion controls
	the chemotactic term and this norm cannot increase.
	Proposition~\ref{prop:density-continuation} then gives
	global existence.
	\begin{lemma}\label{lem:small-L2-density}
		Let $[0,T_*)$ be the maximal existence interval of a strong
		solution of \eqref{eq:physical-scalars-moments} with the boundary conditions \eqref{eq:main-boundary}
		whose initial data satisfy the assumptions of
		Proposition~\ref{prop:density-continuation}.
		There exists $\rho_*>0$, depending only on $\Omega$, such
		that if
		\[
		\|n(t_0)\|_{L^{2}(\Omega)}\leq\rho_*
		\qquad\text{for some }t_0\in[0,T_*),
		\]
		then $T_*=\infty$ and
		\begin{equation}\label{eq:small-L2-density-energy}
			\|n(t)\|_{L^{2}(\Omega)}^2+
			\int_{t_0}^t\|\nabla n(\sigma)\|_{L^{2}(\Omega)}^2 d\sigma
			\leq\|n(t_0)\|_{L^{2}(\Omega)}^2\qquad\text{for}\quad
			t\geq t_0.
		\end{equation}
	\end{lemma}
	\begin{proof}
		In the density equation of
		\eqref{eq:physical-scalars-moments}, the identity
		$e_\theta\cdot\nabla=r^{-1}\partial_\theta$ gives
		\[
		Ar^{-2}\partial_\theta n+u\cdot\nabla n
		=(Ar^{-1}e_\theta+u)\cdot\nabla n.
		\]
		Noting that $Ar^{-1}$ is independent of $\theta$ and
		$\nabla\cdot u=0$, we have
		\[
		\nabla\cdot(Ar^{-1}e_\theta+u)
		=r^{-1}\partial_\theta(Ar^{-1})+\nabla\cdot u=0.
		\]
		Moreover, by Proposition~\ref{prop:local-strong-theory} and the
		embedding $H^2(\Omega)\hookrightarrow L^\infty(\Omega)$, we get
		\[
		\|Ar^{-1}e_\theta+u(t)\|_{L^\infty(\Omega)}
		\leq A+K\|u(t)\|_{H^2(\Omega)},
		\]
		which is bounded on every compact time interval of existence.
		For each fixed $t$, let $n_k\in C_c^\infty(\Omega)$
		converge to $n(t)$ in $H^1(\Omega)$. Integration by parts
		gives
		\[
		\int_\Omega n_k(Ar^{-1}e_\theta+u)\cdot\nabla n_k dx
		=-\frac12\int_\Omega n_k^2
		\nabla\cdot(Ar^{-1}e_\theta+u)dx=0.
		\]
		Since $n_k\to n(t)$ in $H^1(\Omega)$ as $k\to\infty$, we have
		\[
		\begin{aligned}
			\|n_k\nabla n_k-n\nabla n\|_{L^1(\Omega)}
			\leq\|n_k-n\|_{L^2(\Omega)}\|\nabla n_k\|_{L^2(\Omega)}+\|n\|_{L^2(\Omega)}\|\nabla(n_k-n)\|_{L^2(\Omega)}
			\longrightarrow0.
		\end{aligned}
		\]
		The boundedness of $Ar^{-1}e_\theta+u(t)$ therefore permits passage to the limit in the transport integral, yielding
		\[
		\int_\Omega n(Ar^{-1}e_\theta+u)\cdot\nabla n dx=0.
		\]
		Thus testing the density equation with $n$ yields
		\begin{equation}\label{energy density}
			\frac12\frac{d}{dt}\|n\|_{L^{2}}^2+\|\nabla n\|_{L^{2}}^2
			=\int_\Omega n\nabla c\cdot\nabla n dx
		\end{equation}
		for almost every $t<T_*$. 
		To estimate the right-hand side, test
		$-\Delta c=n-c$ with $n-c\in H_0^1(\Omega)$. This gives
		\begin{align*}
			\|n-c\|_{L^{2}}^2=\int_\Omega\nabla c\cdot\nabla n dx-\|\nabla c\|_{L^{2}}^2=\frac14\|\nabla n\|_{L^{2}}^2
			-\left\|\nabla c-\frac12\nabla n\right\|_{L^{2}}^2
			\leq\frac14\|\nabla n\|_{L^{2}}^2.
		\end{align*}
		Since $n\in H^2(\Omega)\cap H_0^1(\Omega)$ for almost every
		$t$, we have $n^2\in H_0^1(\Omega)$ at such times. Testing
		the same elliptic equation with $n^2$ therefore yields
		\[
		\int_\Omega n\nabla c\cdot\nabla n dx
		=\frac12\int_\Omega\nabla c\cdot\nabla(n^2)dx
		=\frac12\int_\Omega(n-c)n^2 dx.
		\]
		By zero extension of $n$ and the two-dimensional
		Gagliardo--Nirenberg inequality, there holds
		\[
		\|n\|_{L^{4}}^2\leq K\|n\|_{L^{2}}\|\nabla n\|_{L^{2}}.
		\]
		Combining the preceding estimates, we obtain
		\begin{equation}\label{eq:small-L2-chemotactic-term}
			\begin{aligned}
				\left|\int_\Omega n\nabla c\cdot\nabla n dx\right|
				\leq \frac12\|n-c\|_{L^{2}}\|n\|_{L^{4}}^2\leq K_*\|n\|_{L^{2}}\|\nabla n\|_{L^{2}}^2,
			\end{aligned}
		\end{equation}
		where $K_*>0$ depends only on $\Omega$. Substituting \eqref{eq:small-L2-chemotactic-term} into \eqref{energy density} gives
		\begin{equation}\label{eq:small-L2-differential}
			\frac12\frac{d}{dt}\|n\|_{L^{2}}^2
			+(1-K_*\|n\|_{L^{2}})\|\nabla n\|_{L^{2}}^2\leq0.
		\end{equation}
		Set $\rho_*=(4K_*)^{-1}$. Suppose that
		$t_1\in(t_0,T_*)$ is the first time at which
		$\|n(t_1)\|_{L^{2}}=2\rho_*$. For $t_0\leq t<t_1$, we have
		$\|n(t)\|_{L^{2}}<2\rho_*$ and $
		1-K_*\|n(t)\|_{L^{2}}\geq\frac12$.
		Integrating \eqref{eq:small-L2-differential} over $[t_0,t]$
		therefore gives
		\[
		\|n(t)\|_{L^{2}}\leq\|n(t_0)\|_{L^{2}}\leq\rho_*\qquad \text{for}\quad t_0\leq t<t_1.
		\]
		By continuity, there holds
		\[
		2\rho_*=\|n(t_1)\|_{L^{2}}
		=\lim_{t\to t_1}\|n(t)\|_{L^{2}}
		\leq\rho_*,
		\]
		a contradiction. Hence $\|n(t)\|_{L^{2}}<2\rho_*$ for all
		$t\in[t_0,T_*)$, and integrating
		\eqref{eq:small-L2-differential} on $[t_0,t]$ yields
		\eqref{eq:small-L2-density-energy} for $t_0\leq t<T_*$.
		
		Combining \eqref{eq:small-L2-density-energy} with
		continuity on $[0,t_0]$, we obtain
		\[
		\sup_{0\leq t<T_*}\|n(t)\|_{L^{2}}
		\leq
		\max\left\{\sup_{0\leq t\leq t_0}\|n(t)\|_{L^{2}}, \rho_*\right\}
		<\infty.
		\]
		It follows from Proposition~\ref{prop:density-continuation}
		that $T_*=\infty$, so \eqref{eq:small-L2-density-energy}
		holds for every $t\geq t_0$.
	\end{proof}
	
	\section{Comparison with the radial solution}
	\label{sec:finite-time} 
	
	The density $n^A$ is now compared with the radial density
	$N$ on fixed finite time intervals, as required in
	Section~\ref{sec:main-proof}. The comparison also yields
	the improved $L^2(\Omega)$ bound determined by the constant
	$B$ in \eqref{eq:bootstrap-level}.
	\begin{proposition}\label{prop:finite-time-averaging}
		Assume the initial and boundary conditions of
		Theorem \ref{thm:main}.
		For every $T>0$, there exists $A_T\geq1$, depending on
		$T$ and the fixed initial data, such that the local strong
		solution extends to $[0,T]$ whenever $A\geq A_T$, and
		$\sup_{0\leq t\leq T}\|n^A(t)\|_{L^{2}(\Omega)}\leq \frac{B}{2}.$
		Moreover,
		\begin{equation}\label{eq:finite-time-averaging}
			\sup_{0\leq t\leq T}
			\left(
			\int_1^\infty
			r|n_0^A(t,r)-N(t,r)|^2 dr
			\right)^{\frac12}
			\to 0\qquad\text{as}\quad A\to\infty,
		\end{equation}
		and, for every $0<\delta<T$,
		\begin{equation}\label{n omega nonzero}
			\sup_{\delta\leq t\leq T}
			\|n_{\neq}^A(t)\|_{L^{2}(\Omega)}
			\to 0\qquad\text{as}\quad A\to\infty.
		\end{equation}
		Here $B$ is defined in
		\eqref{eq:bootstrap-level}, and $(N,C)$ is the radial solution
		of \eqref{eq:radial-reference-problem} with
		$N(0)=(n_{\mathrm{in}})_0$.
	\end{proposition}
	The proof combines decay of the nonzero angular modes with
	an energy estimate for $n_0^A-N$. The estimates are first
	derived on the interval where
	$\|n^A(t)\|_{L^2(\Omega)}<B$. The resulting comparison
	improves this bound and gives existence through time $T$.
	\subsection{Uniform bounds for the density, vorticity, and velocity}
	Decay of the nonzero angular modes requires bounds
	independent of $A$ for the density, vorticity, velocity,
	and the vector field $un+n\nabla c$ in
	\eqref{eq:physical-scalars-moments}. Weighted estimates
	control the density and vorticity at large radii through
	\[
	\|\mathbf{1}_{\{r>R\}}f\|_{L^2(\Omega)}
	\leq R^{-2}\|r^2f\|_{L^2(\Omega)},
	\qquad f=n,w.
	\]
	Under an $L^2(\Omega)$ bound on the density, energy
	estimates give the required weighted bounds and the
	spatial regularity and bounds at large radii for
	$un+n\nabla c$.
	\begin{lemma}
		\label{lem:finite-time-uniform-moments}
		Let $0<T_*\leq T<\infty$, and let $(n,c,u,w)$ solve
		\eqref{eq:physical-scalars-moments} with the boundary conditions
		\eqref{eq:main-boundary}, in the strong solution class of
		Proposition~\ref{prop:local-strong-theory} on every compact subinterval of
		$[0,T_*)$. Suppose that $n\geq0$,
		\begin{equation}\label{eq:finite-time-L2-bootstrap}
			\sup_{0\leq t<T_*}\|n(t)\|_{L^2(\Omega)}\leq B,
		\end{equation}
		and
		$n_{\mathrm{in}}\in L^\infty(\Omega)\cap H_0^1(\Omega)$,
		$w_{\mathrm{in}}\in H_0^1(\Omega)$,
		$r^2n_{\mathrm{in}},\,r^2w_{\mathrm{in}}\in L^2(\Omega)$.
		Then
		\begin{equation}\label{eq:finite-time-velocity-bounds}
			\begin{aligned}
				&\|n\|_{L^\infty(0,T_*;L^\infty(\Omega))}
				+\|w\|_{L^\infty(0,T_*;L^2(\Omega))}\\&\quad+\|\nabla n\|_{L^2(0,T_*;L^2(\Omega))}
				+\|\nabla w\|_{L^2(0,T_*;L^2(\Omega))}
				\leq K_T,\\
				&\|r^2n\|_{L^\infty(0,T_*;L^2(\Omega))}
				+\|r^2w\|_{L^\infty(0,T_*;L^2(\Omega))}\\
				&\quad+\|r^2\nabla n\|_{L^2(0,T_*;L^2(\Omega))}
				+\|r^2\nabla w\|_{L^2(0,T_*;L^2(\Omega))}
				\leq K_T,\\
				&\|u\|_{L^\infty(0,T_*;H^1(\Omega))}
				+\|u\|_{L^2(0,T_*;H^2(\Omega))}
				+\|u\|_{L^2(0,T_*;L^\infty(\Omega))}
				\leq K_T.
			\end{aligned}
		\end{equation}
		For every $0<s<1$, the vector field $un+n\nabla c$ satisfies
		\begin{equation}\label{eq:finite-time-flux-regularity}
			\begin{aligned}
				&\|un+n\nabla c\|_{L^2(0,T_*;H^s(\Omega))}
				\leq K_{T,s},\\
				&\|\mathbf1_{\{r>R\}}(un+n\nabla c)\|_{L^2(0,T_*;L^2(\Omega))}
				\leq K_TR^{-2}\quad \text{for}\quad R\geq2.
			\end{aligned}
		\end{equation}
		Here $\mathbf{1}_{\{r>R\}}$ denotes the indicator function of
		$\{x\in\Omega: |x|>R\}$.
		The constants depend only on $\Omega$, $T$, $B$, and the displayed
		initial norms, and, for $K_{T,s}$, also on $s$. They are independent
		of $A$ and $T_*$.
	\end{lemma}
	\begin{proof}
		The estimates for $n$ and $w$ are established first.
		They then give the velocity bounds in
		\eqref{eq:finite-time-velocity-bounds} and the estimates
		for $un+n\nabla c$ in \eqref{eq:finite-time-flux-regularity}.
		
		\noindent\textbf{Step 1. Unweighted bounds for $n$ and $w$ in
			\eqref{eq:finite-time-velocity-bounds}.}
		Lemma~\ref{lem:exterior-moser} and
		\eqref{eq:finite-time-L2-bootstrap} give
		\begin{equation}\label{eq:finite-time-density-Linfty}
			\|n\|_{L^\infty(0,T_*;L^\infty(\Omega))}
			\leq K\bigl(\|n_{\mathrm{in}}\|_{L^\infty(\Omega)}+B\bigr)
			(1+B^2).
		\end{equation}
		By H\"older's inequality,
		$\|n(t)\|_{L^4(\Omega)}\leq
		\|n(t)\|_{L^2(\Omega)}^{\frac12}\|n(t)\|_{L^\infty(\Omega)}^{\frac12}$.
		Thus \eqref{eq:finite-time-L2-bootstrap} and
		\eqref{eq:finite-time-density-Linfty} also bound $n$ in
		$L^\infty(0,T_*;L^4(\Omega))$.
		Applying the $H^2(\Omega)$ and $W^{2,4}(\Omega)$ estimates of
		Lemma~\ref{lem:screened-regularity} and the embedding
		$W^{2,4}(\Omega)\hookrightarrow W^{1,\infty}(\Omega)$, we obtain
		\begin{equation}\label{eq:finite-time-chemical-gradient}
			\|c\|_{L^\infty(0,T_*;W^{1,\infty}(\Omega))}
			+\|\nabla c\|_{L^\infty(0,T_*;H^1(\Omega))}
			\leq K_T.
		\end{equation}
		By Proposition~\ref{prop:local-strong-theory},
		$n_t,w_t\in L^2(0,\tau;L^2(\Omega))$ for every $\tau<T_*$,
		where the regularity of $w_t$ follows from
		$w=\operatorname{curl}u$. For $f=n,w$, the angular term vanishes
		by periodicity. Since $\nabla\cdot u=0$ and
		$u\in H^2(\Omega)\hookrightarrow L^\infty(\Omega)$, approximation
		of $f\in H_0^1(\Omega)$ by compactly supported smooth functions
		also eliminates the transport by $u$. Thus
		\begin{equation}\label{eq:finite-time-transport-cancellation}
			\int_\Omega Ar^{-2}f\partial_\theta fdx=0,
			\qquad
			\int_\Omega f u\cdot\nabla f dx=0\qquad\text{for}\quad
			f=n,w.
		\end{equation}
		For almost every time, $n\in H^2(\Omega)\cap H_0^1(\Omega)$.
		The embedding $H^2(\Omega)\hookrightarrow L^\infty(\Omega)$
		and the Sobolev chain rule give $\frac{1}{2}n^{2}\in H_0^1(\Omega)$.
		Testing the elliptic equation $-\Delta c=n-c$ with $\frac{1}{2}n^{2}$
		therefore yields
		\begin{equation}\label{eq:finite-time-chemotaxis-identity}
			\int_\Omega n\nabla c\cdot\nabla ndx
			=\int_\Omega\nabla c\cdot\nabla\left(\frac{n^{2}}{2}\right)dx
			=\frac12\int_\Omega(n-c)n^2 dx.
		\end{equation}
		Testing the density and vorticity equations in
		\eqref{eq:physical-scalars-moments} with $n$ and $w$, respectively,
		and applying \eqref{eq:finite-time-transport-cancellation} and
		\eqref{eq:finite-time-chemotaxis-identity}, we obtain the following
		inequalities by \eqref{eq:finite-time-L2-bootstrap},
		\eqref{eq:finite-time-density-Linfty},
		\eqref{eq:finite-time-chemical-gradient}, and Young's inequality:
		\begin{equation}\label{eq:finite-time-basic-energy-inequalities}
			\begin{aligned}
				\frac12\frac{d}{dt}\|n\|_{L^2(\Omega)}^2
				+\|\nabla n\|_{L^2(\Omega)}^2
				&=\frac12\int_\Omega(n-c)n^2 dx\\
				&\leq
				\frac12\bigl(\|n\|_{L^\infty(\Omega)}
				+\|c\|_{L^\infty(\Omega)}\bigr)\|n\|_{L^2(\Omega)}^2
				\leq K_T B^2,\\
				\frac12\frac{d}{dt}\|w\|_{L^2(\Omega)}^2
				+\|\nabla w\|_{L^2(\Omega)}^2
				&=\int_\Omega n e_\theta\cdot\nabla w dx\leq
				\frac12\|n\|_{L^2(\Omega)}^2
				+\frac12\|\nabla w\|_{L^2(\Omega)}^2.
			\end{aligned}
		\end{equation}
		Integrating the density inequality in
		\eqref{eq:finite-time-basic-energy-inequalities} on $[0,t]$
		and letting $t\to T_*$ gives
		\begin{equation}\label{eq:finite-time-density-dissipation}
			\int_0^{T_*}\|\nabla n(t)\|_{L^2(\Omega)}^2 dt
			\leq\frac12\|n_{\mathrm{in}}\|_{L^2(\Omega)}^2+K_TB^2T.
		\end{equation}
		Integrating the vorticity inequality in
		\eqref{eq:finite-time-basic-energy-inequalities} and using
		\eqref{eq:finite-time-L2-bootstrap} yields
		\begin{equation}\label{eq:finite-time-g-bound}
			\sup_{0\leq t<T_*}\|w(t)\|_{L^2(\Omega)}^2
			+\int_0^{T_*}\|\nabla w(t)\|_{L^2(\Omega)}^2 dt
			\leq 2\bigl(\|w_{\mathrm{in}}\|_{L^2(\Omega)}^2+B^2T\bigr).
		\end{equation}
		Combining \eqref{eq:finite-time-density-Linfty},
		\eqref{eq:finite-time-density-dissipation}, and
		\eqref{eq:finite-time-g-bound} proves the unweighted estimate
		in \eqref{eq:finite-time-velocity-bounds}.
		
		\noindent\textbf{Step 2. Weighted bounds for $n$ and $w$ in
			\eqref{eq:finite-time-velocity-bounds}.}
		Fix $0<\tau<T_*$ and $R>1$. To justify the energy identity
		with weight $r^4$, first set
		\begin{equation}\label{eq:finite-time-bounded-weight}
			\rho_R(r)=\frac{r^2}{1+r^{2}R^{-2}},
			\qquad
			0<\rho_R\leq\min\{r^2,R^2\},
			\qquad
			0\leq\rho_R'\leq\frac{2\rho_R}{r}.
		\end{equation}
		By \eqref{eq:finite-time-bounded-weight}, $\rho_R$ and $\rho_R'$
		are bounded, so $\rho_R^2 f\in H_0^1(\Omega)$ for $f=n,w$.
		By periodicity in $\theta$
		and integration by parts, the diffusion and transport terms satisfy
		\begin{equation}\label{eq:finite-time-weighted-testing-identities}
			\begin{aligned}
				\int_\Omega \nabla f\cdot\nabla(\rho_R^2 f)dx
				&=\|\rho_R\nabla f\|_{L^2(\Omega)}^2
				+2\int_\Omega \rho_R\rho_R' f\partial_r f dx,\\
				\int_\Omega \rho_R^2 f
				(Ar^{-2}\partial_\theta f+&u\cdot\nabla f)dx
				=-\int_\Omega \rho_R\rho_R' u_r f^2 dx.
			\end{aligned}
		\end{equation}
		Testing the density and vorticity equations in
		\eqref{eq:physical-scalars-moments} with $\rho_R^2n$ and
		$\rho_R^2w$, respectively, and using
		\eqref{eq:finite-time-weighted-testing-identities} and
		$e_\theta\cdot\nabla\rho_R=0$, we obtain
		\begin{equation}\label{eq:finite-time-bounded-weight-identity}
			\begin{aligned}
				&\frac12\frac{d}{dt}\sum_{f=n,w}\|\rho_R f\|_{L^2(\Omega)}^2
				+\sum_{f=n,w}\|\rho_R\nabla f\|_{L^2(\Omega)}^2\\
				=~&-2\sum_{f=n,w}\int_\Omega
				\rho_R\rho_R' f\partial_r f dx
				+\int_\Omega\rho_R\rho_R' u_r(n^2+w^2)dx\\
				&+\int_\Omega\rho_R^2 n\nabla c\cdot\nabla n dx
				+2\int_\Omega\rho_R\rho_R' n^2\partial_rc dx
				+\int_\Omega\rho_R^2 n e_\theta\cdot\nabla w dx.
			\end{aligned}
		\end{equation}
		By \eqref{eq:finite-time-bounded-weight} and $r\geq1$, we have
		$\rho_R'\leq2\rho_R$. Hence the Cauchy--Schwarz and Young
		inequalities give, for $f=n,w$,
		\begin{align*}
			2\left|\int_\Omega \rho_R\rho_R' f\partial_r f dx\right|\leq 4\|\rho_R f\|_{L^2(\Omega)}
			\|\rho_R\nabla f\|_{L^2(\Omega)}\leq \frac14\|\rho_R\nabla f\|_{L^2(\Omega)}^2
			+K\|\rho_R f\|_{L^2(\Omega)}^2.
		\end{align*}
		Applying the same inequalities to the terms containing
		$n\nabla c\cdot\nabla n$ and $ne_\theta\cdot\nabla w$, we obtain
		\begin{align*}
			&\left|\int_\Omega \rho_R^2n\nabla c\cdot\nabla n dx\right|
			+\left|\int_\Omega \rho_R^2ne_\theta\cdot\nabla w dx\right|\\
			\leq~&
			\|\nabla c\|_{L^\infty(\Omega)}
			\|\rho_Rn\|_{L^2(\Omega)}\|\rho_R\nabla n\|_{L^2(\Omega)}
			+\|\rho_Rn\|_{L^2(\Omega)}\|\rho_R\nabla w\|_{L^2(\Omega)}\\
			\leq~&
			\frac12\sum_{f=n,w}\|\rho_R\nabla f\|_{L^2(\Omega)}^2
			+K\bigl(1+\|\nabla c\|_{L^\infty(\Omega)}^2\bigr)
			\|\rho_Rn\|_{L^2(\Omega)}^2.
		\end{align*}
		For the two remaining terms, using $\rho_R'\leq2\rho_R$ again gives
		\begin{align*}
			&\left|\int_\Omega \rho_R\rho_R' u_r(n^2+w^2)dx\right|
			+2\left|\int_\Omega \rho_R\rho_R' n^2\partial_r c dx\right|\\
			\leq~&
			2\|u_r\|_{L^\infty(\Omega)}
			\sum_{f=n,w}\|\rho_R f\|_{L^2(\Omega)}^2
			+4\|\nabla c\|_{L^\infty(\Omega)}
			\|\rho_Rn\|_{L^2(\Omega)}^2.
		\end{align*}
		Substituting these estimates into
		\eqref{eq:finite-time-bounded-weight-identity} and using
		$\|\nabla c\|_{L^\infty(\Omega)}
		\leq 1+\|\nabla c\|_{L^\infty(\Omega)}^2$, we obtain
		\begin{equation}\label{eq:finite-time-bounded-weight-inequality}
			\begin{aligned}
				&\frac{d}{dt}\sum_{f=n,w}\|\rho_R f\|_{L^2(\Omega)}^2
				+\frac12\sum_{f=n,w}\|\rho_R\nabla f\|_{L^2(\Omega)}^2\\
				\leq~& K\bigl(1+\|u_r\|_{L^\infty(\Omega)}
				+\|\nabla c\|_{L^\infty(\Omega)}^2\bigr)
				\sum_{f=n,w}\|\rho_R f\|_{L^2(\Omega)}^2.
			\end{aligned}
		\end{equation}
		By Proposition~\ref{prop:local-strong-theory} and the embedding $H^2(\Omega)\hookrightarrow L^\infty(\Omega)$,
		$u_r$ is bounded on $[0,\tau]\times\Omega$.
		Together with \eqref{eq:finite-time-chemical-gradient}, this
		bounds the coefficient on the right-hand side of
		\eqref{eq:finite-time-bounded-weight-inequality}
		independently of $R$. Moreover, \eqref{eq:finite-time-bounded-weight}
		gives
		\[
		\|\rho_R f(0)\|_{L^2(\Omega)}
		\leq \|r^2 f(0)\|_{L^2(\Omega)},
		\qquad f=n,w.
		\]
		Applying Gronwall's inequality to
		\eqref{eq:finite-time-bounded-weight-inequality} first bounds
		the weighted $L^2(\Omega)$ norms. Integrating the same
		inequality over $[0,\tau]$ then bounds the weighted gradients.
		Consequently,
		\[
		\sup_{0\leq t\leq\tau}
		\sum_{f=n,w}\|\rho_R f(t)\|_{L^2(\Omega)}^2
		+\int_0^\tau\sum_{f=n,w}
		\|\rho_R\nabla f(t)\|_{L^2(\Omega)}^2 dt
		\leq K_{A,\tau}
		\sum_{f=n,w}\|r^2 f(0)\|_{L^2(\Omega)}^2,
		\]
		where $K_{A,\tau}$ is independent of $R$ but may depend
		on the local strong solution bound for $u$ on $[0,\tau]$.
		This estimate is used only to justify the limit $R\to\infty$.
		A bound independent of $A$ and $\tau$ is derived below from
		\eqref{eq:moment-energy-preliminary}.
		Since $\rho_R(r)\to r^2$ as $R\to\infty$, monotone
		convergence gives
		\begin{equation}\label{eq:finite-time-weighted-local-integrability}
			r^2n,\,r^2w\in L^\infty(0,\tau;L^2(\Omega)),
			\qquad
			r^2\nabla n,\,r^2\nabla w\in L^2(0,\tau;L^2(\Omega)).
		\end{equation}
		The same argument at each fixed $t\in[0,\tau]$ gives
		$r^2n(t),\,r^2w(t)\in L^2(\Omega)$.
		Integrating \eqref{eq:finite-time-bounded-weight-identity}
		over $[0,t]$, for $0\leq t\leq\tau$, gives
		\begin{align}\label{eq:finite-time-bounded-weight-integrated}
			&\frac12\sum_{f=n,w}
			\left(\|\rho_R f(t)\|_{L^2(\Omega)}^2
			-\|\rho_R f(0)\|_{L^2(\Omega)}^2\right)+\int_0^t\sum_{f=n,w}
			\|\rho_R\nabla f(s)\|_{L^2(\Omega)}^2 ds\nonumber\\
			=~&\int_0^t\int_\Omega\biggl[
			-2\sum_{f=n,w}\rho_R\rho_R'f\partial_r f
			+\rho_R\rho_R'u_r(n^2+w^2)+\rho_R^2n\nabla c\cdot\nabla n
			\\\nonumber&\qquad+2\rho_R\rho_R'n^2\partial_r c
			+\rho_R^2ne_\theta\cdot\nabla w
			\biggr]dxds.
		\end{align}
		By \eqref{eq:finite-time-bounded-weight}, there holds
		\[
		\rho_R^2\to r^4,\quad
		0\leq\rho_R\rho_R'
		=\frac{2r^3}{(1+r^2R^{-2})^3}\leq 2r^3,\quad
		\rho_R\rho_R'\longrightarrow 2r^3\quad\text{as}\quad
		R\to\infty.
		\]
		Thus the two energy terms and the gradient integral on the
		left of \eqref{eq:finite-time-bounded-weight-integrated}
		converge by monotone convergence. The energy terms at fixed
		times are finite by the preceding bound independent of $R$.
		On the right of \eqref{eq:finite-time-bounded-weight-integrated},
		the coefficients $\rho_R^2$ and $\rho_R\rho_R'$ are bounded by
		$r^4$ and $2r^3$, respectively. Since $r\geq1$, we have, for
		$f=n,w$,
		\[
		|\rho_R\rho_R'f\partial_r f|
		\leq 2r^3|f||\nabla f|
		\leq r^4\bigl(|f|^2+|\nabla f|^2\bigr).
		\]
		The other integrands are bounded in the same way: the terms
		containing $u_r$ or $\partial_r c$ are bounded by their
		$L^\infty(\Omega)$ norms times a multiple of $r^4(n^2+w^2)$,
		and the terms containing $\nabla n$ or $\nabla w$ are bounded
		using $2|n||\nabla f|\leq n^2+|\nabla f|^2$.
		Since $u_r$ and $\nabla c$ are bounded on $[0,\tau]\times\Omega$,
		each integrand is bounded in absolute value by a constant
		independent of $R$ times $r^4(n^2+w^2+|\nabla n|^2+|\nabla w|^2)$.
		This function is integrable by
		\eqref{eq:finite-time-weighted-local-integrability}. Passing
		to the limit by dominated convergence on the right gives
		\begin{align*}
			&\frac12\sum_{f=n,w}
			\left(\|r^2f(t)\|_{L^2(\Omega)}^2
			-\|r^2f(0)\|_{L^2(\Omega)}^2\right)
			+\int_0^t\sum_{f=n,w}\|r^2\nabla f(s)\|_{L^2(\Omega)}^2 ds\\
			=~&\int_0^t\int_\Omega\biggl[
			-4\sum_{f=n,w}r^3f\partial_r f
			+2r^3u_r(n^2+w^2)
			+r^4n\nabla c\cdot\nabla n\\
			&\hspace{28mm}+4r^3n^2\partial_r c
			+r^4ne_\theta\cdot\nabla w
			\biggr]dxds.
		\end{align*}
		Moving the gradient integral to the right expresses the
		weighted energy as its initial value plus the integral over
		$[0,t]$ of a function in $L^1(0,\tau)$. Hence the energy is
		absolutely continuous on $[0,\tau]$. By the fundamental theorem
		of calculus for Lebesgue integrals, the following identity
		holds for almost every $t\in(0,\tau)$:
		\begin{align}\label{eq:finite-time-r4-energy-identity}
			&\frac12\frac{d}{dt}\sum_{f=n,w}\|r^2f\|_{L^2(\Omega)}^2
			+\sum_{f=n,w}\|r^2\nabla f\|_{L^2(\Omega)}^2\notag\\
			=~&-4\sum_{f=n,w}\int_\Omega r^3f\partial_r fdx
			+2\int_\Omega r^3u_r(n^2+w^2)dx\\\notag
			&+\int_\Omega r^4n\nabla c\cdot\nabla ndx
			+4\int_\Omega r^3n^2\partial_r cdx
			+\int_\Omega r^4ne_\theta\cdot\nabla w dx.
		\end{align}
		We now use \eqref{eq:finite-time-r4-energy-identity}
		to derive a bound independent of $A$. We first estimate the
		terms other than the integral containing $u_r$. Using $r^3\leq r^4$,
		the Cauchy--Schwarz and Young inequalities give, for $f=n,w$,
		\begin{align*}
			4\left|\int_\Omega r^3 f\partial_r fdx\right|
			\leq 4\|r^2f\|_{L^2(\Omega)}\|r^2\nabla f\|_{L^2(\Omega)}\leq \frac14\|r^2\nabla f\|_{L^2(\Omega)}^2
			+K\|r^2f\|_{L^2(\Omega)}^2.
		\end{align*}
		For the chemotactic terms, the bound on $\nabla c$ in
		\eqref{eq:finite-time-chemical-gradient} and Young's
		inequality give
		\begin{align*}
			&\left|\int_\Omega r^4n\nabla c\cdot\nabla n dx\right|
			+4\left|\int_\Omega r^3n^2\partial_r c dx\right|\\
			\leq~&\|\nabla c\|_{L^\infty(\Omega)}
			\left(\|r^2n\|_{L^2(\Omega)}\|r^2\nabla n\|_{L^2(\Omega)}
			+4\|r^2n\|_{L^2(\Omega)}^2\right)\\\leq~& \frac12\|r^2\nabla n\|_{L^2(\Omega)}^2
			+K_T\|r^2n\|_{L^2(\Omega)}^2.
		\end{align*}
		Since $|e_\theta|=1$, the vorticity source satisfies
		\begin{align*}
			\left|\int_\Omega r^4ne_\theta\cdot\nabla w dx\right|\leq \|r^2n\|_{L^2(\Omega)}\|r^2\nabla w\|_{L^2(\Omega)}\leq \frac12\|r^2\nabla w\|_{L^2(\Omega)}^2
			+K\|r^2n\|_{L^2(\Omega)}^2.
		\end{align*}
		For the transport term, write
		$r^3|f|^2=(r^4|f|^2)^{\frac34}(|f|^2)^{\frac14}$.
		H\"older's inequality with exponents $\frac43$ and $4$ yields
		\begin{equation*}
			\int_\Omega r^3|f|^2 dx
			\leq \|f\|_{L^2(\Omega)}^{\frac12}
			\|r^2f\|_{L^2(\Omega)}^{\frac32},
			\qquad f=n,w.
		\end{equation*}
		Applying it with the
		$L^2(\Omega)$ bounds for $n$ and $w$ in
		\eqref{eq:finite-time-L2-bootstrap} and
		\eqref{eq:finite-time-g-bound}, respectively, gives
		\begin{equation}\label{eq:finite-time-weighted-transport}
			\begin{aligned}
				2\left|\int_\Omega r^3u_r(n^2+w^2)dx\right|
				\leq K_T\|u_r\|_{L^\infty(\Omega)}
				\left(\sum_{f=n,w}\|r^2f\|_{L^2(\Omega)}^2\right)^{\frac34}.
			\end{aligned}
		\end{equation}
		Here we used
		$a^{\frac34}+b^{\frac34}\leq 2^{\frac14}(a+b)^{\frac34}$ for $a,b\geq0$.
		Inserting \eqref{eq:finite-time-weighted-transport} and the
		estimates for diffusion, chemotaxis, and the vorticity source
		into \eqref{eq:finite-time-r4-energy-identity}, we get
		\begin{equation}\label{eq:moment-energy-preliminary}
			\begin{aligned}
				&\frac{d}{dt}\sum_{f=n,w}\|r^2f\|_{L^2(\Omega)}^2
				+\frac12\sum_{f=n,w}\|r^2\nabla f\|_{L^2(\Omega)}^2\\
				\leq~&K_T\sum_{f=n,w}\|r^2f\|_{L^2(\Omega)}^2
				+K_T\|u_r\|_{L^\infty(\Omega)}
				\left(\sum_{f=n,w}\|r^2f\|_{L^2(\Omega)}^2\right)^{\frac34}.
			\end{aligned}
		\end{equation}
		To estimate $u_r$ in \eqref{eq:moment-energy-preliminary}, we use
		the radial velocity estimate \eqref{eq:radial-velocity-interpolation}.
		By Cauchy--Schwarz and $\int_\Omega r^{-4}dx=\pi$, it holds that
		$\|w\|_{L^1(\Omega)}
		\leq \sqrt{\pi}\|r^2w\|_{L^2(\Omega)}$.
		Since $w\in H_0^1(\Omega)$, the Gagliardo--Nirenberg inequality
		applied to its zero extension gives
		\[
		\|w\|_{L^4(\Omega)}
		\leq K\|w\|_{L^2(\Omega)}^{\frac12}
		\|\nabla w\|_{L^2(\Omega)}^{\frac12}.
		\]
		Using these two bounds in \eqref{eq:radial-velocity-interpolation}
		and controlling $\|w\|_{L^2(\Omega)}$ by
		\eqref{eq:finite-time-g-bound}, we obtain
		\begin{equation}\label{eq:finite-time-radial-velocity}
			\|u_r\|_{L^\infty(\Omega)}
			\leq K_T\|r^2w\|_{L^2(\Omega)}^{\frac13}
			\|\nabla w\|_{L^2(\Omega)}^{\frac13}.
		\end{equation}
		Substituting \eqref{eq:finite-time-radial-velocity} into
		\eqref{eq:moment-energy-preliminary} bounds the transport term by
		\[
		K_T\|\nabla w\|_{L^2(\Omega)}^{\frac13}
		\left(\sum_{f=n,w}\|r^2f\|_{L^2(\Omega)}^2\right)^{\frac34+\frac16}.
		\]
		Since $\frac34+\frac16=\frac{11}{12}<1$, the inequalities $x^{\frac{11}{12}}\leq 1+x$
		and $x^{\frac13}\leq1+x^2$ for $x\geq0$ give
		\begin{equation}\label{eq:finite-time-moment-gronwall}
			\begin{aligned}
				&\frac{d}{dt}\sum_{f=n,w}\|r^2f\|_{L^2(\Omega)}^2
				+\frac12\sum_{f=n,w}\|r^2\nabla f\|_{L^2(\Omega)}^2\\\leq~& K_T\bigl(1+\|\nabla w\|_{L^2(\Omega)}^2\bigr)
				\left(1+\sum_{f=n,w}\|r^2f\|_{L^2(\Omega)}^2\right).
			\end{aligned}
		\end{equation}
		It follows from \eqref{eq:finite-time-g-bound} and $\tau\leq T$ that
		$\int_0^\tau \bigl(1+\|\nabla w(s)\|_{L^2(\Omega)}^2\bigr)ds
		\leq K_T$.
		Gronwall's inequality applied to \eqref{eq:finite-time-moment-gronwall}
		therefore bounds the weighted $L^2(\Omega)$ norms on $[0,\tau]$
		independently of $A$ and $\tau$. Integrating
		\eqref{eq:finite-time-moment-gronwall} and using this bound
		also controls the weighted gradients in $L^2(0,\tau;L^2(\Omega))$.
		Letting $\tau\to T_*$ proves the weighted estimate in
		\eqref{eq:finite-time-velocity-bounds}.
		
		\noindent\textbf{Step 3. Velocity estimates in
			\eqref{eq:finite-time-velocity-bounds}.}
		The weighted bound for $w$ from Step~2 and
		\eqref{eq:weighted-vorticity-velocity} imply
		$\|u\|_{L^\infty(0,T_*;L^2(\Omega))}\leq K_T$.
		Since $T_*\leq T$, we also have
		\[
		\|u\|_{L^2(0,T_*;L^2(\Omega))}
		\leq \sqrt{T}\,\|u\|_{L^\infty(0,T_*;L^2(\Omega))}
		\leq K_T.
		\]
		On the other hand, \eqref{eq:finite-time-g-bound} and $T_*\leq T$
		give
		\[
		\|w\|_{L^\infty(0,T_*;L^2(\Omega))}
		+\|w\|_{L^2(0,T_*;H^1(\Omega))}
		\leq K_T.
		\]
		Applying \eqref{eq:exterior-hodge-H1} and
		\eqref{eq:exterior-hodge-H2} in the corresponding time norms, we obtain
		\[
		\begin{aligned}
			\|u\|_{L^\infty(0,T_*;H^1(\Omega))}
			&\leq K\Bigl(
			\|u\|_{L^\infty(0,T_*;L^2(\Omega))}
			+\|w\|_{L^\infty(0,T_*;L^2(\Omega))}
			\Bigr)
			\leq K_T,\\
			\|u\|_{L^2(0,T_*;H^2(\Omega))}
			&\leq K\Bigl(
			\|u\|_{L^2(0,T_*;L^2(\Omega))}
			+\|w\|_{L^2(0,T_*;H^1(\Omega))}
			\Bigr)
			\leq K_T.
		\end{aligned}
		\]
		Finally, the embedding $H^2(\Omega)\hookrightarrow L^\infty(\Omega)$
		yields
		\[
		\|u\|_{L^2(0,T_*;L^\infty(\Omega))}
		\leq K\|u\|_{L^2(0,T_*;H^2(\Omega))}
		\leq K_T.
		\]
		
		\noindent\textbf{Step 4. Estimates for $un+n\nabla c$ in
			\eqref{eq:finite-time-flux-regularity}.}
		We first estimate these terms in $L^2(0,T_*;H^s(\Omega))$,
		where $0<s<1$. By \eqref{eq:finite-time-L2-bootstrap} and the
		unweighted estimates in \eqref{eq:finite-time-velocity-bounds}, we obtain
		\begin{equation}\label{eq:finite-time-density-vorticity-H1}
			\|n\|_{L^2(0,T_*;H^1(\Omega))}\leq K_T.
		\end{equation}
		Applying the product estimate of Lemma~\ref{lem:exterior-product}
		at each time and then taking $L^2$ in time gives
		\begin{align*}
			&\|un\|_{L^2(0,T_*;H^s(\Omega))}
			\leq K_s\|u\|_{L^\infty(0,T_*;H^1(\Omega))}
			\|n\|_{L^2(0,T_*;H^1(\Omega))},\\
			&\|n\nabla c\|_{L^2(0,T_*;H^s(\Omega))}
			\leq K_s\|\nabla c\|_{L^\infty(0,T_*;H^1(\Omega))}
			\|n\|_{L^2(0,T_*;H^1(\Omega))}.
		\end{align*}
		Using \eqref{eq:finite-time-density-vorticity-H1},
		the velocity estimates in \eqref{eq:finite-time-velocity-bounds},
		and \eqref{eq:finite-time-chemical-gradient} in these
		inequalities gives the $H^s(\Omega)$ estimate in
		\eqref{eq:finite-time-flux-regularity}.
		
		It remains to bound the same terms on $\{r>R\}$.
		Using $|n|\leq R^{-2}|r^2n|$ on this region gives
		\[
		\|\mathbf1_{\{r>R\}}un\|_{L^2(0,T_*;L^2(\Omega))}
		\leq R^{-2}\|u\|_{L^2(0,T_*;L^\infty(\Omega))}
		\|r^2n\|_{L^\infty(0,T_*;L^2(\Omega))}.
		\]
		Using the same inequality for $n$ and $T_*\leq T$ yields
		\begin{align*}
			\|\mathbf1_{\{r>R\}}n\nabla c\|_{L^2(0,T_*;L^2(\Omega))}
			\leq R^{-2}\sqrt{T}
			\|\nabla c\|_{L^\infty(0,T_*;L^\infty(\Omega))}
			\|r^2n\|_{L^\infty(0,T_*;L^2(\Omega))}.
		\end{align*}
		Combining these bounds with \eqref{eq:finite-time-chemical-gradient}
		and \eqref{eq:finite-time-velocity-bounds} gives the estimate
		on $\{r>R\}$ in \eqref{eq:finite-time-flux-regularity}.
	\end{proof}
	\subsection{Decay of the nonzero angular modes}
	Decay of the nonzero angular modes is first proved for
	a scalar transport--diffusion equation with a source term
	in divergence form. The bounds in
	Lemma~\ref{lem:finite-time-uniform-moments} verify its
	hypotheses for the density equation.
	\begin{lemma}\label{lem:finite-time-fast-rotation}
		Fix $T>0$ and $0<s<1$. For each $A\geq1$, let
		$0<\tau_A\leq T$ and let $Q^A$ be a vector field on
		$(0,\tau_A)\times\Omega$. Suppose that $f^A$ is a real-valued
		energy solution of
		\begin{equation}\label{eq:fast-rotation-scalar}
			\partial_t f^A+Ar^{-2}\partial_\theta f^A-\Delta f^A
			=-\diver Q^A,
			\qquad f^A|_{r=1}=0,
		\end{equation}
		with
		$f^A\in C([0,\tau_A);L^2(\Omega))\cap L^2(0,\tau_A;H_0^1(\Omega))$.
		For $R>1$, set $\Omega_R=\{x\in\Omega:r<R\}$. Assume that
		\begin{align}
			&\sup_{A\geq1}\Bigl(
			\|f^A(0)\|_{L^2(\Omega)}
			+\|r^2f^A\|_{L^\infty(0,\tau_A;L^2(\Omega))}
			+\|Q^A\|_{L^2(0,\tau_A;L^2(\Omega))}
			\Bigr)<\infty,
			\label{eq:fast-rotation-bounds}\\
			&\sup_{A\geq1}
			\|Q^A\|_{L^2(0,\tau_A;H^s(\Omega_R))}<\infty
			\qquad\text{for every }R>1,
			\label{eq:fast-rotation-local-regularity}\\
			&\lim_{R\to\infty}\sup_{A\geq1}
			\|\mathbf{1}_{\{r>R\}}Q^A\|_{L^2(0,\tau_A;L^2(\Omega))}=0.
			\label{eq:fast-rotation-flux-tail}
		\end{align}
		Then, for every $0<\delta<T$, there holds
		\begin{equation}\label{eq:fast-rotation-limit}
			\lim_{A\to\infty}
			\sup_{\delta\leq t<\tau_A}
			\|f^A_{\neq}(t)\|_{L^2(\Omega)}=0,
		\end{equation}
		where the supremum is taken to be zero if $\tau_A\leq\delta$.
	\end{lemma}
	\begin{proof}
		\noindent\textbf{Step 1. Decomposition into $h^A$ and $z^A$
			and the energy estimate for $z^A$.}
		For each fixed $A$, the velocity field $Ar^{-1}e_\theta$ is
		bounded and divergence free. The standard Galerkin construction
		gives a unique energy solution for $L^2(\Omega)$ initial data
		and forcing in $L^2(0,\tau_A;H^{-1}(\Omega))$. In particular,
		let $h^A$ be the energy solution of
		\begin{equation}\label{eq hA}
			\begin{cases}
				\partial_t h^A+Ar^{-2}\partial_\theta h^A-\Delta h^A=0
				& \text{in}\quad(0,\tau_A)\times\Omega,\\
				h^A|_{r=1}=0
				& \text{for}\quad 0<t<\tau_A,\\
				h^A(0)=f^A(0)
				& \text{in}\quad\Omega.
			\end{cases}
		\end{equation}
		Set $z^A=f^A-h^A$. Subtracting \eqref{eq hA} from
		\eqref{eq:fast-rotation-scalar} gives
		\[
		\begin{cases}
			\partial_t z^A+Ar^{-2}\partial_\theta z^A-\Delta z^A
			=-\diver Q^A
			& \text{in}\quad(0,\tau_A)\times\Omega,\\
			z^A|_{r=1}=0
			& \text{for}\quad 0<t<\tau_A,\\
			z^A(0)=0
			& \text{in}\quad\Omega.
		\end{cases}
		\]
		Since $z^A\in L^2(0,\tau_A;H_0^1(\Omega))$ and
		$Q^A\in L^2(0,\tau_A;L^2(\Omega))$ by
		\eqref{eq:fast-rotation-bounds}, the equation for $z^A$ gives
		$\partial_t z^A\in L^2(0,\tau_A;H^{-1}(\Omega))$.
		Thus the equation can be tested with $z^A$. The transport
		term vanishes by incompressibility and approximation in
		$H_0^1(\Omega)$. Integration by parts on the right-hand side
		and Young's inequality give
		\[
		\begin{aligned}
			\frac12\frac{d}{dt}\|z^A\|_{L^2(\Omega)}^2
			+\|\nabla z^A\|_{L^2(\Omega)}^2
			=\int_\Omega Q^A\cdot\nabla z^A dx\leq \frac12\|Q^A\|_{L^2(\Omega)}^2
			+\frac12\|\nabla z^A\|_{L^2(\Omega)}^2.
		\end{aligned}
		\]
		Integrating in time and using $z^A(0)=0$, we obtain
		\begin{equation}\label{eq:fast-rotation-uniform-flux-map}
			\|z^A\|_{L^\infty(0,\tau_A;L^2(\Omega))}
			\leq \|Q^A\|_{L^2(0,\tau_A;L^2(\Omega))}.
		\end{equation}
		
		\noindent\textbf{Step 2. Decay of $z^A_{\neq}$.}
		For each fixed $\epsilon>0$, we construct an approximation
		$Q_\epsilon^A$ of $Q^A$. The energy estimate
		\eqref{eq:fast-rotation-uniform-flux-map} controls the
		difference between the two forced solutions by $\epsilon$.
		The bound on $r\diver Q_\epsilon^A$ then allows us to apply
		a Fourier mode estimate and bound the nonzero modes of the
		approximating solution by $K_\epsilon A^{-\frac13}$.
		
		\noindent\underline{\textbf{(i) Spatial approximation of $Q^A$.}}
		By \eqref{eq:fast-rotation-flux-tail}, choose $R_\epsilon>2$
		so that
		\begin{equation}\label{eq:fast-rotation-chosen-radius}
			\sup_{A\geq1}
			\|\mathbf{1}_{\{r>R_\epsilon\}}Q^A\|
			_{L^2(0,\tau_A;L^2(\Omega))}
			<\frac{\epsilon}{2}.
		\end{equation}
		On $\Omega_{R_\epsilon}$, let $P_m$ be the $L^2$-orthogonal
		projection onto the span of the first $m$ elements of an
		orthonormal basis of $L^2(\Omega_{R_\epsilon};\mathbb R^2)$
		consisting of smooth, compactly supported fields. Such a basis
		is obtained by applying Gram--Schmidt to an $L^2$-dense
		sequence in $C_c^\infty(\Omega_{R_\epsilon};\mathbb R^2)$
		and discarding linearly dependent elements.
		Completeness of the basis gives $P_mg\to g$ in
		$L^2(\Omega_{R_\epsilon})$ for each $g$. This convergence
		is uniform on the unit ball of $H^s(\Omega_{R_\epsilon})$.
		Indeed, for any $\eta>0$, the compact embedding
		$H^s(\Omega_{R_\epsilon})\hookrightarrow L^2(\Omega_{R_\epsilon})$
		gives finitely many $g_1,\ldots,g_N$ such that each $g$ in
		this ball satisfies $\|g-g_j\|_{L^2(\Omega_{R_\epsilon})}<\frac{\eta}{2}$
		for some $j$. For all sufficiently large $m$, one has
		$\max_{1\leq j\leq N}\|(I-P_m)g_j\|_{L^2(\Omega_{R_\epsilon})}<\frac{\eta}{2}$.
		Since $I-P_m$ is an $L^2$ contraction, there holds
		\[
		\|(I-P_m)g\|_{L^2(\Omega_{R_\epsilon})}
		\leq \|g-g_j\|_{L^2(\Omega_{R_\epsilon})}
		+\|(I-P_m)g_j\|_{L^2(\Omega_{R_\epsilon})}<\eta .
		\]
		The choice of $m$ is independent of $g$ in the unit ball.
		Taking the supremum over the ball and using that $\eta>0$
		is arbitrary, we obtain
		\begin{equation}\label{eq:fast-rotation-projection-limit}
			\lim_{m\to\infty}
			\sup_{\|g\|_{H^s(\Omega_{R_\epsilon})}\leq1}
			\|(I-P_m)g\|_{L^2(\Omega_{R_\epsilon})}=0.
		\end{equation}
		By \eqref{eq:fast-rotation-local-regularity} with $R=R_\epsilon$,
		choose $\eta>0$ such that
		\[
		\eta\sup_{A\geq1}
		\|Q^A\|_{L^2(0,\tau_A;H^s(\Omega_{R_\epsilon}))}
		<\frac{\epsilon}{2}.
		\]
		For this choice of $\eta$, take $m$ sufficiently large in
		\eqref{eq:fast-rotation-projection-limit} so that the projection
		error on the unit ball is at most $\eta$. Applying this bound
		to $g(\|g\|_{H^s(\Omega_{R_\epsilon})})^{-1}$ when $g\neq0$ gives
		\[
		\|(I-P_m)g\|_{L^2(\Omega_{R_\epsilon})}
		\leq \eta\|g\|_{H^s(\Omega_{R_\epsilon})},
		\qquad g\in H^s(\Omega_{R_\epsilon}),
		\]
		which also holds for $g=0$. Denote by
		$Q^A(t)|_{\Omega_{R_\epsilon}}$ the restriction of
		$Q^A(t,\cdot)$ to $\Omega_{R_\epsilon}$. Applying the last
		inequality to this restriction for almost every $t$ and
		taking the $L^2$ norm in time gives
		\begin{equation*}
			\|(I-P_m)(Q^A|_{\Omega_{R_\epsilon}})\|
			_{L^2(0,\tau_A;L^2(\Omega_{R_\epsilon}))}
			\leq \eta\|Q^A\|_{L^2(0,\tau_A;H^s(\Omega_{R_\epsilon}))}.
		\end{equation*}
		By the choice of $\eta$, taking the supremum over $A$ yields
		\begin{equation}\label{eq:fast-rotation-annular-approximation}
			\sup_{A\geq1}
			\|(I-P_m)(Q^A|_{\Omega_{R_\epsilon}})\|
			_{L^2(0,\tau_A;L^2(\Omega_{R_\epsilon}))}
			<\frac{\epsilon}{2}.
		\end{equation}
		In particular, the same $m$ works for every $A$.
		Define $Q_\epsilon^A$ as the zero extension of
		$P_m(Q^A|_{\Omega_{R_\epsilon}})$ to $\Omega$. The first $m$
		basis fields are compactly supported in $\Omega_{R_\epsilon}$,
		so their zero extensions are smooth on $\Omega$. The map
		$g\mapsto r\diver g$ is bounded in $L^2(\Omega)$ on the
		span of these finitely many fields. Together with the
		$L^2$ contraction property of $P_m$, this gives
		\begin{equation}\label{eq:fast-rotation-projected-flux}
			\begin{aligned}
				\|Q_\epsilon^A\|_{L^2(0,\tau_A;L^2(\Omega))}
				&\leq \|Q^A\|_{L^2(0,\tau_A;L^2(\Omega))},\\
				\|r\diver Q_\epsilon^A\|_{L^2(0,\tau_A;L^2(\Omega))}
				&\leq K_\epsilon\|Q^A\|_{L^2(0,\tau_A;L^2(\Omega))}.
			\end{aligned}
		\end{equation}
		The constant $K_\epsilon$ may depend on $R_\epsilon$ and $m$,
		which are fixed independently of $A$. The uniform bound
		\eqref{eq:fast-rotation-bounds} and the second inequality in
		\eqref{eq:fast-rotation-projected-flux} therefore control the
		weighted divergence. For the approximation error, use
		\eqref{eq:fast-rotation-annular-approximation} on
		$\Omega_{R_\epsilon}$ and \eqref{eq:fast-rotation-chosen-radius}
		on its complement, where $Q_\epsilon^A=0$. Thus
		\begin{equation}\label{eq:fast-rotation-spatial-approximation}
			\begin{aligned}
				\sup_{A\geq1}
				\|Q^A-Q_\epsilon^A\|_{L^2(0,\tau_A;L^2(\Omega))}
				<\epsilon,\qquad
				\sup_{A\geq1}
				\|r\diver Q_\epsilon^A\|_{L^2(0,\tau_A;L^2(\Omega))}
				\leq K_\epsilon.
			\end{aligned}
		\end{equation}
		
		\noindent\underline{\textbf{(ii) Decay for the approximating equation.}}
		The first inequality in \eqref{eq:fast-rotation-projected-flux}
		and \eqref{eq:fast-rotation-bounds} give
		$Q_\epsilon^A\in L^2(0,\tau_A;L^2(\Omega))$.
		Hence $-\diver Q_\epsilon^A\in L^2(0,\tau_A;H^{-1}(\Omega))$.
		The Galerkin argument used in Step~1 therefore gives a unique
		energy solution $z_\epsilon^A$ to
		\begin{equation}\label{eq:fast-rotation-approximate-forced-problem}
			\begin{cases}
				\partial_t z_\epsilon^A+Ar^{-2}\partial_\theta z_\epsilon^A
				-\Delta z_\epsilon^A=-\diver Q_\epsilon^A
				& \text{in}\quad(0,\tau_A)\times\Omega,\\
				z_\epsilon^A|_{r=1}=0 & \text{for}\quad 0<t<\tau_A,\\
				z_\epsilon^A(0)=0 & \text{in}~\quad\Omega.
			\end{cases}
		\end{equation}
		In particular,
		$z_\epsilon^A\in C([0,\tau_A);L^2(\Omega))
		\cap L^2(0,\tau_A;H_0^1(\Omega))$.
		Taking the $k$th angular Fourier coefficient of
		\eqref{eq:fast-rotation-approximate-forced-problem} and using
		$r^{\frac12}\widehat{(z_\epsilon^A)}_k$ as the radial unknown gives
		\begin{equation}\label{eq:fast-rotation-modal-equation}
			\left(\partial_t-\partial_r^2
			+\frac{k^2-\frac14+ikA}{r^2}\right)
			\bigl(r^{\frac12}\widehat{(z_\epsilon^A)}_k\bigr)
			=-r^{\frac12}\widehat{(\diver Q_\epsilon^A)}_k.
		\end{equation}
		The initial and boundary conditions in
		\eqref{eq:fast-rotation-approximate-forced-problem} imply
		\[
		\left.r^{\frac12}\widehat{(z_\epsilon^A)}_k\right|_{t=0}=0
		\quad\text{in }L^2(1,\infty),
		\qquad
		\left.r^{\frac12}\widehat{(z_\epsilon^A)}_k\right|_{r=1}=0
		\quad\text{for a.e. }t\in(0,\tau_A).
		\]
		To verify the radial regularity, note that
		\[
		\partial_r\bigl(r^{\frac12}\widehat{(z_\epsilon^A)}_k\bigr)
		=r^{\frac12}\partial_r\widehat{(z_\epsilon^A)}_k
		+\frac12r^{-\frac12}\widehat{(z_\epsilon^A)}_k.
		\]
		Using $r^{-1}\leq r$ and applying Parseval's identity to
		$z_\epsilon^A$ and $\partial_r z_\epsilon^A$, we obtain
		\begin{align*}
			\bigl\|r^{\frac12}\widehat{(z_\epsilon^A)}_k
			\bigr\|_{L^2(0,\tau_A;H^1(1,\infty))}^2
			&\leq K\int_0^{\tau_A}\int_1^\infty
			r\left(
			|\widehat{(z_\epsilon^A)}_k|^2
			+|\partial_r\widehat{(z_\epsilon^A)}_k|^2
			\right)drdt\\
			&\leq K\|z_\epsilon^A\|_{L^2(0,\tau_A;H^1(\Omega))}^2
			<\infty.
		\end{align*}
		Together with the zero trace at $r=1$, this yields
		$r^{\frac12}\widehat{(z_\epsilon^A)}_k
		\in L^2(0,\tau_A;H_0^1(1,\infty))$. The continuity of
		$z_\epsilon^A$ in $L^2(\Omega)$ and Parseval's identity also give
		$r^{\frac12}\widehat{(z_\epsilon^A)}_k
		\in C([0,\tau_A);L^2(1,\infty))$.
		For the right-hand side of \eqref{eq:fast-rotation-modal-equation},
		Parseval's identity and the second bound in
		\eqref{eq:fast-rotation-spatial-approximation} give
		\[
		\|r^{\frac32}\widehat{(\diver Q_\epsilon^A)}_k\|
		_{L^2(0,\tau_A;L^2(1,\infty))}
		\leq K\|r\diver Q_\epsilon^A\|_{L^2(0,\tau_A;L^2(\Omega))}
		\leq K_\epsilon.
		\]
		Since $|kA|\geq A$ for $k\neq0$, Lemma~2.5 of~\cite{LZZ2026}
		applies to \eqref{eq:fast-rotation-modal-equation} for all
		nonzero modes once $A$ is sufficiently large. Taking viscosity
		one, $f_1=-r^{\frac12}\widehat{(\diver Q_\epsilon^A)}_k$ and $f_2=0$
		in that lemma, we obtain
		\begin{equation}\label{eq:fast-rotation-modal-forcing}
			\begin{aligned}
				&\sup_{0\leq t<\tau_A}
				\int_1^\infty r|\widehat{(z_\epsilon^A)}_k(t,r)|^2 dr\\
				\leq~&
				\frac{K}{\max\{k^2,|kA|^{\frac23}\}}
				\int_0^{\tau_A}\int_1^\infty
				r^3|\widehat{(\diver Q_\epsilon^A)}_k(t,r)|^2 drdt,
				\qquad k\neq0,
			\end{aligned}
		\end{equation}
		where $K$ is independent of $A$, $k$, and $\tau_A$. The
		time continuity of the transformed mode allows the use of
		the supremum in place of the essential supremum.
		Using $\max\{k^2,|kA|^{\frac23}\}\geq A^{\frac23}$, we sum
		\eqref{eq:fast-rotation-modal-forcing} over $k\neq0$ and use
		Parseval's identity. Taking the square root and using the
		second bound in \eqref{eq:fast-rotation-spatial-approximation},
		we obtain
		\begin{equation}\label{eq:fast-rotation-approximate-forced-bound}
			\begin{aligned}
				\|(z_\epsilon^A)_{\neq}\|_{L^\infty(0,\tau_A;L^2(\Omega))}
				\leq K A^{-\frac13}
				\|r\diver Q_\epsilon^A\|_{L^2(0,\tau_A;L^2(\Omega))}\leq K_\epsilon A^{-\frac13}.
			\end{aligned}
		\end{equation}
		
		\noindent\underline{\textbf{(iii) Comparison with $z^A$.}}
		Subtracting \eqref{eq:fast-rotation-approximate-forced-problem}
		from the equation for $z^A$ derived in Step~1 from
		\eqref{eq:fast-rotation-scalar}, we find that $z^A-z_\epsilon^A$
		has zero initial and boundary data and forcing
		$-\diver(Q^A-Q_\epsilon^A)$.
		Applying \eqref{eq:fast-rotation-uniform-flux-map} to this
		difference and using the first bound in
		\eqref{eq:fast-rotation-spatial-approximation}, we obtain
		\[
		\|z^A-z_\epsilon^A\|_{L^\infty(0,\tau_A;L^2(\Omega))}
		\leq
		\|Q^A-Q_\epsilon^A\|_{L^2(0,\tau_A;L^2(\Omega))}<\epsilon.
		\]
		For any $g\in L^2(\Omega)$, write $g=g_0+g_{\neq}$, where
		$g_0$ is the angular mean of $g$. Since $g_0$ is independent
		of $\theta$ and $\int_0^{2\pi}g_{\neq}(r,\theta)d\theta=0$
		for almost every $r$, the two terms are orthogonal in $L^2(\Omega)$. Consequently,
		\begin{equation}\label{eq:nonzero-mode-L2}
			\|g\|_{L^2(\Omega)}^2
			=\|g_0\|_{L^2(\Omega)}^2+\|g_{\neq}\|_{L^2(\Omega)}^2,
			\qquad
			\|g_{\neq}\|_{L^2(\Omega)}\leq\|g\|_{L^2(\Omega)}.
		\end{equation}
		Applying \eqref{eq:nonzero-mode-L2} to $g=z^A(t)-z_\epsilon^A(t)$,
		taking the essential supremum in time, and using the preceding
		bound for $z^A-z_\epsilon^A$ give
		\[
		\|(z^A-z_\epsilon^A)_{\neq}\|_{L^\infty(0,\tau_A;L^2(\Omega))}
		\leq
		\|z^A-z_\epsilon^A\|_{L^\infty(0,\tau_A;L^2(\Omega))}<\epsilon.
		\]
		The triangle inequality for
		$z^A_{\neq}=(z^A-z_\epsilon^A)_{\neq}+(z_\epsilon^A)_{\neq}$
		and \eqref{eq:fast-rotation-approximate-forced-bound} now give
		\[
		\|z^A_{\neq}\|_{L^\infty(0,\tau_A;L^2(\Omega))}
		\leq \epsilon+K_\epsilon A^{-\frac13}.
		\]
		Taking the upper limit as $A\to\infty$ for fixed $\epsilon$,
		and then letting $\epsilon\to 0$, proves
		\begin{equation}\label{eq:fast-rotation-forced-limit}
			\lim_{A\to\infty}
			\|z^A_{\neq}\|_{L^\infty(0,\tau_A;L^2(\Omega))}=0.
		\end{equation}
		
		\noindent\textbf{Step 3. Proof of \eqref{eq:fast-rotation-limit}.}
		Fix $0<\delta<T$ and $R>1$. By the convention in the statement of the
		lemma, the supremum in \eqref{eq:fast-rotation-limit} is zero
		when $\tau_A\leq\delta$. It therefore suffices to consider
		$\tau_A>\delta$. The argument combines the decay of
		$h^A_{\neq}$ on $\Omega_R$ with a bound for $f^A_{\neq}$ on
		$\{r>R\}$.
		
		Taking the $k$th angular Fourier coefficient of \eqref{eq hA} gives
		\[
		\left(\partial_t-\partial_r^2-\frac1r\partial_r
		+\frac{k^2+ikA}{r^2}\right)\widehat h^A_k=0,\qquad
		\left.\widehat h^A_k\right|_{r=1}=0,
		\qquad
		\left.\widehat h^A_k\right|_{t=0}
		=\widehat f^A_k(0,r).
		\]
		The same change of unknown as in
		\eqref{eq:fast-rotation-modal-equation} yields
		\[
		\begin{aligned}
			\left(\partial_t-\partial_r^2
			+\frac{k^2-\frac14+ikA}{r^2}\right)
			\bigl(r^{\frac12}\widehat h^A_k\bigr)=0,\quad
			\left.r^{\frac12}\widehat h^A_k\right|_{r=1}=0,
			\quad
			\left.r^{\frac12}\widehat h^A_k\right|_{t=0}
			=r^{\frac12}\widehat f^A_k(0,r).
		\end{aligned}
		\]
		The transformed initial data belong to $L^2(1,\infty)$ by
		Parseval's identity and the initial-data bound in
		\eqref{eq:fast-rotation-bounds}. As for $z_\epsilon^A$ in Step~2,
		the energy regularity of $h^A$ gives
		$r^{\frac12}\widehat h^A_k\in C([0,\tau_A);L^2(1,\infty))
		\cap L^2(0,\tau_A;H_0^1(1,\infty))$.
		Since $|kA|\geq A$ for
		$k\neq0$, the condition $|kA|^{-1}\ll 2^{-3}$ holds for all
		nonzero modes once $A$ is sufficiently large. Theorem~1.1
		of~\cite{LZZ2026}, with $q=1$ and viscosity one, then gives
		\begin{equation}\label{eq:fast-rotation-homogeneous-modal-bound}
			\begin{aligned}
				\sup_{0\leq t<\tau_A}
				\int_1^\infty r
				\left(1+\frac{|kA|^{\frac23}t}{r^2}\right)^2
				|\widehat h^A_k(t,r)|^2 dr\leq K\int_1^\infty
				r|\widehat f^A_k(0,r)|^2 dr,
				\qquad k\neq0,
			\end{aligned}
		\end{equation}
		where $K$ is independent of $A$, $k$, and $\tau_A$. For
		fixed $A$ and $k$, the time-dependent weight is continuous
		in $L^\infty(1,\infty)$. The weighted transformed mode is
		therefore continuous in $L^2(1,\infty)$, which justifies
		the supremum over time in
		\eqref{eq:fast-rotation-homogeneous-modal-bound}.
		On $\Omega_R$, the weight in
		\eqref{eq:fast-rotation-homogeneous-modal-bound} satisfies
		\[
		1+\frac{|kA|^{\frac23}t}{r^2}
		\geq 1+\frac{A^{\frac23}\delta}{R^2}
		\qquad\text{for}\quad
		\delta\leq t<\tau_A,\ 1<r<R,\ k\neq0.
		\]
		Restricting the radial integral in
		\eqref{eq:fast-rotation-homogeneous-modal-bound} to $(1,R)$,
		summing over $k\neq0$, and applying Parseval's identity give
		\begin{equation}\label{eq:fast-rotation-homogeneous-local-bound}
			\sup_{\delta\leq t<\tau_A}
			\|h^A_{\neq}(t)\|_{L^2(\Omega_R)}
			\leq
			\frac{K\|f^A(0)\|_{L^2(\Omega)}}
			{1+A^{\frac23}\delta R^{-2}}.
		\end{equation}
		This gives decay on $\Omega_R$ for each fixed $R$. Outside
		$\Omega_R$, the weighted bound in \eqref{eq:fast-rotation-bounds}
		controls $f^A_{\neq}$ directly. Indeed, since
		$\mathbf{1}_{\{r>R\}}$ is independent of $\theta$, we have
		$(\mathbf{1}_{\{r>R\}}f^A)_{\neq}
		=\mathbf{1}_{\{r>R\}}f^A_{\neq}$. Applying
		\eqref{eq:nonzero-mode-L2} to $g=\mathbf{1}_{\{r>R\}}f^A(t)$
		and using $r^{-2}\leq R^{-2}$ on $\{r>R\}$ gives
		\[
		\|\mathbf{1}_{\{r>R\}}f^A_{\neq}(t)\|_{L^2(\Omega)}
		\leq \|\mathbf{1}_{\{r>R\}}f^A(t)\|_{L^2(\Omega)}
		\leq R^{-2}\|r^2f^A(t)\|_{L^2(\Omega)}
		\]
		for almost every $t\in(0,\tau_A)$.
		Since $f^A\in C([0,\tau_A);L^2(\Omega))$, the norm on the
		left is continuous in time. Its bound by the essential
		supremum of the right-hand side therefore holds at every
		$t\in[0,\tau_A)$:
		\begin{equation}\label{eq:fast-rotation-scalar-tail}
			\sup_{0\leq t<\tau_A}
			\|\mathbf{1}_{\{r>R\}}f^A_{\neq}(t)\|_{L^2(\Omega)}
			\leq R^{-2}\|r^2f^A\|_{L^\infty(0,\tau_A;L^2(\Omega))}.
		\end{equation}
		On $\Omega_R$, the decomposition
		$f^A_{\neq}=h^A_{\neq}+z^A_{\neq}$ and
		\eqref{eq:fast-rotation-homogeneous-local-bound} control
		the interior norm. The time continuity of $z^A=f^A-h^A$ in
		$L^2(\Omega)$ implies that the $L^\infty(0,\tau_A;L^2(\Omega))$
		bound for $z^A_{\neq}$ holds at every time. On $\{r>R\}$,
		\eqref{eq:fast-rotation-scalar-tail} controls the exterior
		norm. Combining these estimates and using
		\eqref{eq:fast-rotation-bounds} gives
		\[
		\sup_{\delta\leq t<\tau_A}\|f^A_{\neq}(t)\|_{L^2(\Omega)}
		\leq \frac{K}{1+A^{\frac23}\delta R^{-2}}
		+\|z^A_{\neq}\|_{L^\infty(0,\tau_A;L^2(\Omega))}
		+KR^{-2},
		\]
		where $K$ is independent of $A$ and $R$. For fixed $R$, the
		first term tends to zero as $A\to\infty$, and the second
		does so by \eqref{eq:fast-rotation-forced-limit}. Thus
		\[
		\limsup_{A\to\infty}
		\sup_{\delta\leq t<\tau_A}\|f^A_{\neq}(t)\|_{L^2(\Omega)}
		\leq KR^{-2}.
		\]
		Letting $R\to\infty$ proves \eqref{eq:fast-rotation-limit}.
	\end{proof}
	\subsection{Proof of Proposition~\ref{prop:finite-time-averaging}}
	The bounds in Lemma~\ref{lem:finite-time-uniform-moments}
	and the decay in Lemma~\ref{lem:finite-time-fast-rotation}
	are now used in an energy comparison of $n_0^A$ with $N$.
	The comparison is carried out while
	$\|n^A(t)\|_{L^2(\Omega)}<B$, and then improves this bound
	to extend the solution to $[0,T]$.
	\begin{proof}[Proof of Proposition~\ref{prop:finite-time-averaging}]
		Fix $T>0$ and take $B$ as in \eqref{eq:bootstrap-level}.
		Let $[0,\tau_A)$ be the maximal initial interval, with
		$\tau_A\leq T$, on which the strong solution exists and
		$\|n^A(t)\|_{L^2(\Omega)}<B$. The initial bound
		$\|n_{\mathrm{in}}\|_{L^2(\Omega)}<\frac{B}{4}$ and
		Proposition~\ref{prop:local-strong-theory} give $\tau_A>0$.
		Throughout the proof, $K_T$ may depend on $T$, $B$, and
		the initial data, but is independent of $A$, $\tau_A$,
		and $\delta$.
		
		\noindent\textbf{Step 1. Estimates for $n_{\neq}^A$ and the
			additional term in the equation for $n_0^A$.}
		On $[0,\tau_A)$, the defining bound gives
		\eqref{eq:finite-time-L2-bootstrap} with $T_*=\tau_A$, and
		Proposition~\ref{prop:local-strong-theory} gives $n^A\geq0$.
		Hence Lemma~\ref{lem:finite-time-uniform-moments} applies.
		By $\diver u^A=0$, the density equation in
		\eqref{eq:physical-scalars-moments} has the form
		\eqref{eq:fast-rotation-scalar} with
		\[
		(f^A,Q^A)=(n^A,u^A n^A+n^A\nabla c^A).
		\]
		For this choice of $(f^A,Q^A)$, the boundary condition
		\eqref{eq:main-boundary}, the strong solution regularity, and
		the estimates \eqref{eq:finite-time-L2-bootstrap} and
		\eqref{eq:finite-time-velocity-bounds} give
		\[
		\left.f^A\right|_{r=1}=0,\qquad
		f^A\in C([0,\tau_A);L^2(\Omega))
		\cap L^2(0,\tau_A;H_0^1(\Omega)).
		\]
		Since the initial data are independent of $A$, the weighted
		estimates in \eqref{eq:finite-time-velocity-bounds} and the
		two bounds in \eqref{eq:finite-time-flux-regularity} yield,
		for each fixed $0<s<1$,
		\begin{align*}
			&\sup_{A\geq1}\Bigl(
			\|f^A(0)\|_{L^2(\Omega)}
			+\|r^2 f^A\|_{L^\infty(0,\tau_A;L^2(\Omega))}
			+\|Q^A\|_{L^2(0,\tau_A;L^2(\Omega))}
			\Bigr)\leq K_{T,s},\\
			&\sup_{A\geq1}
			\|Q^A\|_{L^2(0,\tau_A;H^s(\Omega_R))}<\infty,
			\qquad R>1,\\
			&\sup_{A\geq1}
			\|\mathbf{1}_{\{r>R\}}Q^A\|_{L^2(0,\tau_A;L^2(\Omega))}
			\leq K_TR^{-2},
			\qquad R\geq2.
		\end{align*}
		These estimates verify \eqref{eq:fast-rotation-bounds},
		\eqref{eq:fast-rotation-local-regularity}, and
		\eqref{eq:fast-rotation-flux-tail}, respectively.
		Applying Lemma~\ref{lem:finite-time-fast-rotation} to $n^A$
		therefore gives, for every $0<\delta<T$,
		\begin{equation}\label{eq:varying-interval-mixing}
			\sup_{\delta\leq t<\tau_A}
			\|n^A_{\neq}(t)\|_{L^2(\Omega)}
			\to 0\quad\text{as}\quad A\to\infty,
		\end{equation}
		with the supremum taken as zero if $\tau_A\leq\delta$.
		
		To derive the equation for the angular mean, use
		$(\diver F)_0=r^{-1}\partial_r(r(F\cdot e_r)_0)$ for a
		vector field $F$. Thus $\diver u^A=0$ gives
		$\partial_r(r(u_r^A)_0)=0$. Since $u_r^A=0$ at $r=1$,
		it follows that $(u_r^A)_0=0$. The mixed terms in the angular
		average of $n^A\partial_r c^A$ vanish because
		$(n^A_{\neq})_0=(\partial_r c^A_{\neq})_0=0$. In particular,
		\[
		(u_r^A n^A)_0=(u_r^A n^A_{\neq})_0,\qquad
		(n^A\partial_r c^A)_0
		=n_0^A\partial_r c_0^A+(n^A_{\neq}\partial_r c^A_{\neq})_0.
		\]
		Averaging \eqref{eq:physical-scalars-moments} therefore gives
		\begin{equation}\label{eq:exact-mean-equation}
			\begin{aligned}
				\partial_t n^A_0
				&=\frac1r\partial_r\!\left[
				r\left\{\partial_r n^A_0-n^A_0\partial_r c^A_0
				-(n^A_{\neq}\partial_r c^A_{\neq}+u_r^A n^A_{\neq})_0
				\right\}\right],\\&
				(1-\partial_r^2-r^{-1}\partial_r)c^A_0=n^A_0,
			\end{aligned}
		\end{equation}
		with $n_0^A=c_0^A=0$ at $r=1$.
		The sum of the averaged products in
		\eqref{eq:exact-mean-equation} satisfies
		\begin{equation}\label{eq:interaction-flux-vanishes}
			\left(\int_0^{\tau_A}\int_1^\infty r
			\left|(n^A_{\neq}\partial_r c^A_{\neq}+u_r^A n^A_{\neq})_0
			\right|^2 drdt\right)^{\frac12}
			\to 0\qquad\text{as}\quad A\to\infty.
		\end{equation}
		To prove \eqref{eq:interaction-flux-vanishes}, fix $0<\delta<T$
		and first estimate the contribution from
		$(0,\min\{\delta,\tau_A\})$. By \eqref{eq:nonzero-mode-L2}
		and \eqref{eq:finite-time-L2-bootstrap}, we get
		\[
		\|n^A_{\neq}\|_{L^2(\Omega)}
		\leq \|n^A\|_{L^2(\Omega)}\leq B.
		\]
		Since angular averaging does not increase the $L^\infty$ norm
		and commutes with $\partial_r$, one also has
		\begin{align*}
			\|n^A_{\neq}\|_{L^\infty(\Omega)}
			&\leq \|n^A\|_{L^\infty(\Omega)}
			+\|n^A_0\|_{L^\infty(\Omega)}
			\leq 2\|n^A\|_{L^\infty(\Omega)},\\
			\|\partial_r c^A_{\neq}\|_{L^\infty(\Omega)}
			&\leq \|\partial_r c^A\|_{L^\infty(\Omega)}
			+\|(\partial_r c^A)_0\|_{L^\infty(\Omega)}
			\leq 2\|\nabla c^A\|_{L^\infty(\Omega)}.
		\end{align*}
		Jensen's inequality in $\theta$ and H\"older's inequality in
		space therefore give, for almost every $t\in(0,\tau_A)$,
		\begin{equation}\label{eq:interaction-flux-uniform-bound}
			\begin{aligned}
				&\left(\int_1^\infty r
				\left|(n^A_{\neq}\partial_r c^A_{\neq}
				+u_r^A n^A_{\neq})_0\right|^2 dr\right)^{\frac12}
				\leq (2\pi)^{-\frac12}
				\|n^A_{\neq}\partial_r c^A_{\neq}
				+u_r^A n^A_{\neq}\|_{L^2(\Omega)}
				\\&\leq (2\pi)^{-1/2}\Bigl(
				\|n^A_{\neq}\|_{L^2(\Omega)}
				\|\partial_r c^A_{\neq}\|_{L^\infty(\Omega)}
				+\|u^A\|_{L^2(\Omega)}
				\|n^A_{\neq}\|_{L^\infty(\Omega)}\Bigr)\\
				&\leq (2\pi)^{-\frac12}\Bigl(
				2B\|\nabla c^A\|_{L^\infty(\Omega)}
				+2\|u^A\|_{L^2(\Omega)}
				\|n^A\|_{L^\infty(\Omega)}\Bigr)
				\leq K_T,
			\end{aligned}
		\end{equation}
		where the last inequality follows from
		\eqref{eq:finite-time-velocity-bounds} and
		\eqref{eq:finite-time-chemical-gradient}.
		Squaring \eqref{eq:interaction-flux-uniform-bound} and integrating
		over $(0,\min\{\delta,\tau_A\})$ yields
		\begin{equation}\label{0 delta}
			\left(\int_0^{\min\{\delta,\tau_A\}}\int_1^\infty r
			\left|(n^A_{\neq}\partial_r c^A_{\neq}
			+u_r^A n^A_{\neq})_0\right|^2 drdt\right)^{\frac12}
			\leq K_T\sqrt{\min\{\delta,\tau_A\}}
			\leq K_T\sqrt{\delta}.
		\end{equation}
		Suppose that $\tau_A>\delta$. On $(\delta,\tau_A)$,
		the decay estimate \eqref{eq:varying-interval-mixing}
		allows both products to be estimated with $n^A_{\neq}$ in
		$L^2(\Omega)$. For almost every $t\in(\delta,\tau_A)$,
		Jensen's inequality in $\theta$ and H\"older's inequality
		in space give
		\begin{align*}
			&\left(\int_1^\infty r
			\left|(n^A_{\neq}\partial_r c^A_{\neq}
			+u_r^A n^A_{\neq})_0\right|^2 dr\right)^{\frac12}\leq (2\pi)^{-\frac12}
			\|n^A_{\neq}\partial_r c^A_{\neq}
			+u_r^A n^A_{\neq}\|_{L^2(\Omega)}\\&\leq K\|n^A_{\neq}\|_{L^2(\Omega)}
			\bigl(\|\nabla c^A\|_{L^\infty(\Omega)}
			+\|u^A\|_{L^\infty(\Omega)}\bigr),
		\end{align*}
		where the last inequality uses
		$\|\partial_r c^A_{\neq}\|_{L^\infty(\Omega)}
		\leq 2\|\nabla c^A\|_{L^\infty(\Omega)}$
		and $|u_r^A|\leq |u^A|$.
		Taking the $L^2$ norm in time and using $\tau_A-\delta\leq T$
		then yield
		\begin{align*}
			&\left(\int_\delta^{\tau_A}\int_1^\infty r
			\left|(n^A_{\neq}\partial_r c^A_{\neq}
			+u_r^A n^A_{\neq})_0\right|^2\,dr\,dt\right)^{\frac12}\\
			&\leq K\|n^A_{\neq}\|_{L^\infty(\delta,\tau_A;L^2(\Omega))}
			\Bigl(\sqrt{T}\|\nabla c^A\|_{L^\infty(0,\tau_A;L^\infty(\Omega))}
			+\|u^A\|_{L^2(0,\tau_A;L^\infty(\Omega))}\Bigr)\\&\leq K_T\|n^A_{\neq}\|_{L^\infty(\delta,\tau_A;L^2(\Omega))},
		\end{align*}
		by \eqref{eq:finite-time-chemical-gradient} and
		\eqref{eq:finite-time-velocity-bounds}. For each fixed $0<\delta<T$, combining the above estimate with \eqref{0 delta} gives
		\begin{align*}
			\left(\int_0^{\tau_A}\int_1^\infty r
			\left|(n^A_{\neq}\partial_r c^A_{\neq}
			+u_r^A n^A_{\neq})_0\right|^2 drdt\right)^{\frac12}\leq K_T\sqrt{\delta}
			+K_T\sup_{\delta\leq t<\tau_A}
			\|n^A_{\neq}(t)\|_{L^2(\Omega)}.
		\end{align*}
		When $\tau_A\leq\delta$, the first estimate covers the whole
		interval, and the supremum is taken to be zero as in
		\eqref{eq:varying-interval-mixing}. For fixed $\delta$, the
		second term tends to zero as $A\to\infty$ by
		\eqref{eq:varying-interval-mixing}. Hence
		\[
		\limsup_{A\to\infty}
		\left(\int_0^{\tau_A}\int_1^\infty r
		\left|(n^A_{\neq}\partial_r c^A_{\neq}
		+u_r^A n^A_{\neq})_0\right|^2 drdt\right)^{\frac12}
		\leq K_T\sqrt{\delta}.
		\]
		Since $K_T$ is independent of $\delta$, letting
		$\delta\to 0$ proves \eqref{eq:interaction-flux-vanishes}.
		
		\noindent\textbf{Step 2. Estimate for $n_0^A-N$ on $[0,\tau_A)$.}
		To prove \eqref{eq:finite-time-averaging}, first establish the
		corresponding estimate on $[0,\tau_A)$. Let $(N,C)$
		be the global radial solution of \eqref{eq:radial-reference-problem}
		given by Proposition~\ref{prop:radial-reference-decay}.
		Subtracting the equation for $C$ from the equation for $c^A_0$
		in \eqref{eq:exact-mean-equation} gives
		\[
		(1-\partial_r^2-r^{-1}\partial_r)(c^A_0-C)=n^A_0-N,
		\qquad (c^A_0-C)|_{r=1}=0.
		\]
		Since $c^A_0-C$ is radial and belongs to $H_0^1(\Omega)$,
		the second bound in \eqref{eq:radial-elliptic-estimates} yields
		\begin{equation}\label{eq:mean-chemical-difference}
			\|\partial_r(c^A_0-C)(t)\|_{L^\infty(1,\infty)}
			\leq \sqrt{6}\left(\int_1^\infty r
			|n^A_0(t,r)-N(t,r)|^2 dr\right)^{\frac12}.
		\end{equation}
		Subtracting the equation for $N$ in
		\eqref{eq:radial-reference-problem} from the equation for $n^A_0$
		in \eqref{eq:exact-mean-equation}
		and writing
		$n^A_0\partial_r c^A_0-N\partial_r C
		=(n^A_0-N)\partial_r c^A_0+N\partial_r(c^A_0-C)$
		give
		\begin{equation}\label{eq:mean-density-difference}
			\begin{aligned}
				&\partial_t(n^A_0-N)
				-\frac1r\partial_r\bigl(r\partial_r(n^A_0-N)\bigr)\\
				=~&-\frac1r\partial_r\biggl[r\Bigl(
				(n^A_0-N)\partial_r c^A_0+N\partial_r(c^A_0-C)
				+(n^A_{\neq}\partial_r c^A_{\neq}
				+u_r^A n^A_{\neq})_0\Bigr)\biggr].
			\end{aligned}
		\end{equation}
		Since $\|\partial_r c^A_0\|_{L^\infty(1,\infty)}
		\leq \|\nabla c^A\|_{L^\infty(\Omega)}$, the bounds \eqref{eq:radial-reference-L2-bound},
		\eqref{eq:finite-time-chemical-gradient},
		and
		\eqref{eq:mean-chemical-difference} give
		\begin{align}\label{eq:mean-comparison-chemotaxis-bound}
			&\left(\int_1^\infty r
			\left|(n^A_0-N)\partial_r c^A_0
			+N\partial_r(c^A_0-C)\right|^2 dr\right)^{\frac12}
			\nonumber\\\leq~& \|\nabla c^A\|_{L^\infty(\Omega)}
			\left(\int_1^\infty r|n^A_0-N|^2 dr\right)^{\frac12}+\|\partial_r(c^A_0-C)\|_{L^\infty(1,\infty)}
			\left(\int_1^\infty r|N|^2 dr\right)^{\frac12}\nonumber\\\leq~& K_T
			\left(\int_1^\infty r|n^A_0-N|^2 dr\right)^{\frac12}.
		\end{align}
		Angular averaging preserves $H_0^1(\Omega)$.
		Thus \eqref{eq:finite-time-L2-bootstrap}, the gradient
		bound in \eqref{eq:finite-time-velocity-bounds},
		and the strong regularity of $N$ give
		$n_0^A-N\in L^2(0,\tau_A;H_0^1(\Omega))$.
		Moreover, \eqref{eq:mean-comparison-chemotaxis-bound}
		and \eqref{eq:interaction-flux-uniform-bound} imply
		\[
		\begin{aligned}
			(n_0^A-N)\partial_r c_0^A
			+N\partial_r(c_0^A-C)+(n_{\neq}^A\partial_r c_{\neq}^A
			+u_r^A n_{\neq}^A)_0
			\in L^2(0,\tau_A;L^2(\Omega)),
		\end{aligned}
		\]
		where each radial function is viewed as a function on $\Omega$.
		Since multiplication by $e_r$ preserves the $L^2(\Omega)$
		norm, the difference equation
		\eqref{eq:mean-density-difference} therefore gives
		$\partial_t(n_0^A-N)\in L^2(0,\tau_A;H^{-1}(\Omega))$.
		The $L^2(\Omega)$ energy chain rule consequently permits
		testing \eqref{eq:mean-density-difference} with $n_0^A-N$.
		Dividing the resulting identity by $2\pi$, we obtain,
		for almost every $t\in(0,\tau_A)$,
		\begin{align*}
			&\frac12\frac{d}{dt}\int_1^\infty r|n_0^A-N|^2 dr
			+\int_1^\infty r|\partial_r(n_0^A-N)|^2 dr\\
			=~&\int_1^\infty r\partial_r(n_0^A-N)
			\bigl((n_0^A-N)\partial_r c_0^A
			+N\partial_r(c_0^A-C)\bigr)dr\\
			&+\int_1^\infty r\partial_r(n_0^A-N)
			(n_{\neq}^A\partial_r c_{\neq}^A
			+u_r^A n_{\neq}^A)_0 dr.
		\end{align*}
		For the first integral on the right-hand side of the energy
		identity, Cauchy--Schwarz and
		\eqref{eq:mean-comparison-chemotaxis-bound} give
		\begin{align*}
			&\left|\int_1^\infty r\partial_r(n^A_0-N)
			\bigl((n^A_0-N)\partial_r c^A_0
			+N\partial_r(c^A_0-C)\bigr)dr\right|\\
			\leq~& K_T
			\left(\int_1^\infty r|\partial_r(n^A_0-N)|^2 dr\right)^{\frac12}
			\left(\int_1^\infty r|n^A_0-N|^2 dr\right)^{\frac12}\\
			\leq~& \frac14\int_1^\infty r|\partial_r(n^A_0-N)|^2 dr
			+K_T\int_1^\infty r|n^A_0-N|^2 dr.
		\end{align*}
		For the second integral, by Cauchy--Schwarz and Young's inequality, we obtain
		\begin{align*}
			&\left|\int_1^\infty r\partial_r(n^A_0-N)
			(n^A_{\neq}\partial_r c^A_{\neq}+u_r^A n^A_{\neq})_0
			dr\right|\\
			\leq~&
			\left(\int_1^\infty r|\partial_r(n^A_0-N)|^2 dr\right)^{\frac12}
			\left(\int_1^\infty r
			\left|(n^A_{\neq}\partial_r c^A_{\neq}
			+u_r^A n^A_{\neq})_0\right|^2 dr\right)^{\frac12}\\
			\leq~&\frac14\int_1^\infty r|\partial_r(n^A_0-N)|^2 dr
			+\int_1^\infty r
			\left|(n^A_{\neq}\partial_r c^A_{\neq}
			+u_r^A n^A_{\neq})_0\right|^2 dr.
		\end{align*}
		Substituting these two bounds into the energy identity gives
		\begin{align*}
			&\frac{d}{dt}\int_1^\infty r|n^A_0-N|^2 dr
			+\int_1^\infty r|\partial_r(n^A_0-N)|^2 dr\\
			\leq~& K_T\int_1^\infty r|n^A_0-N|^2 dr
			+K\int_1^\infty r
			\left|(n^A_{\neq}\partial_r c^A_{\neq}
			+u_r^A n^A_{\neq})_0\right|^2 dr.
		\end{align*}
		Since $n^A_0(0)=N(0)$ and $\tau_A\leq T$, Gronwall's
		inequality and \eqref{eq:interaction-flux-vanishes} imply
		\begin{equation}\label{eq:mean-convergence-varying}
			\begin{split}
				&\sup_{0\leq t<\tau_A}
				\left(\int_1^\infty r|n^A_0(t,r)-N(t,r)|^2 dr\right)^{\frac12}\\
				\leq~& K_T\left(\int_0^{\tau_A}\int_1^\infty r
				\left|(n^A_{\neq}\partial_r c^A_{\neq}
				+u_r^A n^A_{\neq})_0\right|^2 drdt\right)^{\frac12}
				\to 0\qquad\text{as}\quad A\to\infty.
			\end{split}
		\end{equation}
		
		\noindent\textbf{Step 3. Continuation to $T$ and proof of
			\eqref{eq:finite-time-averaging} and \eqref{n omega nonzero}.}
		The first inequality in \eqref{eq:finite-time-basic-energy-inequalities}
		gives
		\[
		\frac{d}{dt}\|n^A(t)\|_{L^2(\Omega)}^2\leq K_TB^2
		\qquad\text{for almost every }t\in(0,\tau_A).
		\]
		Set $\delta_0=\min\left\{\frac{T}{2},[16(1+K_T)]^{-1}\right\}$. Since
		$\|n_{\mathrm{in}}\|_{L^2(\Omega)}<\frac{B}{4}$ by
		\eqref{eq:bootstrap-level}, integration yields
		\[
		\|n^A(t)\|_{L^2(\Omega)}^2
		\leq \|n_{\mathrm{in}}\|_{L^2(\Omega)}^2+K_TB^2t
		\leq \frac{B^2}{8}\quad\text{for}\quad
		0\leq t<\min\{\delta_0,\tau_A\}.
		\]
		For $\delta_0\leq t<\tau_A$, using the decomposition
		$n^A=N+(n_0^A-N)+n_{\neq}^A$, we get
		\begin{align*}
			\|n^A(t)\|_{L^2(\Omega)}
			\leq \|N(t)\|_{L^2(\Omega)}
			+\sqrt{2\pi}\left(\int_1^\infty r
			|n_0^A(t,r)-N(t,r)|^2 dr\right)^{\frac12}+\|n_{\neq}^A(t)\|_{L^2(\Omega)}.
		\end{align*}
		The first term is less than $\frac{B}{4}$ by \eqref{eq:bootstrap-level}.
		By \eqref{eq:mean-convergence-varying} and
		\eqref{eq:varying-interval-mixing} with $\delta=\delta_0$,
		there exists $A_T\geq1$ such that the sum of the last two terms
		is at most $\frac{B}{4}$ throughout this interval for $A\geq A_T$.
		Together with the short-time estimate, this gives
		\begin{equation}\label{eq:improved-finite-bootstrap}
			\sup_{0\leq t<\tau_A}\|n^A(t)\|_{L^2(\Omega)}
			\leq \frac B2\quad\text{for}\quad
			A\geq A_T.
		\end{equation}
		Fix $A\geq A_T$. By \eqref{eq:improved-finite-bootstrap}
		and Proposition~\ref{prop:density-continuation}, $\tau_A$
		cannot be the maximal existence time of the strong solution.
		Thus the solution is defined at $\tau_A$, and continuity gives
		\[
		\|n^A(\tau_A)\|_{L^2(\Omega)}\leq \frac B2<B.
		\]
		If $\tau_A<T$, the bound $\|n^A(t)\|_{L^2(\Omega)}<B$ would
		persist beyond $\tau_A$, contradicting its definition.
		Consequently, $\tau_A=T$ and
		\[
		\sup_{0\leq t\leq T}\|n^A(t)\|_{L^2(\Omega)}
		\leq \frac B2.
		\]
		With $\tau_A=T$, the $L^2(\Omega)$ continuity of $n^A$
		and $N$ allows $t=T$ to be included in the suprema in
		\eqref{eq:mean-convergence-varying} and
		\eqref{eq:varying-interval-mixing}. These estimates yield
		\eqref{eq:finite-time-averaging} and \eqref{n omega nonzero},
		respectively.
	\end{proof}
	
	\appendix
	\section{Local existence and continuation}
	\label{app:local}
	The local construction and the continuation argument
	are given in the rotating coordinates of
	Section~\ref{sec:model}. Their estimates may depend on the
	fixed rotation amplitude $A$.
	\subsection{Local existence and uniqueness of strong solutions}
	The proof of Proposition~\ref{prop:local-strong-theory}
	combines linear estimates for the Dirichlet Laplacian
	$A_D$ and the Stokes operator $S$ with a contraction
	argument. These operators were introduced after
	\eqref{eq:local-rotated-scalar}--\eqref{eq:local-rotated-velocity}.
	\begin{proof}
		\noindent{\textbf{Step 1. Estimates for $A_D$, $S$ and $\mathbb P$.}}
		The linear estimates in Step~2 require an $H^1(\Omega)$ bound
		for $\mathbb P$ and spatial regularity estimates for $S$ and $A_D$.
		These estimates are established below.
		
		\noindent\underline{\textbf{(i) Estimate for $\mathbb P$.}}
		Let
		\[
		\begin{gathered}
			\mathcal H=\left\{v\in L^2(\Omega;\mathbb R^2):
			\int_\Omega v\cdot\nabla\psi dx=0
			\text{ for every }\psi\in H^1(\Omega)\right\},\\
			\mathcal V=\mathcal H\cap H^1(\Omega;\mathbb R^2).
		\end{gathered}
		\]
		To verify that $\mathcal H$ is closed, suppose that
		$v_j\in\mathcal H$ and $v_j\to v$ in
		$L^2(\Omega;\mathbb R^2)$ as $j\to\infty$.
		For every $\psi\in H^1(\Omega)$, the definition of $\mathcal H$
		and the Cauchy--Schwarz inequality give
		\[
		\begin{aligned}
			\left|\int_\Omega v\cdot\nabla\psi\,dx\right|
			=\left|\int_\Omega(v-v_j)\cdot\nabla\psi dx\right|\leq\|v-v_j\|_{L^2(\Omega)}
			\|\nabla\psi\|_{L^2(\Omega)}\to 0.
		\end{aligned}
		\]
		Thus $v\in\mathcal H$. The orthogonal projection theorem
		therefore defines the $L^2(\Omega)$ projection $\mathbb P$
		onto $\mathcal H$.
		For every $F\in L^2(\Omega;\mathbb R^2)$, we have
		\[
		\mathbb PF\in\mathcal H,\qquad
		\|F-\mathbb PF\|_{L^2(\Omega)}
		=\min_{z\in\mathcal H}\|F-z\|_{L^2(\Omega)}.
		\]
		If $v\in\mathcal H$, then $\mathbb PF+\lambda v\in\mathcal H$
		for every $\lambda\in\mathbb R$. Hence
		\[
		\begin{aligned}
			0\leq
			\|F-\mathbb PF-\lambda v\|_{L^2(\Omega)}^2
			-\|F-\mathbb PF\|_{L^2(\Omega)}^2=-2\lambda\int_\Omega(F-\mathbb PF)\cdot v dx
			+\lambda^2\|v\|_{L^2(\Omega)}^2.
		\end{aligned}
		\]
		Differentiating at the minimum $\lambda=0$ gives
		\begin{equation}\label{eq:local-projection-orthogonality}
			\int_\Omega(F-\mathbb PF)\cdot v dx=0
			\qquad\text{for every }v\in\mathcal H.
		\end{equation}
		Taking $v=\mathbb PF$ in
		\eqref{eq:local-projection-orthogonality}, we obtain
		\[
		\|\mathbb PF\|_{L^2(\Omega)}^2
		=\int_\Omega F\cdot\mathbb PF dx
		\leq\|F\|_{L^2(\Omega)}\|\mathbb PF\|_{L^2(\Omega)},
		\]
		and consequently
		\begin{equation}\label{eq:local-projection-L2}
			\|\mathbb PF\|_{L^2(\Omega)}\leq\|F\|_{L^2(\Omega)}.
		\end{equation}
		Now let $F\in H^1(\Omega;\mathbb R^2)$.
		For $\eta\in C_c^\infty(\Omega)$, one has
		$\nabla\cdot(\nabla^\perp\eta)
		=-\partial_{x_1}\partial_{x_2}\eta
		+\partial_{x_2}\partial_{x_1}\eta=0$.
		Since $\nabla^\perp\eta$ is compactly supported in $\Omega$,
		integration by parts gives, for every $\psi\in H^1(\Omega)$,
		\[
		\int_\Omega\nabla^\perp\eta\cdot\nabla\psi dx
		=-\int_\Omega\psi\,\nabla\cdot(\nabla^\perp\eta)dx=0.
		\]
		Thus $\nabla^\perp\eta\in\mathcal H$, and
		\eqref{eq:local-projection-orthogonality} yields
		$\int_\Omega(F-\mathbb PF)\cdot\nabla^\perp\eta dx=0$.
		Using the definition of distributional curl, it follows that
		\[
		\begin{aligned}
			\left\langle\operatorname{curl}(\mathbb PF),\eta\right\rangle
			=-\int_\Omega\mathbb PF\cdot\nabla^\perp\eta dx=-\int_\Omega F\cdot\nabla^\perp\eta dx
			=\int_\Omega\operatorname{curl}F \eta dx.
		\end{aligned}
		\]
		Since $\eta$ is arbitrary, we arrive at
		\begin{equation}\label{eq:local-projection-curl}
			\operatorname{curl}(\mathbb PF)=\operatorname{curl}F
			\quad\text{in the sense of distributions in }\Omega.
		\end{equation}
		In particular, $\operatorname{curl}(\mathbb PF)\in L^2(\Omega)$.
		Applying \eqref{eq:exterior-hodge-H1} to $\mathbb PF\in\mathcal H$
		and using \eqref{eq:local-projection-L2}--\eqref{eq:local-projection-curl},
		we conclude that
		\begin{equation}\label{eq:local-helmholtz-H1}
			\begin{aligned}
				\|\mathbb PF\|_{H^1(\Omega)}
				&\leq K\bigl(\|\mathbb PF\|_{L^2(\Omega)}
				+\|\operatorname{curl}(\mathbb PF)\|_{L^2(\Omega)}\bigr)\\
				&\leq K\bigl(\|F\|_{L^2(\Omega)}
				+\|\operatorname{curl}F\|_{L^2(\Omega)}\bigr)
				\leq K\|F\|_{H^1(\Omega)}.
			\end{aligned}
		\end{equation}
		Therefore $\mathbb PF\in\mathcal V$ whenever
		$F\in H^1(\Omega;\mathbb R^2)$.
		
		\noindent\underline{\textbf{(ii) Estimates for $S$.}}
		Let $S$ be the Stokes operator on $\mathcal H$ with
		$v\cdot\mathbf n=0$ and $\operatorname{curl}v=0$ on
		$\partial\Omega$.
		By \cite[Propositions~3.4, 4.9 and~4.10,
		Corollary~4.13]{Kunstmann2025}, $S$ is nonnegative and
		self-adjoint, and
		\begin{equation}\label{eq:local-stokes-domains}
			\begin{gathered}
				D(S)=\left\{v\in H^2(\Omega;\mathbb R^2)\cap\mathcal H:
				\operatorname{curl}v=0\text{ on }\partial\Omega\right\},
				\qquad Sv=-\Delta v,\\
				D(S^{\frac12})=\mathcal V,\qquad
				\|S^{\frac12}v\|_{L^2(\Omega)}^2
				=\|\operatorname{curl}v\|_{L^2(\Omega)}^2\quad\text{for}
				\quad v\in\mathcal V.
			\end{gathered}
		\end{equation}
		Using \eqref{eq:exterior-hodge-H1} and
		$\|\operatorname{curl}v\|_{L^2(\Omega)}\leq K\|v\|_{H^1(\Omega)}$,
		we obtain
		\begin{equation}\label{eq:local-stokes-form-norm}
			K^{-1}\|v\|_{H^1(\Omega)}^2
			\leq\|v\|_{L^2(\Omega)}^2+\|S^{\frac12}v\|_{L^2(\Omega)}^2
			\leq K\|v\|_{H^1(\Omega)}^2,
			\qquad v\in\mathcal V.
		\end{equation}
		For $v\in D(S)$, put $w=\operatorname{curl}v$.
		The boundary condition in \eqref{eq:local-stokes-domains}
		gives $w\in H_0^1(\Omega)$.
		Since $\nabla\cdot v=0$, we arrive at
		\begin{equation}\label{eq:local-stokes-vorticity-identity}
			Sv=-\Delta v=-\nabla^\perp w,
			\qquad
			\|\nabla w\|_{H^j(\Omega)}=\|Sv\|_{H^j(\Omega)},
			\qquad j\in\{0,1\},
		\end{equation}
		where the second identity for $j=1$ requires $Sv\in H^1(\Omega)$.
		Moreover, self-adjointness and \eqref{eq:local-stokes-domains} yield
		\[
		\begin{aligned}
			\|w\|_{L^2(\Omega)}^2
			=\|S^{\frac12}v\|_{L^2(\Omega)}^2
			=\langle Sv,v\rangle_{L^2(\Omega)}\leq\|Sv\|_{L^2(\Omega)}\|v\|_{L^2(\Omega)}
			\leq\frac12\|Sv\|_{L^2(\Omega)}^2
			+\frac12\|v\|_{L^2(\Omega)}^2.
		\end{aligned}
		\]
		Combining this bound with
		\eqref{eq:local-stokes-vorticity-identity}, for $j\in\{0,1\}$, we obtain
		\begin{equation}\label{eq:local-vorticity-regularity}
			\begin{aligned}
				\|w\|_{H^{j+1}(\Omega)}
				&\leq K_j\bigl(\|w\|_{L^2(\Omega)}
				+\|\nabla w\|_{H^j(\Omega)}\bigr)\leq K_j\bigl(\|Sv\|_{H^j(\Omega)}
				+\|v\|_{L^2(\Omega)}\bigr).
			\end{aligned}
		\end{equation}
		For $j=0$, \eqref{eq:exterior-hodge-H2} therefore gives
		\begin{equation}\label{eq:local-stokes-H2}
			\begin{aligned}
				\|v\|_{H^2(\Omega)}
				\leq K\bigl(\|v\|_{L^2(\Omega)}+\|w\|_{H^1(\Omega)}\bigr)\leq K\bigl(\|Sv\|_{L^2(\Omega)}+\|v\|_{L^2(\Omega)}\bigr).
			\end{aligned}
		\end{equation}
		The reverse bound follows from $Sv=-\Delta v$.
		Thus the norm $\|v\|_{L^2(\Omega)}+\|Sv\|_{L^2(\Omega)}$
		on $D(S)$ is equivalent to $\|v\|_{H^2(\Omega)}$.
		To prove the $H^3(\Omega)$ estimate, assume additionally
		$Sv\in H^1(\Omega;\mathbb R^2)$.
		Let $\varphi$ be the stream function in
		\eqref{eq:exterior-stream-function}, so that
		\[
		v=\nabla^\perp\varphi,\qquad
		\Delta\varphi=w\quad\text{in}\quad\Omega,\qquad
		\varphi|_{r=1}=0.
		\]
		The interior and Dirichlet boundary estimates
		\cite[Theorems~8.10 and~8.13]{GT2001}, applied on fixed
		neighborhoods of $\{1\leq r\leq3\}$, give
		\[
		\|\varphi\|_{H^4(\{1<r<3\})}
		\leq K\bigl(\|w\|_{H^2(\{1<r<4\})}
		+\|\varphi\|_{H^1(\{1<r<4\})}\bigr).
		\]
		Using \eqref{eq:exterior-stream-local-bound} and
		\eqref{eq:local-vorticity-regularity} with $j=1$, we obtain
		\begin{equation}\label{eq:local-stokes-near-H3}
			\begin{aligned}
				\|v\|_{H^3(\{1<r<3\})}
				&\leq K\|\varphi\|_{H^4(\{1<r<3\})}\leq K\bigl(\|w\|_{H^2(\Omega)}
				+\|\varphi\|_{H^1(B_4(0))}\bigr)\\
				&\leq K\bigl(\|Sv\|_{H^1(\Omega)}+\|v\|_{L^2(\Omega)}\bigr).
			\end{aligned}
		\end{equation}
		For the remaining region, choose a fixed smooth radial function
		$\chi$, with $\chi=0$ for $r\leq2$ and $\chi=1$ for $r\geq3$.
		Extend $\chi v$ by zero to $\mathbb R^2$.
		Since it vanishes near $\partial\Omega$, the product rule gives
		\[
		\Delta(\chi v)
		=-\chi Sv+2\sum_{i=1}^2
		(\partial_{x_i}\chi)\partial_{x_i}v+(\Delta\chi)v
		\quad\text{in}\quad\mathbb R^2.
		\]
		The derivatives of $\chi$ are supported in $\{2\leq r\leq3\}$.
		By \eqref{eq:local-stokes-H2}, there holds
		\begin{equation}\label{eq:local-stokes-cutoff-H1}
			\begin{aligned}
				\|\Delta(\chi v)\|_{H^1(\mathbb R^2)}
				\leq K\bigl(\|Sv\|_{H^1(\Omega)}+\|v\|_{H^2(\Omega)}\bigr)\leq K\bigl(\|Sv\|_{H^1(\Omega)}+\|v\|_{L^2(\Omega)}\bigr).
			\end{aligned}
		\end{equation}
		Taking the Fourier transform on $\mathbb R^2$ gives
		$\widehat{\Delta(\chi v)}(\xi)=-|\xi|^2\widehat{\chi v}(\xi)$.
		By Plancherel's theorem and
		$(1+|\xi|^2)^3\leq K[1+(1+|\xi|^2)|\xi|^4]$,
		we obtain
		\[
		\begin{aligned}
			\int_{\mathbb R^2}(1+|\xi|^2)^3
			|\widehat{\chi v}(\xi)|^2 d\xi
			&\leq K\|\chi v\|_{L^2(\mathbb R^2)}^2
			+K\|\Delta(\chi v)\|_{H^1(\mathbb R^2)}^2
			<\infty.
		\end{aligned}
		\]
		Thus $\chi v\in H^3(\mathbb R^2;\mathbb R^2)$.
		Using \eqref{eq:local-stokes-cutoff-H1}, we conclude that
		\[
		\begin{aligned}
			\|\chi v\|_{H^3(\mathbb R^2)}
			\leq K\bigl(
			\|\chi v\|_{L^2(\mathbb R^2)}
			+\|\Delta(\chi v)\|_{H^1(\mathbb R^2)}
			\bigr)\leq K\bigl(
			\|Sv\|_{H^1(\Omega)}+\|v\|_{L^2(\Omega)}
			\bigr).
		\end{aligned}
		\]
		Since $1-\chi=0$ for $r\geq3$ and $\chi$ is fixed
		and smooth, the product rule gives
		\[
		\|(1-\chi)v\|_{H^3(\Omega)}
		\leq K\|v\|_{H^3(\{1<r<3\})}.
		\]
		Combining this bound with
		\eqref{eq:local-stokes-near-H3}, the preceding estimate,
		and $v=(1-\chi)v+\chi v$, we obtain
		\[
		\begin{aligned}
			\|v\|_{H^3(\Omega)}
			&\leq
			\|(1-\chi)v\|_{H^3(\Omega)}
			+\|\chi v\|_{H^3(\mathbb R^2)}\leq
			K\|v\|_{H^3(\{1<r<3\})}
			+\|\chi v\|_{H^3(\mathbb R^2)}\\
			&\leq K\bigl(
			\|Sv\|_{H^1(\Omega)}+\|v\|_{L^2(\Omega)}
			\bigr).
		\end{aligned}
		\]
		Together with \eqref{eq:local-stokes-H2}, which gives
		the case $j=0$, this proves
		\begin{equation}\label{eq:local-slip-regularity}
			\|v\|_{H^{j+2}(\Omega)}
			\leq K_j\bigl(
			\|Sv\|_{H^j(\Omega)}+\|v\|_{L^2(\Omega)}
			\bigr),
			\qquad j\in\{0,1\}.
		\end{equation}
		Here $v\in D(S)$, with the additional assumption
		$Sv\in H^1(\Omega;\mathbb R^2)$ when $j=1$.	
		If $v\in D(S^{\frac32})$, the spectral theorem gives
		$v\in D(S)$ and $Sv\in D(S^{\frac12})=\mathcal V$.
		Applying \eqref{eq:local-stokes-form-norm} to $Sv$, followed by
		\eqref{eq:local-slip-regularity} with $j=1$, yields
		\[
		\begin{aligned}
			\|Sv\|_{H^1(\Omega)}
			&\leq K\bigl(\|Sv\|_{L^2(\Omega)}
			+\|S^{\frac12}(Sv)\|_{L^2(\Omega)}\bigr),\\
			\|v\|_{H^3(\Omega)}
			&\leq K\bigl(\|S^{\frac32}v\|_{L^2(\Omega)}
			+\|Sv\|_{L^2(\Omega)}+\|v\|_{L^2(\Omega)}\bigr).
		\end{aligned}
		\]
		
		\noindent\underline{\textbf{(iii) Estimates for $A_D$.}}
		The Dirichlet Laplacian $A_D$ is the nonnegative self-adjoint
		operator associated with the form
		$\int_\Omega\nabla z\cdot\nabla\zeta dx$ for $z,\zeta\in H_0^1(\Omega)$.
		Consequently,
		\[
		\begin{gathered}
			D(A_D)=\{z\in H_0^1(\Omega):-\Delta z\in L^2(\Omega)\},
			\qquad A_Dz=-\Delta z,\\
			D(A_D^{\frac12})=H_0^1(\Omega),\qquad
			\|A_D^{\frac12}z\|_{L^2(\Omega)}^2
			=\int_\Omega|\nabla z|^2 dx\quad\text{for}
			\quad z\in H_0^1(\Omega).
		\end{gathered}
		\]
		For $z\in D(A_D)$, one has
		\[
		(1-\Delta)z=z+A_Dz\in L^2(\Omega),\qquad
		z\in H_0^1(\Omega).
		\]
		Applying \eqref{eq:screened-W2p} with $p=2$ gives
		$z\in H^2(\Omega)$ and
		\[
		\|z\|_{H^2(\Omega)}
		\leq K\|z+A_Dz\|_{L^2(\Omega)}
		\leq K\bigl(\|A_Dz\|_{L^2(\Omega)}+\|z\|_{L^2(\Omega)}\bigr).
		\]
		Conversely, $z\in H^2(\Omega)\cap H_0^1(\Omega)$ implies
		$-\Delta z\in L^2(\Omega)$. Hence
		$D(A_D)=H^2(\Omega)\cap H_0^1(\Omega)$.
		For $f\in H^1(\Omega)$, one has
		\[
		(1-\Delta)((1+A_D)^{-1}f)=f,\qquad
		(1+A_D)^{-1}f\in H_0^1(\Omega).
		\]
		The $H^3(\Omega)$ estimate in
		Lemma~\ref{lem:screened-regularity} therefore gives
		\begin{equation}\label{eq:local-dirichlet-regularity}
			\begin{aligned}
				\|z\|_{H^2(\Omega)}
				&\leq K\bigl(\|A_Dz\|_{L^2(\Omega)}+\|z\|_{L^2(\Omega)}\bigr),
				&&z\in D(A_D),\\
				\|(1+A_D)^{-1}f\|_{H^3(\Omega)}
				&\leq K\|f\|_{H^1(\Omega)},
				&&f\in H^1(\Omega).
			\end{aligned}
		\end{equation}
		
		\noindent{\textbf{Step 2. Estimates for the linear heat and Stokes equations.}}
		Let $0<\tau\leq1$ and consider
		\[
		\partial_tz+A_Dz=f,\qquad z(0)=z_0\in H_0^1(\Omega),
		\qquad f\in L^2(0,\tau;L^2(\Omega)).
		\]
		Testing with $A_Dz$ and $z$, respectively, gives
		\[
		\begin{aligned}
			\frac{d}{dt}\|A_D^{\frac12}z\|_{L^2(\Omega)}^2
			+\|A_Dz\|_{L^2(\Omega)}^2\leq\|f\|_{L^2(\Omega)}^2,\quad
			\frac{d}{dt}\|z\|_{L^2(\Omega)}^2
			\leq\|z\|_{L^2(\Omega)}^2+\|f\|_{L^2(\Omega)}^2.
		\end{aligned}
		\]
		Integration, Gronwall's inequality,
		$\partial_tz=f-A_Dz$, and \eqref{eq:local-dirichlet-regularity}
		yield
		\begin{equation}\label{eq:local-scalar-linear}
			\begin{aligned}
				&\|z\|_{C([0,\tau];H^1(\Omega))}
				+\|z\|_{L^2(0,\tau;H^2(\Omega))}
				+\|\partial_tz\|_{L^2(0,\tau;L^2(\Omega))}\\&\leq K\bigl(\|z_0\|_{H^1(\Omega)}
				+\|f\|_{L^2(0,\tau;L^2(\Omega))}\bigr).
			\end{aligned}
		\end{equation}
		For the velocity, consider
		\[
		\partial_tv+Sv=g,\qquad v(0)=v_0\in D(S),
		\qquad g\in L^2(0,\tau;\mathcal V).
		\]
		Applying $S^{\frac12}$ and testing with $S^{\frac32}v$ gives
		\[
		\begin{aligned}
			\frac12\frac{d}{dt}\|Sv\|_{L^2(\Omega)}^2
			+\|S^{\frac32}v\|_{L^2(\Omega)}^2
			=\langle S^{\frac12}g,S^{\frac32}v\rangle_{L^2(\Omega)}\leq\frac12\|S^{\frac12}g\|_{L^2(\Omega)}^2
			+\frac12\|S^{\frac32}v\|_{L^2(\Omega)}^2.
		\end{aligned}
		\]
		Testing with $v$ also yields
		\[
		\frac{d}{dt}\|v\|_{L^2(\Omega)}^2
		\leq\|v\|_{L^2(\Omega)}^2+\|g\|_{L^2(\Omega)}^2.
		\]
		Using \eqref{eq:local-stokes-domains},
		\eqref{eq:local-stokes-form-norm},
		\eqref{eq:local-slip-regularity}, and $\partial_tv=g-Sv$, we obtain
		\begin{equation}\label{eq:local-velocity-linear}
			\begin{aligned}
				&\|v\|_{C([0,\tau];H^2(\Omega))}
				+\|v\|_{L^2(0,\tau;H^3(\Omega))}
				+\|\partial_tv\|_{L^2(0,\tau;H^1(\Omega))}\\
				&\leq K\bigl(\|v_0\|_{H^2(\Omega)}
				+\|g\|_{L^2(0,\tau;H^1(\Omega))}\bigr).
			\end{aligned}
		\end{equation}
		The constants in \eqref{eq:local-scalar-linear} and
		\eqref{eq:local-velocity-linear} are independent of $\tau\leq1$.
		The solutions are defined by
		\[
		\begin{aligned}
			z(t)=e^{-tA_D}z_0+\int_0^te^{-(t-s)A_D}f(s)ds,\qquad
			v(t)=e^{-tS}v_0+\int_0^te^{-(t-s)S}g(s)ds.
		\end{aligned}
		\]
		To justify the preceding estimates and the continuity in time,
		for $m\geq1$ define
		\[
		\begin{aligned}
			z_m(t)
			&=e^{-(t+m^{-1})A_D}z_0
			+\int_0^t e^{-(t-s+m^{-1})A_D}f(s)ds,\\
			v_m(t)
			&=e^{-(t+m^{-1})S}v_0
			+\int_0^t e^{-(t-s+m^{-1})S}g(s)ds.
		\end{aligned}
		\]
		These functions satisfy
		\[
		\begin{aligned}
			\partial_tz_m+A_Dz_m
			&=e^{-\frac{A_D}{m}}f,
			&\qquad z_m(0)&=e^{-\frac{A_D}{m}}z_0,\\
			\partial_tv_m+Sv_m
			&=e^{-\frac{S}{m}}g,
			&\qquad v_m(0)&=e^{-\frac{S}{m}}v_0.
		\end{aligned}
		\]
		For each fixed $m$, the spectral theorem gives
		\[
		\begin{aligned}
			\|A_De^{-\frac{A_D}{m}}\|_{L^2(\Omega)\to L^2(\Omega)}
			\leq\sup_{\lambda\geq0}\lambda e^{-\frac{\lambda}{m}}
			<\infty,\qquad
			\|S^{\frac32}e^{-\frac{S}{m}}\|_{\mathcal H\to\mathcal H}
			\leq\sup_{\lambda\geq0}\lambda^{\frac32}e^{-\frac{\lambda}{m}}
			<\infty.
		\end{aligned}
		\]
		Consequently, the integral formulas above define
		$z_m\in C([0,\tau];D(A_D))$ and
		$v_m\in C([0,\tau];D(S^{\frac32}))$.
		In particular, the energy calculations leading to
		\eqref{eq:local-scalar-linear} and
		\eqref{eq:local-velocity-linear} apply to $z_m$ and $v_m$.
		The spectral theorem, together with
		\eqref{eq:local-stokes-form-norm},
		\eqref{eq:local-stokes-H2}, and
		\eqref{eq:local-dirichlet-regularity}, also gives
		\[
		\begin{aligned}
			e^{-\frac{A_D}{m}}z_0&\to z_0
			\quad\text{in}\quad H^1(\Omega),\qquad\qquad
			e^{-\frac{S}{m}}v_0\to v_0
			\quad\text{in}\quad H^2(\Omega),\\
			e^{-\frac{A_D}{m}}f&\to f
			~\quad\text{in}\quad L^2(0,\tau;L^2(\Omega)),\quad
			e^{-\frac{S}{m}}g\to g
			~\quad\text{in}\quad L^2(0,\tau;H^1(\Omega))
		\end{aligned}
		\]
		as $m\to\infty$.
		Applying \eqref{eq:local-scalar-linear} to $z_m-z_\ell$,
		we obtain
		\[
		\begin{aligned}
			&\|z_m-z_\ell\|_{C([0,\tau];H^1(\Omega))}
			+\|z_m-z_\ell\|_{L^2(0,\tau;H^2(\Omega))}
			+\|\partial_tz_m-\partial_tz_\ell\|
			_{L^2(0,\tau;L^2(\Omega))}\\
			\leq~& K\Bigl(
			\|(e^{-\frac{A_D}{m}}-e^{-\frac{A_D}{\ell}})z_0\|_{H^1(\Omega)}
			+\|(e^{-\frac{A_D}{m}}-e^{-\frac{A_D}{\ell}})f\|
			_{L^2(0,\tau;L^2(\Omega))}
			\Bigr)
			\to 0
		\end{aligned}
		\]
		as $m,\ell\to\infty$.
		Similarly, \eqref{eq:local-velocity-linear} gives
		\[
		\begin{aligned}
			&\|v_m-v_\ell\|_{C([0,\tau];H^2(\Omega))}
			+\|v_m-v_\ell\|_{L^2(0,\tau;H^3(\Omega))}
			+\|\partial_tv_m-\partial_tv_\ell\|
			_{L^2(0,\tau;H^1(\Omega))}\\
			\leq~& K\Bigl(
			\|(e^{-\frac{S}{m}}-e^{-\frac{S}{\ell}})v_0\|_{H^2(\Omega)}
			+\|(e^{-\frac{S}{m}}-e^{-\frac{S}{\ell}})g\|
			_{L^2(0,\tau;H^1(\Omega))}
			\Bigr)
			\to 0.
		\end{aligned}
		\]
		Thus both sequences converge in all the spaces appearing
		on the left-hand sides. Passing to the limit in their integral
		formulas identifies the limits with $z$ and $v$ defined above.
		In particular,
		\[
		z_m\to z
		\quad\text{in}\quad C([0,\tau];H^1(\Omega)),\qquad
		v_m\to v
		\quad\text{in}\quad C([0,\tau];H^2(\Omega)).
		\]
		Since every approximating function is continuous in the
		corresponding spatial norm, the limits have the asserted
		continuity in time.
		Convergence of the time and spatial derivatives also permits passage to the
		limit in the linear equations and in
		\eqref{eq:local-scalar-linear}--\eqref{eq:local-velocity-linear}.
		Finally, let $(z_1,v_1)$ and $(z_2,v_2)$ be two solutions
		with the same initial data and forcing terms.
		Their differences satisfy
		\[
		\begin{cases}
			\partial_t(z_1-z_2)+A_D(z_1-z_2)=0,\\
			(z_1-z_2)(0)=0,
		\end{cases}
		\qquad
		\begin{cases}
			\partial_t(v_1-v_2)+S(v_1-v_2)=0,\\
			(v_1-v_2)(0)=0.
		\end{cases}
		\]
		Applying \eqref{eq:local-scalar-linear} and
		\eqref{eq:local-velocity-linear} to these differences yields
		\[
		\|z_1-z_2\|_{C([0,\tau];H^1(\Omega))}
		+\|v_1-v_2\|_{C([0,\tau];H^2(\Omega))}=0.
		\]
		This proves uniqueness.
		
		\noindent{\textbf{Step 3. Local existence and uniqueness.}}
		For $c=(1+A_D)^{-1}n$, Lemma~\ref{lem:screened-regularity}
		and Sobolev embedding give
		\begin{equation}\label{eq:local-chemical-H3}
			\|c\|_{H^3(\Omega)}
			+\|\nabla c\|_{L^\infty(\Omega)}
			+\|\Delta c\|_{L^4(\Omega)}
			\leq K\|n\|_{H^1(\Omega)}.
		\end{equation}
		H\"older's inequality and the product rule then imply
		\begin{equation}\label{eq:local-nonlinear-bounds}
			\begin{aligned}
				\|u\cdot\nabla n\|_{L^2(\Omega)}
				&\leq\|u\|_{L^\infty(\Omega)}\|\nabla n\|_{L^2(\Omega)}
				\leq K\|u\|_{H^2(\Omega)}\|n\|_{H^1(\Omega)},\\
				\|\nabla\cdot(n\nabla c)\|_{L^2(\Omega)}
				&\leq\|\nabla n\|_{L^2(\Omega)}\|\nabla c\|_{L^\infty(\Omega)}
				+\|n\|_{L^4(\Omega)}\|\Delta c\|_{L^4(\Omega)}
				\leq K\|n\|_{H^1(\Omega)}^2,\\
				\|u\cdot\nabla u\|_{H^1(\Omega)}
				&\leq K\Bigl(
				\|u\|_{L^\infty(\Omega)}\|\nabla u\|_{L^2(\Omega)}
				+\|\nabla u\|_{L^4(\Omega)}^2
				\\&+\|u\|_{L^\infty(\Omega)}\|D^2u\|_{L^2(\Omega)}\Bigr)\leq K\|u\|_{H^2(\Omega)}^2,\\
				\|ne_r\|_{H^1(\Omega)}&\leq K\|n\|_{H^1(\Omega)}.
			\end{aligned}
		\end{equation}
		The Cartesian formulas for $r^{-1}e_\theta$ and $e_r$ give
		\[
		\|r^{-1}e_\theta\|_{W^{2,\infty}(\Omega)}
		+\|e_r\|_{W^{1,\infty}(\Omega)}\leq K.
		\]
		Together with \eqref{eq:local-helmholtz-H1}, these bounds imply
		\begin{equation}\label{eq:local-background-bounds}
			\begin{aligned}
				\|Ar^{-1}e_\theta\cdot\nabla n\|_{L^2(\Omega)}
				&\leq A\|n\|_{H^1(\Omega)},\\
				\bigl\|\mathbb P\bigl[-Ar^{-1}e_\theta\cdot\nabla u
				-Au\cdot\nabla(r^{-1}e_\theta)-2A(-u_2,u_1)\bigr]\bigr\|_{H^1(\Omega)}
				&\leq KA\|u\|_{H^2(\Omega)}.
			\end{aligned}
		\end{equation}
		For two pairs $(n,u)$ and $(\bar n,\bar u)$, let
		$\bar c=(1+A_D)^{-1}\bar n$. The identities
		\[
		\begin{aligned}
			u\cdot\nabla n-\bar u\cdot\nabla\bar n
			&=(u-\bar u)\cdot\nabla n+\bar u\cdot\nabla(n-\bar n),\\
			n\nabla c-\bar n\nabla\bar c
			&=(n-\bar n)\nabla c+\bar n\nabla(c-\bar c),\\
			u\cdot\nabla u-\bar u\cdot\nabla\bar u
			&=(u-\bar u)\cdot\nabla u+\bar u\cdot\nabla(u-\bar u)
		\end{aligned}
		\]
		and \eqref{eq:local-chemical-H3}, applied to $n-\bar n$, give
		\begin{equation}\label{eq:local-nonlinear-differences}
			\begin{aligned}
				\|u\cdot\nabla n-\bar u\cdot\nabla\bar n\|_{L^2(\Omega)}
				&\leq KR\bigl(\|u-\bar u\|_{H^2(\Omega)}
				+\|n-\bar n\|_{H^1(\Omega)}\bigr),\\
				\|\nabla\cdot(n\nabla c-\bar n\nabla\bar c)\|_{L^2(\Omega)}
				&\leq K\Bigl(\|n-\bar n\|_{H^1(\Omega)}\|c\|_{H^3(\Omega)}
				+\|\bar n\|_{H^1(\Omega)}\|c-\bar c\|_{H^3(\Omega)}\Bigr)\\
				&\leq KR\|n-\bar n\|_{H^1(\Omega)},\\
				\|u\cdot\nabla u-\bar u\cdot\nabla\bar u\|_{H^1(\Omega)}
				&\leq KR\|u-\bar u\|_{H^2(\Omega)}
			\end{aligned}
		\end{equation}
		whenever both pairs have $H^1(\Omega)\times H^2(\Omega)$
		norm at most $R$. Here the product calculation in
		\eqref{eq:local-nonlinear-bounds} also gives
		$\|v\cdot\nabla z\|_{H^1(\Omega)}
		\leq K\|v\|_{H^2(\Omega)}\|z\|_{H^2(\Omega)}$.
		Consider
		\[
		\begin{aligned}
			X_\tau&=C([0,\tau];H_0^1(\Omega))\times C([0,\tau];D(S)),\\
			\|(n,u)\|_{X_\tau}
			&=\|n\|_{C([0,\tau];H^1(\Omega))}
			+\|u\|_{C([0,\tau];H^2(\Omega))}.
		\end{aligned}
		\]
		By \eqref{eq:local-stokes-H2} and $Sv=-\Delta v$, the norm on
		$D(S)$ is equivalent to the $H^2(\Omega)$ norm. Thus $X_\tau$,
		with the displayed norm, is a Banach space.
		Let $\mathcal T(n,u)$ solve the linear problems in Step~2
		with the prescribed initial data and the right-hand sides of
		\eqref{eq:local-rotated-scalar}--\eqref{eq:local-rotated-velocity}
		evaluated at $(n,u)$. By \eqref{eq:local-helmholtz-H1}, the
		velocity right-hand side belongs to $L^2(0,\tau;\mathcal V)$.
		The estimates \eqref{eq:local-scalar-linear},
		\eqref{eq:local-velocity-linear}, and
		\eqref{eq:local-nonlinear-bounds}--\eqref{eq:local-nonlinear-differences}
		give, on the ball of radius $R$ in $X_\tau$,
		\begin{equation}\label{eq:local-fixed-point-bounds}
			\begin{aligned}
				&\|\mathcal T(n,u)\|_{X_\tau}
				\leq K_A\bigl(\|n_{\mathrm{in}}\|_{H^1(\Omega)}
				+\|u_{\mathrm{in}}\|_{H^2(\Omega)}\bigr)
				+\tau^{\frac12}K_A(R+R^2),\\
				&\|\mathcal T(n,u)-\mathcal T(\bar n,\bar u)\|_{X_\tau}
				\leq\tau^{\frac12}K_A(1+2R)
				\|(n-\bar n,u-\bar u)\|_{X_\tau},
			\end{aligned}
		\end{equation}
		where $K_A\geq1$ is independent of $R$ and $\tau\leq1$.
		Choose
		\[
		R=2K_A\bigl(\|n_{\mathrm{in}}\|_{H^1(\Omega)}
		+\|u_{\mathrm{in}}\|_{H^2(\Omega)}\bigr)+1,
		\qquad \tau=\min\{1,[4K_A(1+2R)]^{-2}\}.
		\]
		Then \eqref{eq:local-fixed-point-bounds} yields
		\[
		\|\mathcal T(n,u)\|_{X_\tau}
		\leq\frac{R-1}{2}+\frac R4<R,\qquad
		\|\mathcal T(n,u)-\mathcal T(\bar n,\bar u)\|_{X_\tau}
		\leq\frac14\|(n-\bar n,u-\bar u)\|_{X_\tau}.
		\]
		The contraction mapping theorem gives a fixed point
		$(n,u)\in X_\tau$.
		Evaluated at this fixed point, the right-hand sides of
		\eqref{eq:local-rotated-scalar} and
		\eqref{eq:local-rotated-velocity} belong to
		$L^2(0,\tau;L^2(\Omega))$ and
		$L^2(0,\tau;\mathcal V)$, respectively.
		Applying \eqref{eq:local-scalar-linear} and
		\eqref{eq:local-velocity-linear}, we obtain
		\[
		\begin{aligned}
			n&\in L^2(0,\tau;H^2(\Omega)),
			&\qquad\partial_tn&\in L^2(0,\tau;L^2(\Omega)),\\
			u&\in L^2(0,\tau;H^3(\Omega;\mathbb R^2)),
			&\qquad\partial_tu&\in L^2(0,\tau;H^1(\Omega;\mathbb R^2)).
		\end{aligned}
		\]
		Together with the definition of $X_\tau$ and
		the elliptic estimate \eqref{eq:local-chemical-H3},
		these bounds give the regularity stated in
		Proposition~\ref{prop:local-strong-theory}.
		To prove uniqueness, consider two strong solutions
		with the same initial data.
		On each compact interval of their common existence,
		choose $R$ bounding the
		$H^1(\Omega)\times H^2(\Omega)$ norms of both solutions.
		Choose $0<\delta\leq1$ such that
		$\delta^{\frac12}K_A(1+2R)\leq\frac12$.
		The second estimate in \eqref{eq:local-fixed-point-bounds}
		then implies that the difference of the two solutions
		vanishes on the first interval of length at most $\delta$.
		Their values agree at its right endpoint, so the same
		estimate applies on the next interval.
		Repeating this argument proves uniqueness throughout
		their common interval of existence.
		The choices of the ball radius and $\tau$ in
		\eqref{eq:local-fixed-point-bounds} depend only on
		$A$, $\Omega$, and an upper bound for
		$\|n_{\mathrm{in}}\|_{H^1(\Omega)}
		+\|u_{\mathrm{in}}\|_{H^2(\Omega)}$.
		Since the equations have no explicit dependence on time,
		the change of variable $t\mapsto t-t_0$ leaves
		the estimates unchanged.
		Thus the same $\tau$ applies to initial data
		prescribed at any $t_0\geq0$, with the solution
		defined on $[t_0,t_0+\tau]$.
		
		The component removed by $\mathbb P$ determines the pressure gradient.
		Using \eqref{eq:local-projection-curl} and the displayed curl identities gives
		\[
		\partial_t w+Ar^{-2}\partial_\theta w+u\cdot\nabla w
		=\Delta w-r^{-1}\partial_\theta n,
		\qquad w=\operatorname{curl}u.
		\]
		This is the vorticity equation in \eqref{eq:physical-scalars-moments}.
		
		\noindent{\textbf{Step 4. Nonnegativity and the $L^1(\Omega)$ bound.}}
		Set $n_-=\max\{-n,0\}$. The field $Ar^{-1}e_\theta+u$ is
		bounded, divergence free, and tangent to $\partial\Omega$.
		Testing \eqref{eq:local-rotated-scalar} with $-n_-$ therefore gives
		\[
		\begin{aligned}
			&\frac12\frac{d}{dt}\|n_-\|_{L^2(\Omega)}^2
			+\|\nabla n_-\|_{L^2(\Omega)}^2
			\\=~&\int_\Omega n_-\nabla c\cdot\nabla n_- dx\leq\frac12\|\nabla n_-\|_{L^2(\Omega)}^2
			+\frac12\|\nabla c\|_{L^\infty(\Omega)}^2
			\|n_-\|_{L^2(\Omega)}^2.
		\end{aligned}
		\]
		By \eqref{eq:local-chemical-H3},
		$\int_0^\tau\|\nabla c(t)\|_{L^\infty(\Omega)}^2 dt<\infty$.
		Since $n_-(0)=0$, Gronwall's inequality implies $n\geq0$.
		Assume additionally $n_{\mathrm{in}}\in L^1(\Omega)$.
		For almost every $t$, $n(t)\in H^2(\Omega)\cap H_0^1(\Omega)$,
		and
		\[
		0\leq\frac{n(t,1+h,\cdot)}h
		=\frac1h\int_1^{1+h}\partial_rn(t,r,\cdot)dr
		\to\partial_rn(t,1,\cdot)
		\quad\text{in}\quad L^2(\mathbb T)\quad\text{as}\quad h\to 0^{+}.
		\]
		Hence $\partial_rn(t,1,\theta)\geq0$ almost everywhere.
		The boundary trace estimate also gives
		\[
		\|\partial_rn\|_{L^2(0,\tau;L^2(\partial\Omega))}
		\leq K\|n\|_{L^2(0,\tau;H^2(\Omega))}<\infty.
		\]
		For $R\geq2$, choose radial cutoffs $0\leq\chi_R\leq1$
		that increase monotonically to $1$ as $R\to\infty$, equal to one
		for $r\leq R$ and zero for $r\geq2R$, with
		\[
		\|\nabla\chi_R\|_{L^\infty(\Omega)}\leq\frac{K}{R},\qquad
		\|\Delta\chi_R\|_{L^2(\Omega)}\leq \frac{K}{R},\qquad
		r^{-1}e_\theta\cdot\nabla\chi_R=0.
		\]
		Integrating \eqref{eq:local-rotated-scalar} against $\chi_R$
		and using $n=0$ at $r=1$, we obtain
		\begin{equation}\label{eq:local-mass-cutoff}
			\begin{aligned}
				\frac{d}{dt}\int_\Omega n\chi_R dx
				&=-\int_{\partial\Omega}\partial_rn dS
				+\int_\Omega n\Delta\chi_R dx
				+\int_\Omega n(u+\nabla c)\cdot\nabla\chi_R dx\\
				&\leq\frac K R\|n\|_{L^2(\Omega)}
				\bigl(1+\|u\|_{L^2(\Omega)}+\|\nabla c\|_{L^2(\Omega)}\bigr).
			\end{aligned}
		\end{equation}
		For fixed $R$, $\chi_R\in L^2(\Omega)$ and
		$\partial_tn\in L^2(0,\tau;L^2(\Omega))$ give
		$t\longmapsto\int_\Omega n(t)\chi_R dx\in H^1(0,\tau)$.
		Integrating \eqref{eq:local-mass-cutoff} over $[0,t]$ therefore gives
		\[
		\begin{aligned}
			\int_\Omega n(t)\chi_R dx
			\leq\|n_{\mathrm{in}}\|_{L^1(\Omega)}
			+\frac K R\int_0^t\|n(s)\|_{L^2(\Omega)}
			\bigl(1+\|u(s)\|_{L^2(\Omega)}+\|\nabla c(s)\|_{L^2(\Omega)}\bigr)ds.
		\end{aligned}
		\]
		The time integral is finite by the strong solution bounds and
		\eqref{eq:local-chemical-H3}. Letting $R\to\infty$ and applying
		the monotone convergence theorem proves
		$\|n(t)\|_{L^1(\Omega)}\leq\|n_{\mathrm{in}}\|_{L^1(\Omega)}$
		for every $0\leq t\leq\tau$.
	\end{proof}
	
	\subsection{An $L^2$ blow-up criterion for strong solutions}
	The local solution can be extended whenever the norms
	required in Proposition~\ref{prop:local-strong-theory}
	remain bounded. A uniform $L^2(\Omega)$ bound for the
	density supplies this control and therefore rules out a
	finite maximal existence time.
	
	\begin{proposition}
		\label{prop:density-continuation}
		Fix $A\geq1$. Assume that the initial data satisfy the hypotheses of
		Proposition~\ref{prop:local-strong-theory} and that
		$n_{\mathrm{in}}\in L^\infty(\Omega)$. Let $T_*$ denote the maximal
		existence time of the corresponding strong solution. If
		\begin{equation}\label{eq:continuation-L2-bound}
			\sup_{0\leq t<T_*}\|n(t)\|_{L^2(\Omega)}<\infty,
		\end{equation}
		then $T_*=\infty$.
	\end{proposition}
	\begin{proof}
		Suppose that $T_*<\infty$ and that
		\eqref{eq:continuation-L2-bound} holds.
		Fix $0<\tau<T_*$.
		The constants below may depend on $\Omega$, $A$, $T_*$,
		the initial norms, and the bound in
		\eqref{eq:continuation-L2-bound}, but are independent of $\tau$.
		The initial velocity assumptions imply
		\begin{equation}\label{eq:continuation-initial-vorticity}
			w_{\mathrm{in}}=\operatorname{curl}u_{\mathrm{in}}
			\in H_0^1(\Omega),\qquad
			\|w_{\mathrm{in}}\|_{H^1(\Omega)}
			\leq K\|u_{\mathrm{in}}\|_{H^2(\Omega)}.
		\end{equation}
		
		\noindent{\textbf{Step 1. Bounds for $n$ and $c$ and energy
				estimates for $u$ and $w$.}}
		By Proposition~\ref{prop:local-strong-theory} and Sobolev embedding, one deduces
		\[
		\|Ar^{-1}e_\theta+u\|_{L^\infty(0,\tau;L^\infty(\Omega))}
		\leq A+K\|u\|_{C([0,\tau];H^2(\Omega))}<\infty.
		\]
		The field $Ar^{-1}e_\theta+u$ is divergence free and tangent
		to $\partial\Omega$. Applying Lemma~\ref{lem:exterior-moser}
		and using \eqref{eq:continuation-L2-bound}, we obtain
		\[
		\begin{aligned}
			\sup_{0\leq t\leq\tau}\|n(t)\|_{L^\infty(\Omega)}
			\leq K\left(
			\|n_{\mathrm{in}}\|_{L^\infty(\Omega)}
			+\sup_{0\leq t<T_*}\|n(t)\|_{L^2(\Omega)}
			\right)\left(
			1+\sup_{0\leq t<T_*}\|n(t)\|_{L^2(\Omega)}^2
			\right)\leq K.
		\end{aligned}
		\]
		Using \eqref{eq:screened-W2p} with $p=2,4$,
		the embedding $W^{2,4}(\Omega)\hookrightarrow W^{1,\infty}(\Omega)$,
		and interpolation, we also have
		\[
		\begin{aligned}
			\|c(t)\|_{L^2(\Omega)}
			+\|\nabla c(t)\|_{L^\infty(\Omega)}
			&\leq K\bigl(\|n(t)\|_{L^2(\Omega)}
			+\|n(t)\|_{L^4(\Omega)}\bigr)\\
			&\leq K\bigl(\|n(t)\|_{L^2(\Omega)}
			+\|n(t)\|_{L^2(\Omega)}^{\frac12}
			\|n(t)\|_{L^\infty(\Omega)}^{\frac12}\bigr)\leq K.
		\end{aligned}
		\]
		Since $\tau<T_*$ is arbitrary, these estimates give
		\begin{equation}\label{eq:continuation-density-chemical}
			\sup_{0\leq t<T_*}\bigl(
			\|n(t)\|_{L^\infty(\Omega)}
			+\|c(t)\|_{L^2(\Omega)}
			+\|\nabla c(t)\|_{L^\infty(\Omega)}\bigr)\leq K.
		\end{equation}
		To estimate the velocity, note that
		\begin{equation}\label{eq:continuation-background-coefficient}
			\begin{gathered}
				\nabla\cdot(r^{-1}e_\theta)=0,\qquad
				(r^{-1}e_\theta)\cdot\mathbf n=0
				\quad\text{on}\quad\partial\Omega,\qquad
				|\nabla(r^{-1}e_\theta)|\leq Kr^{-2}\leq K.
			\end{gathered}
		\end{equation}
		Incompressibility and the boundary conditions therefore give
		\[
		\int_\Omega
		\bigl((Ar^{-1}e_\theta+u)\cdot\nabla u\bigr)\cdot u dx
		=\frac12\int_\Omega
		(Ar^{-1}e_\theta+u)\cdot\nabla|u|^2 dx=0.
		\]
		Moreover,
		\[
		\langle Su,u\rangle_{L^2(\Omega)}
		=\|w\|_{L^2(\Omega)}^2,\qquad
		(-u_2,u_1)\cdot u=0,
		\]
		where the first identity follows from
		\eqref{eq:local-stokes-domains}.
		Since $u\in\mathcal H$, the orthogonality relation
		\eqref{eq:local-projection-orthogonality} allows us to remove
		$\mathbb P$ when testing \eqref{eq:local-rotated-velocity} with $u$.
		Using \eqref{eq:continuation-background-coefficient} and
		the Cauchy--Schwarz and Young inequalities, we obtain
		\begin{equation}\label{eq:continuation-fluid-energy}
			\begin{aligned}
				&\frac12\frac{d}{dt}\|u\|_{L^2(\Omega)}^2
				+\|w\|_{L^2(\Omega)}^2=-A\int_\Omega
				\bigl(u\cdot\nabla(r^{-1}e_\theta)\bigr)\cdot u dx
				+\int_\Omega ne_r\cdot u dx\\
				&\leq KA\|u\|_{L^2(\Omega)}^2
				+\|n\|_{L^2(\Omega)}\|u\|_{L^2(\Omega)}\leq K(1+A)\|u\|_{L^2(\Omega)}^2
				+\frac12\|n\|_{L^2(\Omega)}^2.
			\end{aligned}
		\end{equation}
		Gronwall's inequality and
		\eqref{eq:continuation-L2-bound} yield
		\begin{equation}\label{eq:continuation-velocity-L2}
			\begin{aligned}
				\sup_{0\leq t\leq\tau}\|u(t)\|_{L^2(\Omega)}^2
				\leq e^{K(1+A)T_*}\left(
				\|u_{\mathrm{in}}\|_{L^2(\Omega)}^2
				+\int_0^\tau\|n(s)\|_{L^2(\Omega)}^2 ds
				\right)\leq K.
			\end{aligned}
		\end{equation}
		Testing the vorticity equation in
		\eqref{eq:physical-scalars-moments} with $w$ gives
		\begin{equation*}
			\begin{aligned}
				\frac12\frac{d}{dt}\|w\|_{L^2(\Omega)}^2
				+\|\nabla w\|_{L^2(\Omega)}^2
				=\int_\Omega ne_\theta\cdot\nabla w dx\leq\frac12\|\nabla w\|_{L^2(\Omega)}^2
				+\frac12\|n\|_{L^2(\Omega)}^2.
			\end{aligned}
		\end{equation*}
		Integrating it over $[0,t]$
		and using \eqref{eq:continuation-L2-bound}, we obtain
		\[
		\begin{aligned}
			\|w(t)\|_{L^2(\Omega)}^2
			+\int_0^t\|\nabla w(s)\|_{L^2(\Omega)}^2 ds
			\leq\|w_{\mathrm{in}}\|_{L^2(\Omega)}^2
			+\int_0^t\|n(s)\|_{L^2(\Omega)}^2 ds\leq K
		\end{aligned}
		\]
		for $0\leq t\leq\tau$.
		Combining this estimate with
		\eqref{eq:continuation-velocity-L2} and
		\eqref{eq:exterior-hodge-H1} gives
		\begin{equation}\label{eq:exterior-low-fluid}
			\begin{aligned}
				&\sup_{0\leq t\leq\tau}\bigl(
				\|u(t)\|_{H^1(\Omega)}^2
				+\|w(t)\|_{L^2(\Omega)}^2\bigr)\\\leq~& K\sup_{0\leq t\leq\tau}\bigl(
				\|u(t)\|_{L^2(\Omega)}^2
				+\|w(t)\|_{L^2(\Omega)}^2\bigr)\leq K.
			\end{aligned}
		\end{equation}
		
		\noindent{\textbf{Step 2. Estimates for
				$\nabla n$ and $\nabla w$.}}
		For $f=n,w$, Proposition~\ref{prop:local-strong-theory} gives
		\[
		\begin{gathered}
			f\in C([0,\tau];H_0^1(\Omega))
			\cap L^2(0,\tau;H^2(\Omega)),\qquad
			\partial_tf\in L^2(0,\tau;L^2(\Omega)).
		\end{gathered}
		\]
		For $w$, these conclusions follow from
		$w=\operatorname{curl}u$ and
		$\partial_tw=\operatorname{curl}\partial_tu$.
		Testing the density and vorticity equations in
		\eqref{eq:physical-scalars-moments} with $-\Delta n$ and $-\Delta w$
		therefore yields, for almost every $t\in(0,\tau)$,
		\begin{equation}\label{eq:continuation-gradient-identities}
			\begin{aligned}
				&\frac12\frac{d}{dt}\|\nabla n\|_{L^2}^2
				+\|\Delta n\|_{L^2}^2
				=\int_\Omega
				\bigl((Ar^{-1}e_\theta+u)\cdot\nabla n\bigr)\Delta n dx
				+\int_\Omega\nabla\cdot(n\nabla c)\Delta n dx,\\
				&\frac12\frac{d}{dt}\|\nabla w\|_{L^2}^2
				+\|\Delta w\|_{L^2}^2
				=\int_\Omega
				\bigl((Ar^{-1}e_\theta+u)\cdot\nabla w\bigr)\Delta w dx
				+\int_\Omega\nabla\cdot(ne_\theta)\Delta w dx.
			\end{aligned}
		\end{equation}
		For $f=n,w$, the boundary condition $f=0$ implies that
		its tangential derivative vanishes on $\partial\Omega$.
		Since $u\cdot\mathbf n=0$, one has
		$(u\cdot\nabla f)|_{\partial\Omega}=0$.
		Integration by parts and $\nabla\cdot u=0$ give
		\[
		\begin{aligned}
			\int_\Omega(u\cdot\nabla f)\Delta f dx
			&=-\sum_{i,j=1}^2\int_\Omega
			(\partial_{x_i}u_j)(\partial_{x_j}f)(\partial_{x_i}f)dx
			-\frac12\int_\Omega u\cdot\nabla|\nabla f|^2 dx\\
			&=-\sum_{i,j=1}^2\int_\Omega
			(\partial_{x_i}u_j)(\partial_{x_j}f)(\partial_{x_i}f)dx.
		\end{aligned}
		\]
		These identities follow by spatial cutoffs and approximation
		in the stated strong solution class.
		Applying the Gagliardo--Nirenberg inequality to an
		$H^1$ extension of $\nabla f$, and then applying
		\eqref{eq:screened-W2p} with $p=2$ to
		$(1-\Delta)f=f-\Delta f$, we obtain
		\begin{equation}\label{eq:continuation-gradient-interpolation}
			\begin{aligned}
				\|\nabla f\|_{L^4(\Omega)}^2
				&\leq K\|\nabla f\|_{L^2(\Omega)}
				\bigl(\|D^2f\|_{L^2(\Omega)}
				+\|\nabla f\|_{L^2(\Omega)}\bigr)\\
				&\leq K\|\nabla f\|_{L^2(\Omega)}
				\bigl(\|\Delta f\|_{L^2(\Omega)}
				+\|f\|_{L^2(\Omega)}\bigr).
			\end{aligned}
		\end{equation}
		By \eqref{eq:continuation-L2-bound} and
		\eqref{eq:exterior-low-fluid}, both $\|f\|_{L^2(\Omega)}$
		and $\|\nabla u\|_{L^2(\Omega)}$ are bounded independently of $\tau$.
		Thus \eqref{eq:continuation-gradient-interpolation} and
		Young's inequality imply
		\begin{equation}\label{eq:continuation-gradient-transport}
			\begin{aligned}
				\left|\int_\Omega(u\cdot\nabla f)\Delta f dx\right|
				&\leq\|\nabla u\|_{L^2(\Omega)}
				\|\nabla f\|_{L^4(\Omega)}^2\leq K\|\nabla f\|_{L^2(\Omega)}
				\bigl(\|\Delta f\|_{L^2(\Omega)}+1\bigr)\\
				&\leq\frac18\|\Delta f\|_{L^2(\Omega)}^2
				+K\|\nabla f\|_{L^2(\Omega)}^2+K.
			\end{aligned}
		\end{equation}
		The same integration by parts applies to $r^{-1}e_\theta$.
		Using \eqref{eq:continuation-background-coefficient}, we get
		\begin{equation}\label{eq:continuation-angular-transport}
			\begin{aligned}
				&\left|\int_\Omega
				(Ar^{-1}e_\theta\cdot\nabla f)\Delta f dx\right|\\=~&A\left|\sum_{i,j=1}^2\int_\Omega
				\partial_{x_i}(r^{-1}e_\theta)_j
				(\partial_{x_j}f)(\partial_{x_i}f)dx\right|\leq KA\|\nabla f\|_{L^2(\Omega)}^2.
			\end{aligned}
		\end{equation}
		The elliptic equation in \eqref{eq:physical-scalars-moments}
		gives $\Delta c=c-n$.
		By \eqref{eq:continuation-density-chemical} and
		\eqref{eq:continuation-L2-bound}, there holds
		\[
		\begin{aligned}
			\|\nabla\cdot(n\nabla c)\|_{L^2(\Omega)}
			&\leq\|\nabla c\|_{L^\infty(\Omega)}
			\|\nabla n\|_{L^2(\Omega)}
			+\|n\|_{L^\infty(\Omega)}
			\bigl(\|c\|_{L^2(\Omega)}+\|n\|_{L^2(\Omega)}\bigr)\\
			&\leq K\bigl(1+\|\nabla n\|_{L^2(\Omega)}\bigr).
		\end{aligned}
		\]
		Moreover,
		\[
		\nabla\cdot(ne_\theta)=r^{-1}\partial_\theta n,
		\qquad
		\|\nabla\cdot(ne_\theta)\|_{L^2(\Omega)}
		\leq\|\nabla n\|_{L^2(\Omega)}.
		\]
		The Cauchy--Schwarz and Young inequalities therefore give
		\begin{equation}\label{eq:continuation-chemical-source}
			\begin{aligned}
				\left|\int_\Omega\nabla\cdot(n\nabla c)\Delta n dx\right|
				&\leq K\bigl(1+\|\nabla n\|_{L^2(\Omega)}\bigr)
				\|\Delta n\|_{L^2(\Omega)}\\
				&\leq\frac18\|\Delta n\|_{L^2(\Omega)}^2
				+K\|\nabla n\|_{L^2(\Omega)}^2+K,\\
				\left|\int_\Omega\nabla\cdot(ne_\theta)\Delta w dx\right|
				&\leq\|\nabla n\|_{L^2(\Omega)}\|\Delta w\|_{L^2(\Omega)}\leq\frac18\|\Delta w\|_{L^2(\Omega)}^2
				+K\|\nabla n\|_{L^2(\Omega)}^2.
			\end{aligned}
		\end{equation}
		Substituting \eqref{eq:continuation-gradient-transport},
		\eqref{eq:continuation-angular-transport}, and
		\eqref{eq:continuation-chemical-source} into
		\eqref{eq:continuation-gradient-identities}, and absorbing the
		Laplacian terms, yields
		\begin{equation}\label{eq:continuation-gradient-inequalities}
			\begin{aligned}
				\frac{d}{dt}\|\nabla n\|_{L^2(\Omega)}^2
				+\|\Delta n\|_{L^2(\Omega)}^2
				&\leq K\|\nabla n\|_{L^2(\Omega)}^2+K,\\
				\frac{d}{dt}\|\nabla w\|_{L^2(\Omega)}^2
				+\|\Delta w\|_{L^2(\Omega)}^2
				&\leq K\|\nabla w\|_{L^2(\Omega)}^2
				+K\|\nabla n\|_{L^2(\Omega)}^2+K.
			\end{aligned}
		\end{equation}
		Applying Gronwall's inequality to the first inequality in
		\eqref{eq:continuation-gradient-inequalities} and using
		\eqref{eq:continuation-L2-bound}, we obtain
		\begin{equation}\label{eq:continuation-density-H1}
			\begin{aligned}
				\sup_{0\leq t\leq\tau}\|n(t)\|_{H^1(\Omega)}^2
				\leq\sup_{0\leq t<T_*}\|n(t)\|_{L^2(\Omega)}^2
				+e^{KT_*}\bigl(\|\nabla n_{\mathrm{in}}\|_{L^2(\Omega)}^2
				+KT_*\bigr)\leq K.
			\end{aligned}
		\end{equation}
		Applying Gronwall's inequality to the second inequality, and
		using \eqref{eq:continuation-initial-vorticity},
		\eqref{eq:exterior-low-fluid}, and
		\eqref{eq:continuation-density-H1}, gives
		\begin{equation}\label{eq:continuation-vorticity-H1}
			\begin{aligned}
				\sup_{0\leq t\leq\tau}\|w(t)\|_{H^1(\Omega)}^2
				\leq~&\sup_{0\leq t\leq\tau}\|w(t)\|_{L^2(\Omega)}^2+e^{KT_*}\Big(
				\|\nabla w_{\mathrm{in}}\|_{L^2(\Omega)}^2
				\\&+K\int_0^\tau\bigl(1+\|\nabla n(s)\|_{L^2(\Omega)}^2\bigr)ds
				\Big)\leq K.
			\end{aligned}
		\end{equation}
		
		\noindent{\textbf{Step 3. Continuation beyond $T_*$.}}
		By \eqref{eq:exterior-hodge-H2},
		\eqref{eq:exterior-low-fluid}, and
		\eqref{eq:continuation-vorticity-H1}, we get
		\[
		\|u(t)\|_{H^2(\Omega)}
		\leq K\bigl(\|u(t)\|_{L^2(\Omega)}+\|w(t)\|_{H^1(\Omega)}\bigr)
		\leq K
		\]
		for $0\leq t\leq\tau$.
		Together with \eqref{eq:continuation-density-H1}, this gives
		\begin{equation}\label{eq:continuation-restart-bound}
			\sup_{0\leq t<T_*}\bigl(
			\|n(t)\|_{H^1(\Omega)}+\|u(t)\|_{H^2(\Omega)}\bigr)\leq K,
		\end{equation}
		because $\tau<T_*$ is arbitrary and $K$ is independent of $\tau$.
		At every $0\leq t_0<T_*$, the strong solution satisfies
		\[
		n(t_0)\geq0,\qquad n(t_0)\in H_0^1(\Omega),\qquad
		u(t_0)\in D(S).
		\]
		Thus Proposition~\ref{prop:local-strong-theory} and
		\eqref{eq:continuation-restart-bound} give
		$\delta=\delta(A,\Omega,K)>0$, independent of $t_0$, such that
		these data generate a strong solution on $[t_0,t_0+\delta]$.
		Choose $t_0$ so that
		$0<T_*-t_0<\frac{\delta}{2}$.
		By uniqueness in Proposition~\ref{prop:local-strong-theory}, this local
		solution agrees with the original solution on $[t_0,T_*)$.
		Since $t_0+\delta>T_*$, it extends the original solution
		beyond $T_*$, contradicting the definition of $T_*$.
		Therefore $T_*=\infty$.
	\end{proof}
	\section*{Declarations}
	\noindent\textbf{Acknowledgments}
	The author thanks Professor Wendong Wang for his guidance
	and helpful discussions.
	
	\noindent\textbf{Competing interests.}
	The author has no relevant financial or non-financial interests to disclose.
	
	\noindent\textbf{Data availability.}
	Data sharing is not applicable to this article as no datasets were generated or analysed during the current study.
	
	\noindent\textbf{Use of artificial intelligence.}
	During the preparation of this manuscript, the author used ChatGPT (OpenAI) to assist with language editing, organization, and the checking of selected mathematical estimates. The author independently verified the mathematical arguments and calculations and reviewed and revised all AI-assisted material as necessary. The author takes full responsibility for the content and conclusions of the manuscript.

\end{document}